\input amstex
\magnification=1200
\loadmsbm
\loadeufm
\loadeusm
\UseAMSsymbols
\input amssym.def

\font\BIGtitle=cmr10 scaled\magstep3

\font\boldsectionfont=cmb10 scaled\magstep1
\font\section=cmsy10 scaled\magstep1

\def\scr#1{{\fam\eusmfam\relax#1}}

\def\scrA{{\scr A}}
\def\scrB{{\scr B}}
\def\scrC{{\scr C}}
\def\scrD{{\scr D}}
\def\scrE{{\scr E}}
\def\scrF{{\scr F}}

\def\scrI{{\scr I}}

\def\scrJ{{\scr J}}
\def\scrL{{\scr L}}
\def\scrM{{\scr M}}
\def\scrN{{\scr N}}

\def\scrP{{\scr P}}

\def\scrS{{\scr S}}
\def\scrU{{\scr U}}
\def\scrV{{\scr V}}
\def\scrW{{\scr W}}
\def\scrZ{{\scr Z}}
\def\scrR{{\scr R}}
\def\scrT{{\scr T}}
\def\gr#1{{\fam\eufmfam\relax#1}}

\def\grC{{\gr C}}	
	
\def\grE{{\gr E}}	
	
	\def\grg{{\gr g}}
	
\def\grI{{\gr I}}

\def\grL{{\gr L}}

\def\grN{{\gr N}}	\def\grn{{\gr n}}
	
\def\grP{{\gr P}}	 
	
\def\grR{{\gr R}}	
\def\grS{{\gr S}}	
\def\grT{{\gr T}}

\def\db#1{{\fam\msbfam\relax#1}}

\def\dbA{{\db A}} 
\def\dbC{{\db C}} 
 \def\dbF{{\db F}}
\def\dbG{{\db G}}

 \def\dbN{{\db N}}
 
\def\dbQ{{\db Q}} \def\dbR{{\db R}}
\def\dbS{{\db S}}

 \def\dbZ{{\db Z}}

\def\kbar{\bar{k}}

\def\Ker{\text{Ker}}
\def\der{\text{der}}
\def\Sh{\hbox{\rm Sh}}
\def\Sp{\pmb{\text{Sp}}}
\def\GL{\pmb{\text{GL}}}
\def\SL{\pmb{\text{SL}}}
\def\GSp{\pmb{\text{GSp}}}

\def\sc{\text{sc}}
\def\Res{\text{Res}}
\def\ab{\text{ab}}
\def\ad{\text{ad}}

\def\Gal{\text{Gal}}
\def\Hom{\text{Hom}}
\def\End{\text{End}}

\def\Spec{\text{Spec}}

\def\Lie{\text{Lie}}

\def\leaderfill{\leaders\hbox to 1em
     {\hss.\hss}\hfill}
\def\nspace{\lineskip=1pt\baselineskip=12pt\lineskiplimit=0pt}

\def\finishproclaim{\par\rm
     \ifdim\lastskip<\medskipamount\removelastskip
     \penalty55\medskip\fi}
\def\endproof{$\hfill \square$}
\def\proof{\par\noindent {\it Proof:}\enspace}
\def\references#1{\par
  \centerline{\boldsectionfont References}\smallskip
     \parindent=#1pt\nspace}
\def\Ref[#1]{\par\hang\indent\llap{\hbox to\parindent
     {[#1]\hfil\enspace}}\ignorespaces}
\def\Item#1{\par\smallskip\hang\indent\llap{\hbox to\parindent
     {#1\hfill$\,\,$}}\ignorespaces}
\def\ItemItem#1{\par\indent\hangindent2\parindent
     \hbox to \parindent{#1\hfill\enspace}\ignorespaces}

\def\Le{{\mathchoice{\,{\scriptstyle\le}\,}
  {\,{\scriptstyle\le}\,}
  {\,{\scriptscriptstyle\le}\,}{\,{\scriptscriptstyle\le}\,}}}
\def\Ge{{\mathchoice{\,{\scriptstyle\ge}\,}
  {\,{\scriptstyle\ge}\,}
  {\,{\scriptscriptstyle\ge}\,}{\,{\scriptscriptstyle\ge}\,}}}

\def\arrowsim{\,\smash{\mathop{\to}\limits^{\lower1.5pt
  \hbox{$\scriptstyle\sim$}}}\,}

\def\doublemaprights#1#2#3#4{\raise3pt\hbox{$\mathop{\,\,\hbox to
     #1pt{\rightarrowfill}\kern-30pt\lower3.95pt\hbox to
     #2pt{\rightarrowfill}\,\,}\limits_{#3}^{#4}$}}

\def\rightcapdownarrow{\raise9pt\hbox{$\ssize\cap$}\kern-7.75pt
     \Big\downarrow}

\def\rcapmapdown#1{\rightcapdownarrow\kern-1.0pt\vcenter{
     \hbox{$\scriptstyle#1$}}}

\def\rmapdown#1{\Big\downarrow\kern-1.0pt\vcenter{
     \hbox{$\scriptstyle#1$}}}
\def\rightsubsetarrow#1{{\ssize\subset}\kern-4.5pt\lower2.85pt
     \hbox to #1pt{\rightarrowfill}}
\def\longtwoheadedrightarrow#1{\raise2.2pt\hbox to #1pt{\hrulefill}
     \!\!\!\twoheadrightarrow}

\def\Gal{\operatorname{\hbox{Gal}}}
\def\Hom{\operatorname{\hbox{Hom}}}

\def\im{\hbox{Im}}

\NoBlackBoxes
\parindent=25pt
\document
\footline={\hfil}

\null
\noindent
{\BIGtitle Manin problems for Shimura varieties of Hodge type}

\vskip 0.3 in
\noindent
\centerline{Adrian Vasiu, Binghamton University}
\centerline{J. Ramanujan Math. Soc. {\bf 26}, No. 1 (2011), 31--84}

\footline={\hfill}
\vskip 0.3 in
\noindent
{\bf ABSTRACT.} Let $k$ be a perfect field of characteristic $p>0$. We prove the existence of ascending and descending slope filtrations for Shimura $p$-divisible objects over $k$. We use them to classify rationally these objects over $\bar k$. Among geometric applications, we mention two. First we formulate Manin problems for Shimura varieties of Hodge type. We solve them if either $p\Ge 3$ or $p=2$ and two mild conditions hold. Second we formulate integral Manin problems. We solve them for certain Shimura varieties of PEL type. 
\bigskip\smallskip\noindent
{\bf MSC (2000).} 11G10, 14G15, 11G18, 11G35, 14F30, 14L05, 14K22, and 20G25. 
\bigskip\smallskip
\noindent
{\bf KEY WORDS.} $F$-crystals, $F$-isocrystals, Newton polygons, $p$-divisible objects and groups, reductive group schemes, abelian and Shimura varieties, Hodge cycles, and stratifications. 
\vskip0.1in

\footline={\hss\tenrm \folio\hss}
\pageno=1

\bigskip\smallskip
\noindent
{\boldsectionfont 1. Introduction}

\bigskip
Let $p\in\dbN$ be a prime. Let $k$ be a perfect field of characteristic $p$. Let $W(k)$ be the ring of Witt vectors with coefficients in $k$. Let $B(k)$ be the field of fractions of $W(k)$. Let $\sigma:=\sigma_k$ be the Frobenius automorphism of $k$, $W(k)$, or $B(k)$. An {\it $F$-isocrystal} over $k$ is a pair $(M[{1\over p}],\phi)$, where $M$ is a free $W(k)$-module of finite rank and $\phi:M[{1\over p}]\arrowsim M[{1\over p}]$ is a $\sigma$-linear automorphism of $M[{1\over p}]$. If we have $\phi(M)\subseteq M$, then the pair $(M,\phi)$ is called an {\it $F$-crystal} over $k$. We denote also by $\phi$ the $\sigma$-linear automorphism of $\End(M[{1\over p}])$ that takes $x\in\End(M[{1\over p}])$ to $\phi\circ x\circ\phi^{-1}\in\End(M[{1\over p}])$. 

The $F$-isocrystals were introduced by Dieudonn\'e in his work on finite, flat, commutative group schemes of $p$-power order over $k$. In [11, Thms. 1 and 2] and [29, Ch. 2, Sect. 4] the $F$-isocrystals over $\kbar$ are classified: an $F$-isocrystal over $\bar k$ is uniquely determined up to isomorphisms by its {\it Newton polygon}. Manin also brought the topic into the context of {\it abelian varieties} as follows (see [29, Ch. 4]). Let $r\in\dbN$ and let $(D,\lambda)$ be a principally quasi-polarized {\it $p$-divisible group} over $k$ of height $2r$. Let $(M_0,\phi_0,\psi_0)$ be its principally quasi-polarized {\it Dieudonn\'e module} (see [1, Ch. 4]). The pair $(M_0,\phi_0)$ is an $F$-crystal over $k$, $\psi_0$ is a perfect alternating form on $M_0$, and we have $pM_0\subseteq\phi_0(M_0)\subseteq M_0$ and $\psi_0(\phi_0(x),\phi_0(y))=p\sigma(\psi_0(x,y))$, where $x$, $y\in M_0$.  The Newton polygon of $(M_0,\phi_0)$ has the following three properties: 

\medskip\noindent
{\bf ($\sharp$)} its slopes belong to the interval $[0,1]$, its starting and ending points are $(0,0)$ and $(2r,r)$, and the multiplicity of a slope $\gamma$ is the same as the multiplicity of the slope $1-\gamma$.

\medskip
The first two properties of ($\sharp$) are particular cases of a theorem of Mazur (see [19, Thm. 1.4.1]). The third property of ($\sharp$) is a consequence of the existence of $\lambda$ and thus of $\psi_0$; in the geometric context of abelian varieties, it appears for the first time in [29]. The original Manin problem conjectured that each Newton polygon that satisfies ($\sharp$), is the Newton polygon of an abelian variety over $\bar k$ of dimension $r$. It was first solved in [41] (see also [34] and [35] for two more recent proofs). In all that follows, expressions of the form $\GL_M$ and $\GSp(M_0,\psi_0)$ are viewed as {\it reductive group schemes} over $W(k)$. Thus $\GL_M(W(k))$ is the group of $W(k)$-linear automorphisms of $M$, etc. 
The principally quasi-polarized Dieudonn\'e module of any other principally quasi-polarized $p$-divisible group over $k$ of height $2r$ is isomorphic to $(M_0,g\phi_0,\psi_0)$, where $g\in\Sp(M_0,\psi_0)(W(k))$ can be arbitrary.

Let $T_0$ be a split, maximal torus of $\GSp(M_0,\psi_0)$ whose fibre over $k$ normalizes the kernel of the reduction mod $p$ of $\phi_0$. Let $N_0$ be the normalizer of $T_0\cap \Sp(M_0,\psi_0)$ in $\Sp(M_0,\psi_0)$. It is easy to see that there exists $g_0\in \Sp(M_0,\psi_0)(W(k))$ such that $(g_0\phi_0)(\Lie(T_0)))=\Lie(T_0)$. The three properties of ($\sharp$) are equivalent to:

\medskip\smallskip\noindent
{\bf 1.1. Fact.} {\it For each $g\in \Sp(M_0,\psi_0)(W(k))$, there exists $w_0\in N_0(W(k))$ such that the Newton polygon of $(M_0,g\phi_0)$ is the same as the Newton polygon of $(M_0,w_0g_0\phi_0)$.} 

\medskip
This equivalence can be checked easily by considering the actions of $w_0g_0\phi_0$'s on the rank $1$ direct summands of $M_0$ normalized by $T_0$. Since many years (see the paragraph before Subsection 1.3 below) it was expected that Fact 1.1 also holds if $\GSp(M_0,\psi_0)$ is replaced by an arbitrary reductive, closed subgroup scheme of $\GL_{M_0}$ related in a natural way to $\phi_0$. The interest in the resulting problems stems from the study of {\it Shimura varieties}. A general study of such problems was started by Kottwitz in [23]. See [8], [32], [31], and [45, Subsect. 2.5] for different types of Shimura varieties. For instance, the Shimura varieties of PEL type are {\it moduli spaces} of polarized abelian varieties endowed with endomorphisms (see [40] and [7]). Also the Shimura varieties of abelian (resp. of Hodge) type are moduli spaces of polarized abelian motives (resp. of polarized abelian varieties) endowed with {\it Hodge cycles} (see [7] and [31]). The Shimura varieties of abelian type are the main testing ground for many parts of the Langlands program (like zeta functions, local correspondences, etc.). The deep understanding of their zeta functions depends on the Langlands--Rapoport conjecture (see [27], [28], and [30]) on $\overline{\dbF_p}$-valued points of special fibres of their good {\it integral models} in mixed characteristic $(0,p)$. To solve this conjecture and to aim in extending [38] and [18] to all Shimura varieties, one needs a good theory of isomorphism classes of $F$-isocrystals with additional structures that are {\it crystalline realizations} of abelian motives associated naturally to these $\overline{\dbF_p}$-valued points. 

This paper and [46] to [49] are part of a sequence meant to contribute to such a theory. The below notion Shimura $p$-divisible object axiomatizes all crystalline realizations one can (or hopes to) associate to points with values in perfect fields of ``good" integral models in mixed characteristic $(0,p)$ of arbitrary (quotients of) Shimura varieties. The main goals of the paper are to prove the existence of slope filtrations of Shimura $p$-divisible objects over $k$, to classify rationally such $p$-divisible objects over $\bar k$, and to generalize the original Manin problem to contexts related to Shimura varieties of Hodge type.

\medskip\smallskip\noindent
{\bf 1.2. The language.} In what follows we will use an integral language that:

\medskip
$\bullet$ is closer in spirit to the works [26], [14], and [54] which precede [23] and [37];

\smallskip
$\bullet$ matches naturally with the crystalline languages used in the references [2], [12], [13], [19], [20], [36], and [45] to [50] whose results are often playing key roles in this paper;

\smallskip
$\bullet$ is a version of the rational language used in [38, Subsects. 1.7, 1.8, and 1.18];

\smallskip
$\bullet$ works over $W(k)$ with $k$ an arbitrary perfect field without assuming the existence of suitable $\dbZ_p$-structures;

\smallskip
$\bullet$ is the only language that can be naturally adapted over smooth algebras over $W(k)$ whose $p$-adic completions are equipped with Frobenius lifts. 

\medskip
For reader's convenience, very often we will make as well connections to the works [23], [37], and [38].
For $a$, $b\in\dbZ$, $b\Ge a$, we define $S(a,b):=\{a,a+1,\ldots,b\}$.  

\medskip\noindent
{\bf 1.2.1. Definition.} A $p$-{\it divisible object with a reductive group} over $k$ is a triple $(M,\phi,G),$
where $M$ is a free $W(k)$-module of finite rank, the pair $(M[{1\over p}],\phi)$ is an $F$-isocrystal over $k$, and $G$ is a reductive, closed subgroup scheme of $\GL_M$, such that there exists a direct sum decomposition 
$$M=\oplus_{i\in S(a,b)} \tilde F^i(M)\leqno (1)$$ 
for which the following two axioms hold:

\medskip
{\bf (a)} we have $\phi^{-1}(M)=\oplus_{i=a}^b p^{-i}\tilde F^i(M)$ and $\phi(\Lie(G_{B(k)}))=\Lie(G_{B(k)})$;

\smallskip
{\bf (b)} the cocharacter $\mu$ of $\GL_M$ such that $\beta\in \dbG_m(W(k))$ acts through $\mu$ on $\tilde F^i(M)$ as the multiplication by $\beta^{-i}$, factors through $G$.

\medskip
Following [36, Sect. 2] we refer to $\mu:\dbG_m\to G$ as a {\it Hodge cocharacter} of $(M,\phi,G)$ and to (1) as its {\it Hodge decomposition}. If $G=\GL_M$, then often we do not mention $G$ and so implicitly ``with a reductive group". For $i\in S(a,b)$ let $F^i(M):=\oplus_{j=i}^b \tilde F^j(M)$ and let $\phi_i:F^i(M)\to M$ be the restriction of $p^{-i}\phi$ to $F^i(M)$. Triples of the form $(M,(F^i(M))_{i\in S(a,b)},\phi)$  show up in [26], [14], [54], [46], etc., and are objects of a $\dbZ_p$-linear category $p-\scrM\scrF_{[a,b]}(W(k))$ (the morphisms being $W(k)$-linear maps that respect after inverting $p$ the Frobenius endomorphisms and the filtrations). If $n\in\dbN$, then the reduction mod $p^n$ of $(M,(F^i(M))_{i\in S(a,b)},(\phi_i)_{i\in S(a,b)})$ is an object of the abelian category $\scrM\scrF_{[a,b]}(W(k))$ used in [26], [14], and [12]. This and the fact that such a reduction is a natural generalization of a truncated Barsotti--Tate group of level $n$ over $W(k)$ (e.g., see [12, Thm. 7.1] for $p\ge 3$) justifies our terminology ``$p$-divisible object". Let $b_L\in S(0,b-a)$ be the smallest number with the property that we have a direct sum decomposition
$$\Lie(G)=\oplus_{i=-b_L}^{b_L} \tilde F^i(\Lie(G))\leqno (2)$$
such that $\beta\in \dbG_m(W(k))$ acts through $\mu$ on $\tilde F^i(\Lie(G))$ as the multiplication by $\beta^{-i}$. 

We have $b_L=0$ if and only if $\mu$ factors through the center of $G$. If $b_L=1$, then $\mu$ is called a {\it minuscule cocharacter} of $G$. If $b_L\Le 1$, then we say $(M,\phi,G)$ is a {\it Shimura $p$-divisible object} over $k$. If $(a,b)=(0,1)$, then we say $(M,\phi,G)$ is a {\it Shimura $F$-crystal} over $k$; they were extensively used in [45, Sect. 5], in [48], and in many previous works on Shimura varieties of PEL type (see [51], [28], [24], etc.).

For $g\in G(W(k))$ let 
$$\scrC_{g}:=(M,g\phi,G).$$ 
We have $\phi^{-1}(M)=(g\phi)^{-1}(M)$ and thus $\scrC_g$ is also a $p$-divisible object with a reductive group over $k$ that has $\mu$ as a Hodge cocharacter. By the extension of $\scrC_{g}$ to a perfect field $k_1$ that contains $k$, we mean the triple $\scrC_{g}\otimes_k k_1:=(M\otimes_{W(k)} W(k_1),g\phi\otimes\sigma_{k_1},G_{W(k_1)})$.

\medskip\noindent
{\bf 1.2.2. Definition.} Let $g$, $g_1\in G(W(k))$. By an {\it inner} (resp. by a {\it rational inner) isomorphism} between $\scrC_g$ and $\scrC_{g_1}$ we mean an element $h\in G(W(k))$ that normalizes $\phi^{-1}(M)$ and (resp. an element $h\in G(B(k))$) such that we have $hg\phi=g_1{\phi}h$.

\medskip\noindent
{\bf 1.2.3. Remark.} Let $\sigma_0:=\phi\mu(p)$; it is a $\sigma$-linear automorphism of $M$ (cf. axiom 1.2.1 (a)). Thus we get $g\phi=g\sigma_0\mu({1\over p})$ and this differs slightly from the works [23], [37], and [38] (for instance, in [38, Subsects. 1.7 and 1.8] one considers expressions of the form $b\mu(p)\sigma_0$ with $b\in G(B(k))$). In other words, $\mu$ plays here the same role played by $\sigma_0(\mu):=\sigma_0 \mu \sigma_0^{-1}$ in [23], [37], and [38]. Therefore $\scrC_g$ and $\scrC_{g_1}$ are rational isomorphic if and only if $g\sigma_0(\mu)([{1\over p}])$ and $g_1\sigma_0(\mu)([{1\over p}])$ belong to the same $\sigma_0$-conjugacy class of elements of $G(B(k))$ (i.e., if and only if there exists $h\in G(B(k))$ such that we have an identity $h g\sigma_0(\mu)([{1\over p}]) \sigma_0(h)^{-1}=g_1\sigma_0(\mu)([{1\over p}])$, where $\sigma_0(h):=\sigma_0 h\sigma_0^{-1}$). Extra direct connections to the  ``$B(G_{B(k)})$-language" of [23] and [37] are made in Subsections 2.6 and 4.5.

\medskip\noindent
{\bf 1.2.4. Notations.} For $g\in G(W(k))$ let $\grS(g)$ be the set of Newton polygon slopes of $(M[{1\over p}],g\phi)$. Let $M[{1\over p}]=\oplus_{\gamma\in \grS(g)} M_{\gamma}(g)$ be the direct sum decomposition stable under $\phi$ and such that all Newton polygon slopes of $(M_{\gamma}(g),\phi)$ are $\gamma$. If $m\in\dbN$ is such that $m\gamma\in\dbZ$ and if $k=\bar k$, then there exists a $B(k)$-basis for $M_{\gamma}(g)$ formed by elements fixed by $p^{-m\gamma}\phi^m$. Let $\gamma_{1,g}<\gamma_{2,g}<\cdots <\gamma_{n_g,g}$ be the numbers in $\grS(g)$ listed increasingly. For $\gamma\in\grS(g)$ let 
$$W^{\gamma}(M,g\phi):=M\cap (\oplus_{\alpha\in \grS(g)\cap [\gamma,\infty)} M_{\alpha}(g))\;\;\text{and}\;\;W_{\gamma}(M,g\phi):=M\cap (\oplus_{\alpha\in \grS(g)\cap (-\infty,\gamma]} M_{\alpha}(g)).$$ 
Due to the axioms 1.2.1 (a) and (b), the {\it Newton quasi-cocharacter} of $\GL_{M[{1\over p}]}$ defined by $g\phi$ (its definition is reviewed in Subsubsection 2.2.1) factors through $G_{B(k)}$ (see Claim 2.2.2). Let $\nu_{g}$ be this factorization. 

\medskip\smallskip\noindent
{\bf 1.3. The basic results.} In Sections 3 and 4 we mainly study the case when $b_L\Le 1$ and we deal with two aspects of the classification of $\scrC_g$'s up to inner isomorphisms. The two aspects are: slope filtrations and $\nu_{g}$'s. We list the basic results.

\medskip\noindent
{\bf 1.3.1. Theorem.} {\it Let $g\in G(W(k))$. Then there exists a unique parabolic subgroup scheme $P_G^+(g\phi)$ of $G$ such that all Newton polygon slopes of the $F$-isocrystal $(\Lie(P_G^+(g\phi)_{B(k)}),g\phi)$ (resp. $(\Lie(G_{B(k)})/\Lie(P_G^+(g\phi)_{B(k)}),g\phi)$) are non-negative (resp. are negative). If $b_L\le 1$, then there exists a Hodge cocharacter of $(M,g\phi,G)$ that factors through $P_G^+(g\phi)$.}

\medskip\noindent
{\bf 1.3.2. Corollary.} {\it Let $g\in G(W(k))$. We assume $b_L\Le 1$. We have:

\medskip
{\bf (a)} There exists a Hodge cocharacter of $(M,g\phi,G)$ that normalizes $W^{\gamma}(M,g\phi)$ for all $\gamma\in \grS(g)$.

\smallskip
{\bf (b)} Up to a rational inner isomorphism, we can assume that we have a direct sum decomposition $M:=\oplus_{\gamma\in \grS(g)} M\cap M_{\gamma}(g)$ and that $\mu$ normalizes $M\cap M_{\gamma}(g)$ for all $\gamma\in \grS(g)$.}
\medskip

If $b_L\le 1$, then from Corollary 1.3.2 (a) and axiom 1.2.1 (a) we get that the filtration
$$(W^{\gamma_{n_g,g}}(M,g\phi),g\phi)\subseteq (W^{\gamma_{n_g-1,g}}(M,g\phi),g\phi)\subseteq \cdots \subseteq (W^{\gamma_{1,g}}(M,g\phi),g\phi)=(M,g\phi)\leqno (3)$$
is a filtration in the category $p-\scrM(W(k))$ of $p$-divisible objects over $k$ (the morphisms being $W(k)$-linear maps that respect after inverting $p$ the Frobenius endomorphisms) which can be extended to a filtration in the category $p-\scrM\scrF(W(k))$. We refer to it as the {\it descending slope filtration} of $(M,g\phi)$. Replacing ``negative" by ``positive", we get another parabolic subgroup scheme $P_G^-(g\phi)$ of $G$ and the {\it ascending slope filtration}
$$(W_{\gamma_{1,g}}(M,g\phi),g\phi)\subseteq (W_{\gamma_{2,g}}(M,g\phi),g\phi)\subseteq \cdots \subseteq (W_{\gamma_{n_g,g}}(M,g\phi),g\phi)=(M,g\phi)\leqno (4)$$ 
of $(M,g\phi)$ in the category $p-\scrM(W(k))$. We can define Shimura $p$-divisible objects over any field $l$ of characteristic $p$ and we can always speak about their ascending slope filtrations over a Cohen ring $K(l)$ of $l$. Thus we view (4) and its analogue over $K(l)$ as a natural extension of Grothendieck's slope filtrations of $p$-divisible groups over $l$ (see [52]).

We shift our attention to the rational classification of $\scrC_{g}$'s. For the remaining part of Section 1, we will assume that $G$ is split. Let $T$ be a split, maximal torus of $G$ such that $\mu$ factors through it. Until Section 2 we also assume that $\phi(\Lie(T))=\Lie(T)$ (in Subsection 2.5 we check that we can always achieve this by replacing $\phi$ with $g_0\phi$ for some $g_0\in G(W(k))$).

\medskip\noindent
{\bf 1.3.3. Theorem.} {\it Let $g\in G(W(k))$. We assume that $b_L\Le 1$, that $G$ is split, and that $\phi(\Lie(T))=\Lie(T)$. Let $N$ be the normalizer of $T$ in $G$. We have:

\medskip
{\bf (a)} There exists $w\in N(W(k))$ such that $\nu_w$ and $\nu_g$ are $G(B(k))$-conjugate. 

\smallskip
{\bf (b)} If $k=\bar k$, then there exist elements $w\in N(W(k))$ and $h\in G(B(k))$ such that $hw\phi=g{\phi}h$ (thus $\scrC_w$ and $\scrC_g$ are rational inner isomorphic).}

\medskip
Section 2 gathers standard properties needed in Sections 3 and 4. See Subsections 2.3 and 4.1 for the proofs of Theorem 1.3.1 and Corollary 1.3.2. The proof of Theorem 1.3.3 is in two steps (see Subsection 4.2). The first step works for all $b_L\in\dbN\cup\{0\}$ and shows two things. First, using the classification of adjoint group schemes over $\dbZ_p$, it shows that there exists $w\in N(W(k))$ such that all Newton polygon slopes of $(\Lie(G_{B(k)}),w\phi)$ are 0 (see Subsubsection 4.2.2). Second, if $k=\bar k$ and if all Newton polygon slopes of $(\Lie(G_{B(k)}),g\phi)$ and of $(\Lie(G_{B(k)}),g_1\phi)$ are 0, then standard arguments developed in [23] and [37] show that there exist rational inner isomorphisms  between $\scrC_{g}$ and $\scrC_{g_1}$ (see Subsections 2.6 and 2.7). The second step is an inductive one. It works only if $b_L\le 1$ (see Subsubsection 4.2.3). The idea is: if $(\Lie(G_{B(k)}),g\phi)$ has non-zero Newton polygon slopes, then $\scrC_{g}$ has standard forms that reduce the situation to a context in which there exists a Levi subgroup scheme $L$ of $P_G^+(g\phi)$ such that the triple $(M,g\phi,L)$ is a Shimura $p$-divisible object. The standard forms (see Section 3) are also the essence of Theorem 1.3.1 and Corollary 1.3.2. They are rooted on the simple property 2.4 (c) and on the fact that each intersection of two parabolic subgroups of $G_k$ contains a maximal split torus of $G_k$. See Subsections 4.3 to 4.7 for examples and complements to Corollary 1.3.2 and Theorem 1.3.3.

The notion inner isomorphism is a natural extension of the classification ideas of [29]. Parabolic subgroup schemes as $P_G^+(g\phi)$ were first used in [44] for $(a,b)=(0,1)$. To our knowledge, Theorem 1.3.1 is a new result. The Newton polygon translation of Theorem 1.3.3 (a) was indirectly hinted at by the Langlands--Rapoport conjecture. Chai extrapolated this conjecture and stated rather explicitly that Theorem 1.3.3 (a) ought to hold (see [6]). Theorem 1.3.1 and versions of Corollary 1.3.2 and Theorem 1.3.3 were first part of our manuscripts math.NT/0104152 and math.NT/0209410 (only few particular cases of Theorem 1.3.3 (b) were known before math.NT/0104152 and most of them could be deduced from [29]). The paper [25] was written after the mentioned two manuscripts; one can use [25, Thm. 4.3] to recover Theorem 1.3.3 (b). See Corollary 4.4 for an interpretation of Theorem 1.3.3 (b) in terms of equivalence classes. 

The last part of Theorem 1.3.1 and Theorem 1.3.3 (and thus also (3) and (4)) do not hold in general if $b_L\Ge 2$ (see Example 2.3.4).

\medskip\smallskip\noindent
{\bf 1.4. Geometric applications.}  Subsection 5.1 introduces standard Hodge situations. They give rise to good moduli spaces in mixed characteristic $(0,p)$ of principally polarized abelian varieties endowed with (specializations of) Hodge cycles, that generalize the moduli spaces of principally polarized abelian varieties endowed with endomorphisms considered in [51], [24], and [28]. Subsection 5.2 formulates Manin problems for standard Hodge situations. The Main Theorem 5.2.3 solves them if either $p\Ge 3$ or $p=2$ and two mild conditions hold. The fact that the two mild conditions always hold if $p\Ge 3$ is implied by either a particular case of [49, Thm. 1.2 and Lemma 2.5.2 (a)] or [20, Cor. (1.4.3)]. The proof of the Main Theorem 5.2.3 relies on Theorem 1.3.3 (b), on Fontaine comparison theory, and on properties of Shimura varieties and reductive group schemes. In particular, we get a new solution to the original Manin problem mentioned before Fact 1.1 (cf. Subsubsection 5.2.4 (a)). Subsection 5.3 redefines in a more direct way the {\it rational stratifications} of the special fibres of the mentioned good moduli spaces one gets based on [37, Thm. 3.6]. 

In Subsection 5.4 we formulate {\it integral Manin problems} for standard Hodge situations. Theorem 5.4.2 solves them in many cases that pertain to Shimura varieties of PEL type; the simplest example implies that each principally quasi-polarized $p$-divisible group over $\bar k$ of height $2r$ is the one of a principally polarized abelian variety over $\bar k$ of dimension $r$ (cf. Example 5.4.3 (a)). These integral problems are natural extrapolations of the ``combination" between the Manin problems and a motivic conjecture of Milne (see [45, Conj. 5.6.6] and see [49, Thm. 1.2] and [20, Cor. (1.4.3)] for refinements and proofs of it). 
                                                                      
\bigskip\smallskip
\noindent                                             
{\boldsectionfont 2. Preliminaries}

\bigskip
See Subsection 2.1 for our conventions and notations. In Subsections 2.2 and 2.3 we include complements on Newton polygons. We deal with two aspects: Newton (quasi-) cocharacters and parabolic subgroup schemes that correspond to either non-negative or to non-positive Newton polygon slopes. Our approach to Newton quasi-cocharacters is slightly different from the standard one that uses the pro-torus of character group $\dbQ$ (see [23], [37], and [36]). In Subsections 2.4 to 2.7 we list different properties of $\scrC_g$'s. 

\medskip\smallskip\noindent
{\bf 2.1. Notations and conventions.} 
Reductive group schemes have connected fibres. If $\Spec(R)$ is an affine scheme and if $H$ is a reductive group scheme over $R$, let $H^{\der}$, $Z(H)$, $H^{\ab}$, and $H^{\ad}$ be the derived group scheme, the center, the maximal abelian quotient, and the adjoint group scheme of $H$ (respectively). We have $H/Z(H)=H^{\ad}$ and $H/H^{\der}=H^{\ab}$. Let $Z^0(H)$ be the maximal torus of $Z(H)$. Let $\Lie(F)$ be the Lie algebra over $R$ of a smooth, closed subgroup scheme $F$ of $H$. If $R=W(k)$ (like $R=\dbZ_p$), then $H(W(k))$ is called a hyperspecial subgroup of $H(B(k))$ (see [43]). Let $R_0\to R$ be a homomorphism. If it is finite and flat, let $\Res_{R/R_0}H$ be the affine group scheme over $\Spec(R_0)$ that is the Weil restriction of scalars of $H$ (see [4, Subsect. 7.6]). In general, the pullback of an $R_0$-scheme $X$ or $X_{R_0}$ (resp. $X_\backprime$ with $\backprime$ an index) to $R$ is denoted by $X_R$ (resp. $X_{\backprime R}$). Let $\bar E$ be an algebraic closure of a field $E$. 

For an $R$-module $N$ let $N^*:=\Hom_R(N,R)$. Let $N^{\otimes s}\otimes_R N^{*\otimes t}$, with $s,t\in\dbN\cup\{0\}$, be the tensor product of $s$-copies of $N$ with $t$-copies of $N^*$ taken in this order. Let
$$
\scrT(N):=\oplus_{s,t\in\dbN\cup\{0\}} N^{\otimes s}\otimes_R N^{*\otimes t}.
$$ 
A family of tensors of $\scrT(N)$ is denoted in the form $(v_{\alpha})_{\alpha\in\scrJ}$, with $\scrJ$ as the set of indices. We emphasize that we use the same notation for two tensors or bilinear forms obtained one from another by an extension of scalars. Let $N_1$ be another $R$-module. Each isomorphism $f:N\arrowsim N_1$ extends naturally to an isomorphism $\scrT(N)\arrowsim\scrT(N_1)$ and therefore we speak about $f$ taking $v_{\alpha}$ to some specific element of $\scrT(N_1)$. A bilinear form on $N$ is called perfect if it induces an isomorphism $N\arrowsim N^*$. If $N$ is a projective, finitely generated $R$-module, then we view $\GL_N$ as a reductive group scheme over $R$. If $f_1$ and $f_2$ are two $\dbZ$-endomorphisms of $N$, then $f_1f_2:=f_1\circ f_2$. 

Until Section 5, whenever we consider a $p$-divisible object with a reductive group $(M,\phi,G)$ over $k$, the following notations $a$, $b$, $S(a,b)$, $\mu$, $b_L$, $\scrC_{g}$'s, $\tilde F^i(M)$ and $F^i(M)$ with $i\in S(a,b)$, and $\tilde F^i(\Lie(G))$ with $i\in S(-b_L,b_L)$ will be as in Subsubsection 1.2.1. Often we do not mention ``over $k$". Let $P$ be the parabolic subgroup scheme of $G$ that normalizes $F^i(M)$ for all $i\in S(a,b)$. From Subsection 2.3 onward, $\gamma\in\grS(g)$, $M_{\gamma}(g)$, $W^{\gamma}(M,g\phi)$, $W_{\gamma}(M,g\phi)$, and $\nu_g$, will be as in Subsections 1.2.4 and 2.2.3. If all Newton polygon slopes of $(\Lie(G_{B(k)}),\phi)$ are $0$, then we say $(M,\phi,G)$ is {\it basic}. Corollary 2.3.2 shows that $(M,\phi,G)$ is basic if and only if the Newton quasi-cocharacter of $\GL_{M[{1\over p}]}$ defined by $\phi$ factors through $Z^0(G_{B(k)})$ (i.e., the notion basic is compatible with the one introduced first in [23]). We denote also by $\phi$ the $\sigma$-linear automorphism of $M^*[{1\over p}]$ that takes $x\in M^*[{1\over p}]$ to $\sigma x\phi^{-1}\in M^*[{1\over p}]$. Thus $\phi$ acts on $\scrT(M[{1\over p}])$ in the natural tensor way. The identification $\End(M)=M\otimes_{W(k)} M^*$ is compatible with the two $\phi$ actions (defined here and before Fact 1.1). If $(a,b)=(0,1)$, then we refer to $(M,F^1(M),\phi,G)$ as a Shimura filtered $F$-crystal.

\medskip\smallskip\noindent
{\bf 2.2. Quasi-cocharacters.} Let $H$ be a reductive group scheme over a connected scheme $S$. Let $\chi(H)$ be the set of cocharacters of $H$. Let $\Lambda(H):=\chi(H)\times \dbQ\setminus\{0\}$. Let $R(H)$ be the smallest equivalence relation on $\Lambda(H)$ that has the following property: two pairs $(\mu_1,r_1)$, $(\mu_2,r_2)\in \Lambda(H)$ are in relation $R(H)$ if there exists $n\in\dbN$ such that $nr_1$, $nr_2\in\dbZ$, $g.c.d.(nr_1,nr_2)=1$, and $\mu_1^{nr_1}=\mu_2^{nr_2}$. Let $\Xi(H)$ be the set of equivalence classes of $R(H)$. An element of $\Xi(H)$ is called a {\it quasi-cocharacter} of $H$. Identifying $\chi(H)$ with the subset $\chi(H)\times \{1\}$ of $\Lambda(H)$, we view naturally $\chi(H)$ as a subset of $\Xi(H)$. If $H$ is a split torus, then $\chi(H)=X_*(H)$ has a natural structure of a free $\dbZ$-module and we can identify naturally $\Xi(H)=X_*(H)\otimes_{\dbZ}\dbQ$. The group $H(S)$ acts on $\Xi(H)$ via its inner conjugation action on $\chi(H)$. If $f:H\to H_1$ is a homomorphism of reductive group schemes over $S$, then $f_*(\mu_1,r_1):=(f\circ\mu_1,r_1)\in \Lambda(H_1)$. The resulting map $f_*:\Lambda(H)\to \Lambda(H_1)$ is compatible with the $R(H)$ and $R(H_1)$ relations; the quotient map $f_*:\Xi(H)\to \Xi(H_1)$ is compatible with the $H(S)$- and $H_1(S)$-actions. We say $\top\in \Xi(H_1)$ factors through $H$ if $\top\in f_*(\Xi(H))$.

\medskip\noindent
{\bf 2.2.1. The slope context.} Let $\scrC_{1_M}=(M,\phi,G)$ be a $p$-divisible object with a reductive group over $k$. Let $\grS(1_M)$ be as in Subsubsection 1.2.4. For $\gamma\in \grS(1_M)$ we write $\gamma={a_{\gamma}\over {b_{\gamma}}}$, where $a_{\gamma}\in\dbZ$, $b_{\gamma}\in\dbN$, and $g.c.d.(a_{\gamma},b_{\gamma})=1$. Let $d:=l.c.m.(b_{\gamma}|\gamma\in \grS(1_M))$. Each $F$-isocrystal over $\bar k$ is a direct sum of simple $F$-isocrystals which have only one Newton polygon slope, cf. Dieudonn\'e--Manin's classification of $F$-isocrystals over $\bar k$. Let
$$(M\otimes_{W(k)} B(\bar k),\phi\otimes\sigma_{\bar k})=\oplus_{\gamma\in \grS(1_M)} (\bar M_{\gamma}(1_M),\phi\otimes\sigma_{\bar k})\leqno (5)$$ 
be the direct sum decomposition of $F$-isocrystals over $\bar k$ such that $(\bar M_{\gamma}(1_M),\phi\otimes\sigma_{\bar k})$ has only one Newton polygon slope $\gamma$. Let $\bar e\in \End(M)\otimes_{W(k)} B(\bar k)$ be the semisimple element that acts on $\bar M_{\gamma}(1_M)$ as the multiplication by ${da_{\gamma}}\over {b_{\gamma}}$. The Galois group $\Gal(\bar k/k)$ acts on $B(\bar k)$ having $B(k)$ as its fixed field and acts on $\scrT(M\otimes_{W(k)} B(\bar k))$ having $\scrT(M[{1\over p}])$ as its set of fixed elements. For each automorphism $\tau\in\text{Aut}(B(\bar k)/B(k))$ defined by an element of $\Gal(\bar k/k)$, $\phi\otimes\sigma_{\bar k}$ and $1_{M\otimes_{W(k)} B(\bar k)}\otimes\tau$ commute. Thus $\bar e$ is fixed by $\Gal(\bar k/k)$ i.e., $\bar e\in\End(M[{1\over p}])=\End(M)\otimes_{W(k)} B(k)$. Therefore (5) is the tensorization with $B(\bar k)$ of a direct sum decomposition $(M[{1\over p}],\phi)=\oplus_{\gamma\in \grS(1_M)} (M_{\gamma}(1_M),\phi)$. 

Let $\tilde\nu_{1_M}\in \chi(\GL_{M[{1\over p}]})$ be such that $\tilde\nu_{1_M}(p)$ acts on $M_{\gamma}(1_M)$ as the multiplication by $p^{{da_{\gamma}}\over {b_{\gamma}}}$. Let $\nu_{1_M}:=[(\tilde\nu_{1_M},{1\over d})]\in \Xi(\GL_{M[{1\over p}]})$. By abuse of language, we say that $\nu_{1_M}(p)$ acts on $M_{\gamma}(1_M)$ as the multiplication by $p^{\gamma}$. As the decompositions (5) are compatible with morphisms and tensor products of $F$-isocrystals, $\tilde\nu_{1_M}$ and $\nu_{1_M}$ factor through the subgroup of $\GL_{M[{1\over p}]}$ that fixes all tensors of $\scrT(M[{1\over p}])$ fixed by $\phi$.

\medskip\noindent
{\bf 2.2.2. Claim.} {\it Both $\tilde\nu_{1_M}$ and $\nu_{1_M}$ factor through $G_{B(k)}$.}

\medskip
\proof  Let $\sigma_0:=\phi\circ\mu(p)$. We have $\sigma_0(M)=\phi(\oplus_{i=1}^b p^{-1}\tilde F^i(M))=\phi(\phi^{-1}(M))=M$, cf. axiom 1.2.1 (a). Thus $\sigma_0$ is a $\sigma$-linear automorphism of $M$. As $\sigma_0$ normalizes $\Lie(G_{B(k)})$ (cf. axioms 1.2.1 (a) and (b)), it also normalizes $\Lie(G)=\Lie(G_{B(k)})\cap\End(M)$. Let
$$M_{\dbZ_p}:=\{x\in M|\sigma_0(x)=x\}$$
and $\grg_{\dbZ_p}:=\{x\in\Lie(G)|\sigma_0(x)=x\}$. To prove the Claim we can assume $k=\bar k$. As $k=\bar k$, we have $M=M_{\dbZ_p}\otimes_{\dbZ_p} W(k)$ and $\Lie(G)=\grg_{\dbZ_p}\otimes_{\dbZ_p} W(k)\subseteq \End(M_{\dbZ_p})\otimes_{\dbZ_p} W(k)=\End(M)$. 

In this paragraph we follow [49, proof of Lemma 2.5.3] to check that there exists a unique reductive subgroup $G_{\dbQ_p}$ of $\GL_{M_{\dbZ_p}[{1\over p}]}$ whose Lie algebra is $\grg_{\dbZ_p}[{1\over p}]$.  The uniqueness part is implied by [3, Ch. II, Subsect. 7.1]. To check the existence part, we consider commutative $\dbQ_p$-algebras $A$ such that  there exists a reductive, closed subgroup scheme $G_A$ of $\GL_{M_{\dbZ_p}\otimes_{\dbZ_p} A}$ whose Lie algebra is $\grg_{\dbZ_p}\otimes_{\dbZ_p} A$. For instance, $A$ can be $B(k)$ itself and thus we can assume $A$ is a finitely generated $\dbQ_p$-subalgebra of $B(k)$ (cf. [10, Vol. III, Exp. XIX, Rm. 2.9]). By replacing $A$ with $A/J$, where $J$ is a maximal ideal of $A$, we can assume $A$ is a finite field extension of $\dbQ_p$. Even more, we can assume that $A$ is a finite Galois extension of $\dbQ_p$. As $\Lie(G_A)=\grg_{\dbZ_p}\otimes_{\dbZ_p} A$, from [3, Ch. II, Subsection 7.1] we get that the natural action of the Galois group $\Gal(A/\dbQ_p)$ on $\grg_{\dbZ_p}\otimes_{\dbZ_p} A$ is defined naturally by an action of $\Gal(A/\dbQ_p)$ on the subgroup $G_A$ of $\GL_{M_{\dbZ_p}\otimes_{\dbZ_p} A}$. This last action is free. As $G_A$ is an affine scheme, the quotient $G_{\dbQ_p}$ of $G_A$ by $\Gal(A/\dbQ_p)$ exists (cf. [4, Ch. 6, 6.1, Thm. 5]) and it is a reductive subgroup of $\GL_{M_{\dbZ_p}[{1\over p}]}$ whose Lie algebra is $\grg_{\dbZ_p}[{1\over p}]$.

Our notations match i.e., $G_{B(k)}$ is the pullback of $G_{\dbQ_p}$ to $B(k)$ (cf. [3, Ch. II, Subsect. 7.1]). Let $G_{\dbZ_p}$ be the Zariski closure of $G_{\dbQ_p}$ in $\GL_{M_{\dbZ_p}}$; its pullback to $W(k)$ is $G$. Thus $G_{\dbZ_p}$ is a reductive, closed subgroup scheme of $\GL_{M_{\dbZ_p}}$; its Lie algebra is $\grg_{\dbZ_p}$. We note that the existence of $G_{\dbZ_p}$ is also a particular case of [48, Prop. 3.2].

Let $(t_{\alpha})_{\alpha\in\scrJ}$ be a family of tensors of $\scrT(M_{\dbZ_p})$ such that $G_{\dbQ_p}$ is the subgroup of $\GL_{M_{\dbZ_p}[{1\over p}]}$ that fixes $t_{\alpha}$ for all $\alpha\in\scrJ$, cf. [9, Prop. 3.1 (c)]. Each $t_{\alpha}$ is fixed by both $\sigma_0$ and $\mu(p)$. Thus we have $\phi(t_{\alpha})=t_{\alpha}$ for all $\alpha\in\scrJ$. Therefore $\tilde\nu_{1_M}$ fixes each $t_{\alpha}$, cf. end of Subsubsection 2.2.1. Thus both $\tilde\nu_{1_M}$ and $\nu_{1_M}$ factor through $G_{B(k)}$.\endproof

\medskip\noindent
{\bf 2.2.3. Definition.} By the Newton cocharacter (resp. quasi-cocharacter) of $\scrC_{1_M}$ we mean the factorization of $\tilde\nu_{1_M}$ (resp. of $\nu_{1_M}$) through $G_{B(k)}$. Similarly, for $g\in G(W(k))$ let $\nu_g\in \Xi(G_{B(k)})$ be the  Newton quasi-cocharacter of $\scrC_g$. 

\medskip\smallskip\noindent
{\bf 2.3. Sign parabolic subgroup schemes.} Let $\grL\grS(1_M)$ be the set of Newton polygon slopes of $(\Lie(G_{B(k)}),\phi)$. The composite of $\nu_{1_M}$ with the homomorphisms $G_{B(k)}\to \GL_{M[{1\over p}]}\to \GL_{\End(M)[{1\over p}]}$, is the Newton quasi-cocharacter of $(\End(M)[{1\over p}],\phi)$. Thus the Newton polygon slope decomposition 
$$(\Lie(G_{B(k)}),\phi)=\oplus_{\gamma\in \grL\grS(1_M)} (L_{\gamma},\phi_{\gamma})$$is such that $\nu_{1_M}(p)$ acts via inner conjugation on $L_{\gamma}$ as the multiplication by $p^{\gamma}$. As $\im(\tilde\nu_{1_M})$ is a split torus of $G_{B(k)}$, the centralizer $C_G(\phi)$ of $\tilde\nu_{1_M}$ in $G_{B(k)}$ is a reductive group of the same rank as $G_{B(k)}$ (see [10, Vol. III, Exp. XIX, Subsect. 2.8]). We have $\Lie(C_G(\phi))=L_0$. Reductive groups as  $C_G(\phi)$ were first considered in [23, Sect. 6]. 

\medskip\noindent
{\bf 2.3.1. Lemma.} {\it There exists a unique parabolic subgroup scheme $P_G^+(\phi)$ of $G$ such that we have $\Lie(P_G^+(\phi)_{B(k)})=\oplus_{\gamma\in \grL\grS(1_M)\cap [0,\infty)} L_{\gamma}$. 
The group $C_G(\phi)$ is a Levi subgroup of $P_G^+(\phi)_{B(k)}$ and $L_{>0}:=\oplus_{\gamma\in \grL\grS(1_M)\cap (0,\infty)} L_{\gamma}$ is the nilpotent radical of $\Lie(P_G^+(\phi)_{B(k)})$.}

\medskip
\proof By replacing $k$ with a finite Galois extension of $k$, we can assume $G$ is split. Thus the group $C_G(\phi)$ is also split. Let $\Lie(G_{B(k)})=\Lie(T_G(\phi))\bigoplus_{\alpha\in\Phi}\grg_{\alpha}$ be the root decomposition with respect to a split, maximal torus $T_G(\phi)$ of $C_G(\phi)$. Here $\Phi$ is a root system of characters of $T_G(\phi)$ whose irreducible factors are indexed by the simple factors of $G^{\ad}_{B(k)}$.
Let $\gamma_1$, $\gamma_2\in \grL\grS(1_M)$. If $\gamma_1+\gamma_2\notin \grL\grS(1_M)$, let $L_{\gamma_1+\gamma_2}:=0$. As $\phi$ is a $\sigma$-linear Lie automorphism of $\Lie(G_{B(k)})$, we have $[L_{\gamma_1},L_{\gamma_2}]\subseteq L_{\gamma_1+\gamma_2}$. Thus the set $\Phi_+:=\{\alpha\in\Phi|\grg_{\alpha}\subseteq L_0\oplus L_{>0}\}$ is a closed subset of $\Phi$. Let $\Phi_0$ be a simple factor of $\Phi$. The intersection $\Phi_0\cap \Phi_+$ is a closed subset of $\Phi_0$. As $\tilde\nu_{1_M}$ factors through $T_G(\phi)$, each $L_{\gamma}$ is a direct sum of some $\grg_{\alpha}$'s and we have $\grg_{\alpha}\subseteq L_{\gamma}$ if and only if $\grg_{-\alpha}\subseteq L_{-\gamma}$. Thus the intersection $\Phi_0\cap\Phi_+$ is a parabolic subset of $\Phi_0$ in the sense of [5, Ch. VI, 7, Def. 4] i.e., it is closed and we have $\Phi_0=-(\Phi_0\cap\Phi_+)\cup\Phi_0\cap\Phi_+$. There exists a unique parabolic subgroup $P_G^+(\phi)_{B(k)}$ of $G_{B(k)}$ whose Lie algebra is $L_0\oplus L_{>0}$, cf. [10, Vol. III,  Exp. XXVI, Prop. 1.4]. Obviously $L_{>0}$ is a nilpotent ideal of $L_0\oplus L_{>0}$ and thus also of the nilpotent radical $\grn$ of $L_0\oplus L_{>0}$. As $C_G(\phi)$ is reductive and $\Lie(C_G(\phi))=L_0$, we have $\grn\cap L_0=0$. Thus $\grn=L_{>0}$. 

Therefore $\Lie(P_G^+(\phi)_{B(k)})=L_0\oplus \grn=\Lie(C_G(\phi))\oplus\grn$ and thus $C_G(\phi)$ is a Levi subgroup of $P_G^+(\phi)_{B(k)}$. As the $W(k)$-scheme that parametrizes  parabolic subgroup schemes of $G$ is projective (cf. [10, Vol. III, Exp. XXVI, Cor. 3.5]), the Zariski closure $P_G^+(\phi)$ of $P_G^+(\phi)_{B(k)}$ in $G$ is a parabolic subgroup scheme of $G$.\endproof

\medskip
A similar argument shows that there exists a unique parabolic subgroup scheme $P_G^-(\phi)$ of $G$ such that 
we have $\Lie(P_G^-(\phi)_{B(k)})=\oplus_{\gamma\in \grL\grS(1_M)\cap (-\infty,0]} L_{\gamma}.$ The group $C_G(\phi)$ is also a Levi subgroup of $P^-_G(\phi)_{B(k)}$. 

Replacing the role of $\tilde\nu_{1_M}$ with the one of $\mu$, we similarly get that the parabolic subgroup scheme $P$ of Subsection 2.1 exists and is uniquely determined by the equality $\Lie(P)=\oplus_{i=0}^{b_L} \tilde F^i(\Lie(G))$. The Lie algebra of the unipotent radical of $P$ is $\oplus_{i=1}^{b_L} \tilde F^i(\Lie(G))$.

The next Corollary is only a variant of [23, Subsect. 5.2].

\medskip\noindent
{\bf 2.3.2. Corollary.} {\it The following three statements are equivalent:

\medskip
{\bf (a)} $P_G^+(\phi)=G$ (or $P_G^-(\phi)=G$);

\smallskip
{\bf (b)} $\scrC_{1_M}$ is basic;

\smallskip
{\bf (c)} $\tilde\nu_{1_M}$ (or $\nu_{1_M}$) factors through $Z^0(G_{B(k)})$.}

\medskip
\proof
Each statement is equivalent to the statement that $L_0=\Lie(C_G(\phi))$ is $\Lie(G_{B(k)})$.\endproof  

\medskip\noindent
{\bf 2.3.3. Definition.} We call $P_G^+(\phi)$ (resp. $P_G^-(\phi)$) the {\it non-negative} (resp. {\it non-positive}) {\it parabolic subgroup scheme} of $\scrC_{1_M}$. We call $C_G(\phi)$ the {\it Levi subgroup} of $\scrC_{1_M}$.

\medskip\noindent
{\bf 2.3.4. Example.} Suppose that $k=\bar k$, $(a,b)=(0,2)$, $\text{rk}_{W(k)}(M)=4$, $G=\GL_M$, $\text{rk}_{W(k)}(\tilde F^2(M))=\text{rk}_{W(k)}(\tilde F^0(M))=2$, and there exists a $W(k)$-basis $\{e_1,\ldots,e_4\}$ for $M$ formed by elements of $\tilde F^0(M)\cup\tilde F^2(M)$ and such that $\phi(e_i)=p^{n_i}e_i$, where $n_1=n_2=0$ and $n_3=n_4=2$. Let $T$ be the maximal torus of $G$ that normalizes $W(k)e_i$ for all $i\in S(1,4)$. Let $g\in G(W(k))$ be such that $\grS(g)=\{{1\over 2},{3\over 2}\}$, cf. [16, Thm. 2]. For each element $w\in G(W(k))$ that normalizes $T$, the $F$-crystal $(M,w\phi)$ over $k$ has at least one Newton polygon slope which is an integer (to check this one can assume that $w$ permutes the set $\{e_1,\ldots,e_4\}$). Thus $\nu_g$ and $\nu_w$ are not $G(B(k))$-conjugate i.e., Theorem 1.3.3 does not hold in this case. The parabolic subgroup scheme of $G$ that normalizes the direct summand $M\cap M_{3\over 2}(g)$ of $M$ of rank $2$, is $P_G^+(g\phi)$. If there exists a Hodge cocharacter of $(M,g\phi,G)$ that factors through $P_G^+(g\phi)$, then $(M\cap M_{3\over 2}(g),g\phi)$ is a $p$-divisible object whose Hodge slopes  belong to the set $\{0,2\}$ and this, based on Mazur's theorem (see [19, Thm. 1.4.1]), contradicts the fact that the end point of the Newton polygon of $(M\cap M_{3\over 2}(g),g\phi)$ is $(2,3)$. Thus also the last part of Theorem 1.3.1 does not hold in this case.

\medskip\smallskip\noindent
{\bf 2.4. Lemma.} {\it {\bf (a)} An element $h\in G(W(k))$ normalizes $\phi^{-1}(M)$ if and only if we have
$$h(\tilde F^i(M))\subseteq\oplus_{j=a}^b p^{\max\{0,i-j\}}\tilde F^j(M)\;\;\;\;\forall i\in S(a,b).\leqno (6)$$
In particular, each element $h\in P(W(k))$ normalizes $\phi^{-1}(M)$. 

\smallskip
{\bf (b)} For $h\in G(B(k))$ let $h_1:=\phi h\phi^{-1}\in \GL_M(B(k))$. We have $h_1\in G(B(k))$. Also, $h$ normalizes $\phi^{-1}(M)$ if and only if $h_1\in G(W(k))$. 
\smallskip
{\bf (c)} We have $b_L\Le 1$ if and only if each element $h\in\Ker(G(W(k))\to G(k))$ normalizes $\phi^{-1}(M)$.} 

\medskip
\proof
We prove (a). As $h(M)=M$, the inclusions of (6) are equivalent to the inclusion $h(\phi^{-1}(M))\subseteq \phi^{-1}(M)$. As $h(M)=M$ and as $p^b\phi^{-1}(M)\subseteq M$, by reasons of lengths of artinian $W(k)$-modules we get that $h(\phi^{-1}(M))\subseteq \phi^{-1}(M)$ if and only if $h(\phi^{-1}(M))=\phi^{-1}(M)$. Thus (6) holds if and only if $h$ normalizes $\phi^{-1}(M)$. If $h\in P(W(k))$, then $h(\tilde F^i(M))\subseteq h(F^i(M))=F^i(M)=\oplus_{j=i}^b \tilde F^j(M)$; thus (6) holds and therefore $h$ normalizes $\phi^{-1}(M)$. Thus (a) holds. To prove (b) we can assume $k=\bar k$. Let $\sigma_0=\phi\circ\mu(p):M\arrowsim M$ be as in the proof of Claim 2.2.2. We have $h_1=\sigma_0\mu({1\over p})h\mu(p)\sigma_0^{-1}\in \sigma_0G(B(k))\sigma_0^{-1}=G(B(k))$. Thus $h_1\in G(W(k))$ if and only if $h_1(M)=M$. As $\phi^{-1}(h_1(M))=h(\phi^{-1}(M))$, we have $h_1(M)=M$ if and only if $h(\phi^{-1}(M))=\phi^{-1}(M)$. Thus (b) holds. 

To prove  (c), we first check that there exists a unique connected, smooth, commutative, closed subgroup scheme $U$ of $G$ that has $\tilde F^{-b_L}(\Lie(G))$ as its Lie algebra. The uniqueness part follows from [3, Ch. II, Subsect. 7.1]. Due to this, a standard Galois descent shows that to check the existence of $U$ we can consider pullbacks to $W(k_1)$ with $k_1$ a finite Galois extension of $k$ and thus we can assume that $G$ is split. By considering a split maximal torus $T$ of $G$ through which $\mu$ factors and a root decomposition $\Lie(G)=\Lie(T)\bigoplus_{\alpha\in\Phi} \grg_{\alpha}$, as in the proof of Lemma 2.3.1 we argue that there exists a subset $\Phi_0$ of $\Phi$ such that $\tilde F^{-b_L}(\Lie(G))=\bigoplus_{\alpha\in\Phi_0} \grg_{\alpha}$ and $(\Phi_0+\Phi_0)\cap\Phi=\emptyset=\Phi_0\cap (-\Phi_0)$; thus the existence of $U$ is a particular case of [10, Vol. III, Exp. XXII, Prop. 5.6.1 and Cor. 5.6.5] and in fact $U$ is isomorphic to a product of copies of $\dbG_a$. 

We first assume that $b_L\le 1$ and we check that each $h\in\Ker(G(W(k))\to G(k))$ normalizes $\phi^{-1}(M)$. Based on  (b), we only have to show that $h_1\in G(W(k))$. As $b_L\Le 1$, $U$ is in fact the unipotent radical of the parabolic subgroup scheme of $G$ whose Lie algebra is $\oplus_{i=-b_L}^{0} \tilde F^i(\Lie(G))$). As $\Lie(P)=\oplus_{i=0}^{b_L} \tilde F^i(\Lie(G))$ and as $b_L\le 1$, we have a direct sum decomposition $\Lie(G)=\Lie(U)\oplus\Lie(P)$ of $W(k)$-modules (cf. (2)). Thus the product morphism $U\times_{W(k)} P\to G$ is \'etale around the point $(1_M,1_M)\in U(W(k))\times P(W(k))$. Therefore we can write $h=u_hg_h$, where $u_h\in U(W(k))$ and $g_h\in P(W(k))$ are both congruent to $1_M$ mod $p$. Thus $h_1=h_2h_3$, where $h_2:=\phi u_h\phi^{-1}$ and $h_3:=\phi g_h\phi^{-1}$. From  (a) and (b) we get that $h_3\in G(W(k))$. To check that $h_2\in G(W(k))$ we can assume $k=\bar k$. We have $h_2=\sigma_0\mu({1\over p})u_h\mu(p)\sigma_0^{-1}$. As $u_h\in\Ker(U(W(k))\to U(k))$ and $b_L\Le 1$, we have $\mu({1\over p})u_h\mu(p)\in U(W(k))$. [Argument: using an identification $U=\dbG_a^s$ with $s:=\text{rk}_{W(k)} \tilde F^{-b_L}(\Lie(G))\in\dbN\cup\{0\}$, $u_h$ gets identified with an element $v_h$ of $(pW(k))^s\leqslant W(k)^s=\dbG_a^s(W(k))$ and $\mu({1\over p})u_h\mu(p)$ gets identified with $p^{-b_L}v_h\in W(k)^s=\dbG_a^s(W(k))$]. Thus $h_2\in\sigma_0U(W(k))\sigma_0^{-1}\leqslant G(W(k))$. Therefore $h_1=h_2h_3\in G(W(k))$. 

We now assume that all elements $h\in\Ker(G(W(k))\to G(k))$ normalize $\phi^{-1}(M)$ and we check that $b_L\le 1$. We show that the assumption $b_L\Ge 2$ leads to a contradiction. Both Lie algebras $\tilde F^{-b_L}(\Lie(G))$ and $\tilde F^{b_L}(\Lie(G))$ are non-trivial; thus $U$ is non-trivial.  

If $u_{b_L}\in\Ker(U(W(k))\to U(k))$ is not congruent mod $p^{b_L}$ to $1_M$, then we have $\mu({1\over p})u_{b_L}\mu(p)\in U(B(k))\setminus U(W(k))$. [Argument: this is similar to the above argument involving $u_h$ and $v_h$]. If $h$ is $u_{b_L}$, then $h_1=\phi h\phi^{-1}\notin G(W(k))$ and thus $h$ does not normalize $\phi^{-1}(M)$ (cf.  (b)). Contradiction. Therefore $b_L\Le 1$. Thus (c) holds.\endproof  

\medskip\noindent
{\bf 2.4.1. Fact.} {\it If $h\in G(W(k))$ normalizes $\phi^{-1}(M)$, then $\mu_h:=h\mu h^{-1}:\dbG_m\to G$ is a Hodge cocharacter of $\scrC_{1_M}$.}

\medskip
\proof 
The Hodge decomposition of $M$ produced by $\mu_h$ is $M=\oplus_{i=a}^b h(\tilde F^i(M))$. We have $\oplus_{i=a}^b p^{-i}h(\tilde F^i(M))=h(\oplus_{i=a}^b p^{-i}\tilde F^i(M))=h(\phi^{-1}(M))=\phi^{-1}(M)$, cf. axiom 1.2.1 (a) and our hypothesis. Thus by very definitions, $\mu_h$ is a Hodge cocharacter of $\scrC_{1_M}$. \endproof

\medskip\smallskip\noindent
{\bf 2.5. The quasi-split case.} Let $G$ be quasi-split. We recall from the proof of Claim 2.2.2 that $\sigma_0=\phi\mu(p)$ normalizes $\Lie(G)$. Let $B$ be a Borel subgroup scheme of $G$ that is contained in $P$ and that has a maximal torus $T$ through which $\mu$ factors. We have $\phi(\Lie(T))=\sigma_0(\Lie(T))$ and $\phi(\Lie(B))\subseteq \sigma_0(\Lie(B))$. Thus $\phi(\Lie(T))$ is the Lie algebra of the maximal torus $T_0:=\sigma_0(T)$ of the Borel subgroup scheme $B_0:=\sigma_0(B_0)$ of $G$. Let $\bar g_0\in G(k)$ be such that $\bar g_0(T_{0k})\bar g^{-1}_0=T_k$ and $\bar g_0(B_{0k})\bar g_0^{-1}=B_k$, cf. [3, Ch. IV, Thm. 15.14]. As $W(k)$ is $p$-adically complete, there exists $g_0\in G(W(k))$ that lifts $\bar g_0$ and such that $g_0(T_0)g_0^{-1}=T$ (cf. [10, Vol. II, Exp. IX, Thm. 3.6]). The Borel subgroup schemes $g_0B_0g_0^{-1}$ and $B$ of $G$ contain $T$ and have the same special fibre. This implies that $g_0B_0g_0^{-1}=B$. We have $g_0\phi(\Lie(T))=g_0(\Lie(T_0))=\Lie(g_0Tg_0^{-1})=\Lie(T)$. Thus as $\mu$ factors through $T$, the triple $(M,g_0\phi,T)$ is a $p$-divisible object with a reductive group. Moreover, $(g_0\phi)(\Lie(B))\subseteq g_0(\Lie(B_0))=\Lie(B)$.

\medskip\smallskip\noindent
{\bf 2.6. Some $\dbZ_p$ structures.} Until Section 3 we assume $k=\bar k$. Let $M_{\dbZ_p}$, $G_{\dbZ_p}$, and $(t_{\alpha})_{\alpha\in\scrJ}$ be as in the proof of Claim 2.2.2. The pair $(M_{\dbZ_p},(t_{\alpha})_{\alpha\in\scrJ})$ is a $\dbZ_p$ structure of $(M,(t_{\alpha})_{\alpha\in\scrJ})$. The difference between two such $\dbZ_p$ structures is measured by a torsor $\grT$ of $G_{\dbZ_p}$. As $G_{\dbF_p}$ is connected, Lang theorem (see [39, p. 132]) implies that $\grT_{\dbF_p}$ is trivial. Thus $\grT$ is trivial (as $G_{\dbZ_p}$ is smooth). Thus the isomorphism class of the triple $(M_{\dbZ_p},G_{\dbZ_p},(t_{\alpha})_{\alpha\in\scrJ})$ is an invariant of $(\scrC_g)_{g\in G(W(k))}$. To be short, we refer to $(M_{\dbZ_p},G_{\dbZ_p},(t_{\alpha})_{\alpha\in\scrJ})$ as the $\dbZ_p$ structure of $(M,G,(t_{\alpha})_{\alpha\in\scrJ})$ defined by $\phi\mu(p)$. As $\sigma_0$ acts as $1_{\scrT(M_{\dbZ_p})}\otimes\sigma$ on $\scrT(M):=\scrT(M_{\dbZ_p})\otimes_{\dbZ_p} W(k)$ and as $\sigma$ on $W(k)$-valued points of subgroup schemes of $\GL_{M_{\dbZ_p}}$, from now on we denote $\sigma_0$ by $\sigma$. Also, if $W_{\dbZ_p}$ is an arbitrary free $\dbZ_p$-module, we denote $1_{W_{\dbZ_p}}\otimes\sigma$ by $\sigma$.

If $b\in \GL_M(B(k))$, then $\sigma(b):=\sigma b\sigma^{-1}$ and $\sigma(\mu)$ is the cocharacter of $G$ such that we have $\sigma(\mu)(p^n)=\sigma(\mu(p^n))$ for all $n\in\dbZ$. For $g$ and $g_1\in G(W(k))$, the two triples $\scrC_{g}=(M,g\sigma\mu({1\over p}),G)$ and $\scrC_{g_1}=(M,g_1\sigma\mu({1\over p}),G)$ are rational inner isomorphic if and only if the $\sigma$-conjugacy classes of the elements $g\sigma(\mu)({1\over p})$ and $g_1\sigma(\mu)({1\over p})$ of $G(B(k))$ are equal (i.e., if and only if there exists $h\in G(B(k))$ such that $hg\sigma(\mu)({1\over p})\sigma(h^{-1})=g_1\sigma(\mu)({1\over p})$).

\medskip\noindent
{\bf 2.6.1. Quotients.} We consider an affine, integral group scheme $G_{\dbZ_p}/H_{\dbZ_p}$ that is the geometric quotient of $G_{\dbZ_p}$ through a flat, normal, closed subgroup scheme $H_{\dbZ_p}$ of $G_{\dbZ_p}$. The fibres of $G_{\dbZ_p}/H_{\dbZ_p}$  are reductive group schemes (cf. [10, Vol. III, Exp. XIX, Subsect. 1.7]). Thus $G_{\dbZ_p}/H_{\dbZ_p}$ is smooth as well as a reductive group scheme. Let $H:=H_{W(k)}$. As $G_{\dbZ_p}/H_{\dbZ_p}$ is of finite type over $\dbZ_p$, there exists a free $\dbZ_p$-module $W_{\dbZ_p}$ of finite rank such that we have a closed embedding homomorphism $G_{\dbZ_p}/H_{\dbZ_p}\hookrightarrow \GL_{W_{\dbZ_p}}$ (cf. [10, Vol. I, Exp. VI${}_B$, Rm. 11.11.1]). Let $\mu_{G/H}$ be the composite of $\mu$ with the epimorphism $G\twoheadrightarrow G/H$. Let $W:=W_{\dbZ_p}\otimes_{\dbZ_p} W(k)$. The triple $(W,\sigma\mu_{G/H}({1\over p}),G/H)$ is a $p$-divisible object with a reductive group and we refer to it as a {\it quotient} of $\scrC_{1_M}$.

\smallskip
Let $g\in G(W(k))$. Let $t\in Z^0(G)(W(k))$. Let $t_1\in Z^0(G)(W(k))=Z^0(G_{\dbZ_p})(W(k))$ be such that $t=t_1\sigma(t_1)^{-1}$, cf. [33, Prop. 2.1]. Thus $t_1g{\phi}t_1^{-1}=tg\phi$. We get:

\medskip\noindent
{\bf 2.6.2. Fact.} {\it The element $t_1$ is an inner isomorphism between $\scrC_{tg}$ and $\scrC_g$. If $G$ is a torus, then $\scrC_g$ and $\scrC_{1_M}$ are inner isomorphic and thus the quasi-cocharacters $\nu_g$ and $\nu_{1_M}$ of $G(B(k))$ are $G(B(k))$-conjugate and therefore coincide.}

\medskip\noindent
{\bf 2.6.3. Fact.} {\it The composite of $\nu_g$ with the epimorphism $G_{B(k)}\twoheadrightarrow G^{\ab}_{B(k)}$ does not depend on $g\in G(W(k))$. Moreover, for each $g\in G(W(k))$, there exists $\tilde g\in G^{\der}(W(k))$ such that $\scrC_g$ and $\scrC_{\tilde g}$ are inner isomorphic.}

\medskip
\proof 
Let $g^{\ab}$ be the image of $g$ in $G^{\ab}(W(k))$. By taking $H_{\dbZ_p}=G^{\der}_{\dbZ_p}$ and thus $G/H=G^{\ab}$, the first part follows from  Fact 2.6.2 applied to $(W,\sigma\mu_{G^{\ab}}({1\over p}),G^{\ab})$ and $t=g^{\ab}$. Let $t_1\in G^{\ab}(W(k))$ be such that $g^{\ab}\sigma\mu_{G^{\ab}}({1\over p})=t_1\sigma\mu_{G^{\ab}}({1\over p})t_1^{-1}$. Let $g_1\in P(W(k))$ be an element that maps to $t_1$. The element $\tilde g:=g_1g\phi g_1^{-1}\phi^{-1}\in G(W(k))$ (cf. Lemma 2.4 (a) and (b)) maps to the identity element of $G^{\ab}(W(k))$. Thus we have $\tilde g\in G^{\der}(W(k))$. As $g_1$ is an inner isomorphism between $\scrC_g$ and $\scrC_{\tilde g}$, the last part of the Fact holds as well. \endproof

\medskip
Two quasi-cocharacters in $\Xi(Z^0(G_{B(k)}))$ are equal if and only if their images in $\Xi(G_{B(k)}^{\ab})$ are equal. Thus from Corollary 2.3.2 and Fact 2.6.3 we get that:
 
\medskip\noindent
{\bf 2.6.4. Corollary.} {\it If  $g_1$, $g_2\in G(W(k))$ are such that $\scrC_{g_1}$ and $\scrC_{g_2}$ are basic, then $\nu_{g_1}=\nu_{g_2}$.}

\medskip\noindent
{\bf 2.6.5. Covers.} Let $q\in\dbN$ be such that $G_{\dbZ_p}$ has a maximal torus that splits over $W(\dbF_{p^q})$. A result of Langlands (see [32, pp. 297--299]) implies that there exists an epimorphism $e_{\dbQ_p}:G^\flat_{\dbQ_p}\twoheadrightarrow G_{\dbQ_p}$ such that $G^{\flat\der}_{\dbQ_p}$ is an a priori given isogeny cover of $G^{\der}_{\dbQ_p}$ and $\Ker(e_{\dbQ_p})$ is a product of tori of the form $\Res_{B(\dbF_q)/\dbQ_p} \dbG_m$. Thus $G^\flat_{\dbQ_p}$ is quasi-split and splits over an unramified extension of $\dbQ_p$. Let $G^{\flat}_{\dbZ_p}$ be a reductive group scheme over $\dbZ_p$ that has $G^{\flat}_{\dbQ_p}$ as its generic fibre, cf. [43, Subsubsects. 1.10.2 and 3.8.1]. The Zariski closure $\Ker(e)$ of $\Ker(e_{\dbQ_p})$ in $G^{\flat}_{\dbZ_p}$ is a subtorus of $Z^0(G^\flat_{\dbZ_p})$. As two hyperspecial subgroups of $G_{\dbZ_p}(\dbQ_p)$ are $G_{\dbZ_p}^{\ad}(\dbQ_p)$-conjugate (see [43, p. 47]), we can assume $(G^{\flat}_{\dbZ_p}/\Ker(e))(\dbZ_p)=G_{\dbZ_p}(\dbZ_p)$. Thus we can identify $G^{\flat}_{\dbZ_p}/\Ker(e)=G_{\dbZ_p}$, cf. [43, Subsubsect. 3.8.1]. The epimorphism $e:G^{\flat}_{\dbZ_p}\twoheadrightarrow G_{\dbZ_p}$ extends $e_{\dbQ_p}$. Let $G^\flat:=G^\flat_{W(k)}$. As $\Ker(e)$ is a torus, $\mu$ lifts to a cocharacter $\mu_{G^\flat}:\dbG_m\to G^\flat$. Let $W^\flat_{\dbZ_p}$ be a free $\dbZ_p$-module of finite rank such that we have a closed embedding homomorphism $G^\flat_{\dbZ_p}\hookrightarrow \GL_{W_{\dbZ_p}^\flat}$. Let $W^\flat:=W^\flat_{\dbZ_p}\otimes_{\dbZ_p} W(k)$. We call the $p$-divisible object with a reductive group $(W^\flat,\sigma\mu_{G^\flat}({1\over p}),G^\flat)$ a {\it cover} of $\scrC_{1_M}$. 

There exists a $G^\flat_{\dbQ_p}$-epimorphism $\scrT(W^\flat_{\dbZ_p}[{1\over p}])\twoheadrightarrow M_{\dbZ_p}[{1\over p}]$, cf. [53, Thm. 3.5] and Weyl complete reductibility theorem. Thus $\nu_{1_M}$ is the composite of the Newton quasi-cocharacter of $(W^\flat,\sigma\mu_{G^\flat}({1\over p}),G^\flat)$ with $e_{B(k)}$. 
From Corollary 2.3.2 we get that the cover $(W^\flat,\sigma\mu_{G^\flat}({1\over p}),G^\flat)$ is basic if and only if $\scrC_{1_M}$ is so.

\medskip\smallskip\noindent
{\bf 2.7. Proposition.} {\it We recall that $k=\bar k$. Let $g_1$, $g_2\in G(W(k))$ be such that $\scrC_{g_1}$ and $\scrC_{g_2}$ are basic. Then there exists $h\in G(B(k))$ such that $hg_1\phi=g_2{\phi}h$.}

\medskip
\proof We appeal to Subsubsection 2.6.5. As $\Ker(e)$ is a torus, there exists $g_i^\flat\in G^\flat(W(k))$ that maps to $g_i$. Thus $(W^\flat,g_i^\flat\sigma\mu_{G^\flat}({1\over p}),G^\flat)$ is basic, cf. end of Subsubsection 2.6.5 (applied to $g_i^\flat\sigma$ instead of $\sigma$). So to simplify the presentation we can assume that $G^\flat=G$ and that $G^{\der}$ is simply connected; thus $g_i^\flat=g_i$. Let $G_{1\dbQ_p}$ be the reductive group over $\dbQ_p$ which represents (as in [38, Subsects. 1.12 to 1.15]) the functor that associates to a commutative $\dbQ_p$-algebra $R$ the group 
$$G_{1\dbQ_p}(R)=\{g_{11}\in G(B(k)\otimes_{\dbQ_p} R)|g_{11}g_1(\sigma\otimes 1_R)\mu({1\over p})=g_1(\sigma\otimes 1_R)\mu({1\over p})g_{11}\};\leqno (7)$$
the second equality of (7) is between $\sigma\otimes 1_R$-linear automorphisms of $M[{1\over p}]\otimes_{\dbQ_p} R$.
We have $\nu_{g_1}=\nu_{g_2}$, cf. Corollary 2.6.4. Thus based on [37, Prop. 1.17] we can speak about the class $\delta\in H^1(\dbQ_p,G_{1\dbQ_p})$ that defines the left torsor $G_{12\dbQ_p}$ of $G_{1\dbQ_p}$ defined by the rule:
$G_{12\dbQ_p}(R)=\{g_{12}\in G(B(k)\otimes_{\dbQ_p} R)|g_{12}g_1(\sigma\otimes 1_R)\mu({1\over p})=g_2(\sigma\otimes 1_R)\mu({1\over p})g_{12}\}$. 
The image of $\delta$ in $H^1(\dbQ_p,G^{\ab}_{1\dbQ_p})$ is the identity element (cf. Fact 2.6.2) and the only class in $H^1(\dbQ_p,G^{\der}_{1\dbQ_p})$ is the trivial one (as $G^{\der}_{1\dbQ_p}$ is simply connected, cf. [21]). Thus $\delta$ is the trivial class i.e., there exists $h\in G_{12\dbQ_p}(\dbQ_p)\subseteq G(B(k))$. For such an $h$ we have $hg_1\phi=g_2{\phi}h$. \endproof

\bigskip\smallskip
\noindent
{\boldsectionfont 3. Standard forms}

\bigskip
Let $\scrC_{1_M}=(M,\phi,G)$ be a Shimura $p$-divisible object over $k$; we have $b_L\Le 1$. Let $g\in G(W(k))$ be such that $\scrC_g$ is not basic. Let $P_G^+(g\phi)$ and $P_G^-(g\phi)$ be as in Definition 2.3.3; none of these two parabolic subgroup schemes of $G$ is $G$ itself (cf. Corollary 2.3.2). In this Section we study $\scrC_g$ up to inner isomorphisms. The goal is to present standard forms of $\scrC_g$. Properties of them are listed. In Subsection 3.1 we work with $P_G^+(g\phi)$ and in Subsection 3.2 we translate Subsection 3.1 in the context of $P_G^-(g\phi)$. 

\medskip\smallskip\noindent
{\bf 3.1. Non-negative standard forms.}  
The intersection $P_k\cap P_G^+(g\phi)_k$ contains a maximal split torus $\tilde T_{k}$ of $G_k$, cf. [3, Ch. V, Prop. 20.7 (i)]. Let $\tilde T$ (resp. $\tilde T^+$) be a split torus of $P$ (resp. of $P_G^+(g\phi)$) that lifts $\tilde T_k$, cf. [10, Vol. II, Exp. IX, Thms. 3.6 and 7.1]. Let $g_3\in\Ker(G(W(k))\to G(k))$ be such that $g_3\tilde T^+g_3^{-1}=\tilde T$, cf. loc. cit. Let $g_2:={\phi}g_3^{-1}\phi^{-1}$. It is an element of $G(W(k))$, cf. Lemma 2.4 (b) and (c). By replacing $g\phi$ with $g_3g{\phi}g_3^{-1}=g_3gg_2\phi$ we can assume that $\tilde T=\tilde T^+$. But $\tilde T$ is $P(W(k))$-conjugate to a maximal split torus of $P$ that contains $\im(\mu)$ (cf. loc. cit. and [3, Ch. IV, Thm. 15.14]). Thus by replacing $\mu$ with a $P(W(k))$-conjugate of it (cf. Fact 2.4.1 and Lemma 2.4 (a)), we can also assume that $\mu$ factors through $\tilde T=\tilde T^+$ and therefore also through $P_G^+(g\phi)$. 

\medskip\noindent
{\bf 3.1.1. Unipotent considerations.} Let $U_G^+(g\phi)$ be the unipotent radical of $P_G^+(g\phi)$. Let $L^+$ be a Levi subgroup scheme of $P_G^+(g\phi)$ that contains $\tilde T$. We identify $L^+$ with $P_G^+(g\phi)/U_G^+(g\phi)$. As $\mu(p)$ and $g\phi$ normalize $\Lie(P_G^+(g\phi)_{B(k)})$, the $\sigma$-linear automorphism $g\phi\mu(p)$ of $M$ normalizes $\Lie(P_G^+(g\phi))=\Lie(P_G^+(g\phi)_{B(k)})\cap\End(M)$. Let $L^+_0$ be the Levi subgroup scheme of $P_G^+(g\phi)$ that has $g\phi\mu(p)(\Lie (L^+))$ as its Lie algebra. There exists $g_1\in P_G^+(g\phi)(W(k)) $ such that $g_1(L^+_0)g_1^{-1}$ and $L^+$ coincide mod $p$, cf. [3, Ch. V, Prop. 20.5]. By replacing $g_1$ with a $\Ker(P_G^+(g\phi)(W(k))\to P_G^+(g\phi)(k))$-multiple of it, we can assume that $g_1(L^+_0)g_1^{-1}$ and $L^+$ have a common maximal torus $T^+$ (cf. the infinitesimal liftings of [10, Vol. II, Exp. IX, Thm 3.6]). Thus $g_1(L^+_0)g_1^{-1}=L^+$. [Argument: we can assume $T^+$ is split and it suffices to show that the root systems of the inner conjugation actions of $T^+$ on $\Lie(g_1(L^+_0)g_1^{-1})$ and $\Lie(L^+)$ coincide; but this holds as it holds mod $p$.] 

Both $\Lie(L^+_{B(k)})$ and $\Lie(P_G^+(g\phi)_{B(k)})$ are normalized by $g_1g\phi$. We write $g_1^{-1}=u^+l^+$, where $u^+\in U_G^+(g\phi)(W(k))$ and $l^+\in L^+(W(k))$. Let $\phi_1:=(u^+)^{-1}g\phi=l^+g_1g\phi$. We have
$$g\phi=u^+\phi_1.\leqno (8)$$
We refer to (8) as a {\it non-negative standard form} of $\scrC_g$. Here are its main properties.

\medskip\noindent
{\bf 3.1.2. Theorem.} {\it We recall that $\scrC_g=(M,g\phi,G)$ is not basic, that $L^+$ is a Levi subgroup scheme of $P_G^+(g\phi)$, and that (up to inner isomorphism we can assume) $\mu$ factors through a maximal split torus $\tilde T$ of $L^+\cap P$. Referring to (8), the following four properties hold: 
\medskip
{\bf (a)} Both Lie algebras $\Lie(P_G^+(g\phi)_{B(k)})$ and $\Lie(L^+_{B(k)})$ are normalized by $\phi_1$.

\smallskip
{\bf (b)} The triple $(M,\phi_1,L^+)$ is a basic Shimura $p$-divisible object. 

\smallskip
{\bf (c)} There exists a rational inner isomorphism between $\scrC_g$ and $\scrC_{(u^+)^{-1}g}=(M,\phi_1,G)$ defined by an element $h\in U_G^+(g\phi)(B(k))$. Thus the Newton quasi-cocharacters $\nu_g$ and $\nu_{(u^+)^{-1}g}$ (see Subsubsections 1.2.4 and 2.2.3 for notations) are $U_G^+(g\phi)(B(k))$-conjugate. 

\smallskip
{\bf (d)} The Levi subgroup $C_G((u^+)^{-1}g\phi)$ of $\scrC_{(u^+)^{-1}g}$ is $L^+_{B(k)}$.}

\medskip
\proof 
As $u^+$ and $g\phi$ normalize $\Lie(P_G^+(g\phi)_{B(k)})$, so does $\phi_1$. 
As $l^+$ and $g_1g\phi$ normalize $\Lie(L^+_{B(k)})$, so does $\phi_1=l^+g_1g\phi$. Thus (a) holds. 

As $\phi_1$ normalizes $\Lie(L^+_{B(k)})$ and as $\mu$ factors through the torus $\tilde T$ of $L^+$, the triple $(M,\phi_1,L^+)$ is a Shimura $p$-divisible object. The actions of $g\phi=u^+\phi_1$ and $\phi_1$ on the quotient $B(k)$-vector space $\Lie(P_G^+(g\phi)_{B(k)}))/\Lie(U_G^+(g\phi)_{B(k)})$ are the same. Thus the Newton polygon slopes of $(\Lie(L^+_{B(k)}),\phi_1)$ are equal to the Newton polygon slopes of the action of $g\phi=u^+\phi_1$ on $\Lie(P_G^+(g\phi)_{B(k)}))/\Lie(U_G^+(g\phi)_{B(k)})$ and therefore are all $0$. Thus $(M,\phi_1,L^+)$ is basic i.e., (b) holds.

We prove (c). We have to show the existence of an element $h\in U_G^+(g\phi)(B(k))$ such that $\phi_1=hu^+{\phi_1}h^{-1}$. As $\scrC_g$ is not basic, the set $\grL\grS^+$ of Newton polygon slopes of $(\Lie(U_G^+(g\phi)_{B(k)}),\phi_1)$ is non-empty. As the group scheme $U_G^+(g\phi)$ is unipotent, $\Lie(U_G^+(g\phi))$ has a characteristic series whose factors are fixed by $U_G^+(g\phi)$. This implies that the Newton polygons of $(\Lie(U_G^+(g\phi)_{B(k)}),\phi_1)$ and $(\Lie(U_G^+(g\phi)_{B(k)}),u^+\phi_1)$ coincide. Thus $\grL\grS^+\in (0,\infty)$, cf. Lemma 2.3.1 applied to $\scrC_g=(M,u^+\phi_1,G)$. 

Let $(\Lie(U_G^+(g\phi)_{B(k)}),\phi_1)=\oplus_{\eta\in \grL\grS^+} (N_{\eta},\phi_1)$ be the Newton polygon slope decomposition. Let $\eta_1<\eta_2<\cdots<\eta_q$ be the elements in $\grL\grS^+$ ordered increasingly. As $\grL\grS^+\in (0,\infty)$, we have $0<\eta_1$. For $i\in S(1,q)$ we check that there exists a connected, unipotent subgroup $U_i$ of $U_G^+(g\phi)_{B(k)}$ such that $\Lie(U_i)=\oplus_{\eta\in \grL\grS^+\cap [\eta_i,\infty)} N_{\eta}$. 

The uniqueness of $U_i$ is implied by [3, Ch. II, Subsect. 7.1]. As $\oplus_{\eta\in \grL\grS^+\cap [\eta_i,\infty)} N_{\eta}$ is an ideal of $\Lie(P_G^+(g\phi)_{B(k)})$, the direct sum $S_i:=\oplus_{\eta\in \grL\grS^+\cap [\eta_i,\infty)} N_{\eta}\oplus\Lie(T^+_{B(k)})$ is a Lie subalgebra of $\Lie(P_G^+(g\phi)_{B(k)})$. As $[\Lie(T^+_{B(k)}),N_{\eta}]=N_{\eta}$, we have $[S_i,S_i]=\oplus_{\eta\in \grL\grS^+\cap [\eta_i,\infty)} N_{\eta}$. Thus the existence of $U_i$ is implied by [3, Ch. II, Cor. 7.9]. Let $U_{q+1}$ be the trivial subgroup of $U_G^+(g\phi)_{B(k)}$. For $i\in S(1,q)$ we have $U_{i+1}\vartriangleleft U_i\vartriangleleft P_G^+(g\phi)_{B(k)}$. 

By induction on $i\in S(1,q)$ we show that there exists an element $h_i\in U_i(B(k))$ such that 
we have $h_ih_{i_1}\cdots h_1u^+{\phi_1}h_1^{-1}\cdots h^{-1}_i\phi_1^{-1}\in U_{i+1}(B(k))$. To ease notations we will only show the existence of $h_1$ (the existence of $h_2,\ldots,h_q$ is argued entirely in the same way). The quotient group $U_1/U_2$ is unipotent and commutative and thus we can identify it with the affine vector group  defined by $\Lie(U_1/U_2)$. We also identify $\Lie(U_1/U_2)=N_{\eta_1}$. Thus $u^+\in U_1(B(k))$ modulo $U_2(B(k))$ is identified with an element $\tilde u_1\in N_{\eta_1}$. We fix a finite sequence $(\dbG_a(s))_{s\in I_1}$ of $\dbG_a$ subgroups of $U_1$ normalized by $\tilde T_{B(k)}$ and such that $\oplus_{s\in I_1} \Lie(\dbG_a(s))=N_{\eta_1}$. For $j\in\dbN\cup\{0\}$ let $h_{1j}\in U_1(B(k))$ be such that we have:
\medskip

{\bf (i)} its image in $(U_1/U_2)(B(k))$ is identified with $(-1)^{j+1}\phi_1^j(\tilde u_1)$, and

\smallskip
{\bf (ii)} it is a product of $B(k)$-valued points of these $\dbG_a(s)$ subgroups of $U_1$ (taken with respect to a fixed total ordering of $I_1$).

\medskip
As $\eta_1>0$ the sequence $(\phi_1^j(\tilde u_1))_{j\in\dbN\cup\{0\}}$ converges to $0$ in the $p$-adic topology of $N_{\eta_1}$. Due to this and (ii) it makes sense to define $h_1:=(\prod_{j=0}^{\infty} h_{1j}^{-1})^{-1}\in U_1(B(k))$.
But $h_{1j}\cdots h_{11}h_{10}u^+\phi_1h_{10}^{-1}\cdots h_{1j}^{-1}\phi_1^{-1}\in U_1(B(k))$ modulo $U_2(B(k))$
is $(-1)^{j}\phi_1^{j+1}(\tilde u_1)$ for all $j\in\dbN\cup\{0\}$. Therefore $h_1u^+\phi_1h_1^{-1}\phi_1^{-1}\in U_1(B(k))$ modulo $U_2(B(k))$ is the $p$-adic limit of the sequence $((-1)^{j}\phi_1^{j+1}(\tilde u_1))_{j\in\dbN\cup\{0\}}$ and thus it is $0$. Thus $h_1u^+\phi_1h_1^{-1}\phi_1^{-1}\in U_2(B(k))$. 

This ends the induction. If $h:=h_qh_{q-1}\cdots h_1\in U_G^+(g\phi)(B(k))$, then we have $hu^+{\phi_1}h^{-1}\phi_1^{-1}=1_M$ and therefore $h$ defines a rational inner isomorphism between $\scrC_g$ and $\scrC_{(u^+)^{-1}g}$. Thus (c) holds.

We prove (d). We have $L^+_{B(k)}\leqslant C_G((u^+)^{-1}g\phi)$, cf. (b). As $C_G((u^+)^{-1}g\phi)$ and $C_G(g\phi)$ are isomorphic (cf. (c)) and as $L^+_{B(k)}$ and $C_G(g\phi)$ are isomorphic, $C_G((u^+)^{-1}g\phi)$ and $L^+_{B(k)}$ are isomorphic. By reasons of dimensions we get $C_G((u^+)^{-1}g\phi)=L^+_{B(k)}$. \endproof

\medskip\smallskip\noindent
{\bf 3.2. Non-positive standard forms.} Let $U_G^-(g\phi)$ be the unipotent radical of $P_G^-(g\phi)$. There exists a maximal split torus $\tilde T^-$ of $P_G^-(g\phi)$ and a Levi subgroup scheme $L^-$ of $P_G^-(g\phi)$ that contains $\tilde T^-$, such that up to inner isomorphisms we can assume that $\mu$ factors through $\tilde T^-$ and that we can write $g\phi=u^-\phi_2$, where $u^-\in U_G^-(g\phi)(W(k))$ and $(M,\phi_2,L^-)$ is a basic Shimura $p$-divisible object. The proof of this is entirely the same as the proof of Theorem 3.1.2 (c). The only significant difference is: in connection to property 3.1.2 (i) we have to work in ``reverse order'' i.e., we have to replace $(-1)^{j+1}\phi_1^j(\tilde u_1)$ by $(-1)^{j}\phi_2^{-j-1}(\tilde u_1)$.

\bigskip\smallskip
\noindent
{\boldsectionfont 4. Proofs of the basic results and complements}

\bigskip
In Subsection 4.1 we prove Theorem 1.3.1 and Corollary 1.3.2. In Subsection 4.2 we prove Theorem 1.3.3. Subsections 4.3 to 4.7 present some conclusions and complements to Theorem 1.3.3. See Subsection 4.6 for examples. See Subsection 4.7 for the quasi-polarized context. Let $\scrC_{1_M}=(M,\phi,G)$ be a Shimura $p$-divisible object over $k$; thus $b_L\Le 1$. We use the notations of the last paragraph of Subsection 2.1. Let $g\in G(W(k))$. Let $P_G^+(g\phi)$ be as in Definition 2.3.3.

\medskip\smallskip\noindent
{\bf 4.1. Proofs of 1.3.1 and 1.3.2.}
We know that the first part of Theorem 1.3.1 holds, cf. Lemma 2.3.1. We prove the second part of Theorem 1.3.1 (for $b_L\le 1$). If $\scrC_g$ is basic, then $P_G^+(g\phi)=G$ (cf. Corollary 2.3.2) and thus there exist Hodge cocharacters of $\scrC_g$ that factor through $P_G^+(g\phi)$. If $\scrC_g$ is not basic, then there exists a Hodge cocharacter of $\scrC_g$ that factors through $P_G^+(g\phi)$, cf. first paragraph of Subsection 3.1. This proves Theorem 1.3.1.

We prove Corollary 1.3.2 (a). Let $P_{\GL_M}^+(g\phi)$ be the non-negative parabolic subgroup scheme of $(M,g\phi,\GL_M)$. It is the parabolic subgroup scheme of $\GL_M$ that normalizes $W^{\gamma}(M,g\phi)$ for all $\gamma\in \grS(g)$. Obviously $\Lie(P_G^+(g\phi)_{B(k)})\subseteq\Lie(P_{\GL_M}^+(g\phi)_{B(k)})$. Thus $P_G^+(g\phi)_{B(k)}$ is a subgroup of $P_{\GL_M}^+(g\phi)_{B(k)}$, cf. [3, Ch. II, Subsect. 7.1]. This implies that $P_G^+(g\phi)$ is a subgroup scheme of $P_{\GL_M}^+(g\phi)$. Therefore each Hodge cocharacter of $\scrC_g$ that factors through $P^+_G(g\phi)$, normalizes $W^{\gamma}(M,g\phi)$ for all $\gamma\in \grS(g)$. Thus Corollary 1.3.2 (a) follows from Theorem 1.3.1. 

We prove Corollary 1.3.2 (b). If $\scrC_g$ is basic, then the Newton cocharacter of $\scrC_g$ factors through $Z^0(G_{B(k)})$ (cf. Corollary 2.3.2) and therefore it extends to a cocharacter of $Z^0(G)$. Thus we have a direct sum decomposition $M=\oplus_{\gamma\in \grS(g)} M\cap M_{\gamma}(g)$ normalized by all cocharacters of $\GL_M$ that commute with $Z^0(G)$ and therefore also by all Hodge cocharacters  of $\scrC_g$. 

If $\scrC_g$ is not basic, then up to a rational inner isomorphism we can assume there exists a Levi subgroup scheme $L^+$ of $P_G^+(g\phi)$ such that $(M,g\phi,L^+)$ is a basic Shimura $p$-divisible object (cf. Theorem 3.1.2 (b) and (c)). As a Hodge cocharacter of $(M,g\phi,L^+)$ is also a Hodge cocharacter of $\scrC_g$, we can apply the basic part of the previous paragraph to $(M,g\phi,L^+)$. This proves Corollary 1.3.2 (b) and ends the proofs of Theorem 1.3.1 and Corollary 1.3.2.\endproof

\medskip\noindent
{\bf 4.1.1. Remark.} Replacing non-negative by non-positive and the references to Subsection 3.1 by ones to Subsection 3.2, as in the first paragraph of Subsection 4.1 we argue that there exist Hodge cocharacters of $\scrC_g$ that factor through $P_G^-(g\phi)$.

\medskip\smallskip\noindent
{\bf 4.2. Proof of 1.3.3.} Until Corollary 4.3 we assume that $G$ is split and that $b_L\Le 1$. Let $T$ and $B$ be as in Subsection 2.5. Thus $\mu$ factors through the maximal torus $T$. We can assume that $\phi(\Lie(T))=\Lie(T)$ and that $\phi(\Lie(B))\subseteq\Lie(B)$, cf. Subsection 2.5. Let $N$ be the normalizer of $T$ in $G$. If $\tilde w\in N(W(k))$, then $\tilde w\phi$ normalizes $\Lie(T)$. Thus $(M,\tilde w\phi,T)$ is a Shimura $p$-divisible object. From Fact 2.6.2 we get:

\medskip\noindent
{\bf 4.2.1. Fact.} {\it If $k=\bar k$, then for each $t\in T(W(k))$ and $\tilde w\in N(W(k))$, there exists an inner isomorphism between $\scrC_{t\tilde w}$ and $\scrC_{\tilde w}$ defined by a suitable element $t_1\in T(W(k))$.}

\medskip
We prove Theorem 1.3.3 in two steps. The first step deals with the basic context and introduces notations we will also use in Subsections 4.5 and 4.6 and in Section 5 (see Subsubsection 4.2.2). The second step achieves the reduction to the basic context (see Subsubsection 4.2.3).

\medskip\noindent
{\bf 4.2.2. Step 1.} To show that there exists $w\in N(W(k))$ such that $\scrC_w$ is basic, we can assume that $k=\bar k$ (cf. Fact 4.2.1); let $M_{\dbZ_p}$ and $G_{\dbZ_p}$ be obtained as in the proof of Claim 2.2.2. We can also assume that $G_{\dbZ_p}$ is adjoint (cf. Subsubsection 2.6.1 applied with $H_{\dbZ_p}=Z(G_{\dbZ_p})$) and $\dbZ_p$-simple and that $b_L=1$. Let $\scrL$
be the Lie type of an arbitrary simple factor of $G$. It is well known that $\scrL$ is not $E_8$, $F_4$, or $G_2$ (the maximal roots of these Lie types have all coefficients different from 1; see [5, plates VII to IX] and [8, Table 1.3.9]). Let $r$ be the rank of $\scrL$. We write $G_k=\prod_{i\in I} G_{ik}$ as a product of absolutely simple, adjoint groups over $k$ (cf. [42, Subsubsect. 3.2.1]). This decomposition lifts to a product decomposition 
$$G=\prod_{i\in I} G_i$$
into absolutely simple, adjoint group schemes over $W(k)$ (cf. [10, Vol. III, Exp. XXIV, Prop. 1.21]). Let $T_i:=T\cap G_i$. Let $B_i:=B\cap G_i$. Let $\Lie(G_i)=\Lie(T_i)\bigoplus_{\alpha\in\Phi_i}\grg_{i\,\alpha}$ be the root decomposition with respect to the split, maximal torus $T_i$ of $\Lie(G_i)$. As $G_{\dbZ_p}$ is $\dbZ_p$-simple, there exists $n\in\dbN$ and an absolutely simple adjoint group scheme $J$ over $W(\dbF_{p^n})$ such that $G_{\dbZ_p}=\Res_{W(\dbF_{p^n})/\dbZ_p}J$, cf. [42, Subsubsect. 3.1.2] and [10, Vol. III, Exp. XXIV, Prop. 1.21]. As $\sigma$ permutes transitively the simple factors of $G$, $\phi$ permutes transitively the $\Lie(G_i)[{1\over p}]$'s. Thus we can assume that $I=S(1,n)$ and that for all $i\in I$ we have $\phi(\Lie(G_i))\subseteq G_{i+1}[{1\over p}]$, where $G_{n+1}:=G_1$. 

Let $\Phi_i^+:=\{\alpha\in\Phi_i|\grg_{i\,\alpha}\in\Lie(B_i)\}$. Let $\Delta_i\subseteq\Phi_i^+$ be the basis of $\Phi_i$ contained in $\Phi_i^+$. Let $\alpha_1(i),\ldots,\alpha_r(i)\in\Delta_i$ be denoted as in [5, plates I to VII] except that we put a lateral right index $(i)$. As $\phi(\Lie(B))\subseteq\Lie(B)$ we have $\phi^n(\Lie(B_i))\subseteq\Lie(B_i)$. Thus there exists a permutation $\pi_i$ of $\Phi_i^+$ that normalizes $\Delta_i$ and such that we have $\phi^n(\grg_{i\,\alpha})\subseteq \grg_{i\,\pi_i(\alpha)}[{1\over p}]$ for all $\alpha\in\Phi_i$. The order of $\pi_i$ does not depend on $i\in I$ and therefore we denote it by $o$. We have $o\in\{1,2,3\}$, cf. the structure of the group of automorphisms of the Dynkin diagram of $G_i$ (it is either trivial or $\dbZ/2\dbZ$ or $S_3$; see loc. cit.). Let $w_0\in (N\cap G_1)(W(k))$ be such it takes $B_1$ to its opposite with respect to $T_1$; see item XI of [5, plates I to VII]. The following six Cases recall the ``most practical'' elliptic element of $(N\cap G_1)(W(k))$. 

\smallskip\noindent
{\bf Case 1.} Suppose that $\scrL$ is $B_r$, $C_r$, or $E_7$ ($r\in\dbN$). Thus $o=1$ and $w_0$ takes $\alpha\in\Phi_1$ to $-\alpha$, cf. loc. cit. Let $n_1:=n$. For $\alpha\in\Phi_1$ we have $(w_0\phi)^{n_1}(\grg_{1\,\alpha})=p^{s_{\alpha}}\grg_{1\,-\alpha}$, with $s_{\alpha}\in\dbZ$. As $w_0\phi$ normalizes $\Lie(T)$, all Newton polygon slopes of $(\Lie(T),w_0\phi)$ are $0$. As we have $s_{-\alpha}=-s_{\alpha}$ for all $\alpha\in\Phi_1$, all Newton polygon slopes of $(\Lie(G_{B(k)})/\Lie(T_{B(k)}),w_0\phi)$ are also $0$. Thus all Newton polygon slopes of $(\Lie(G_{B(k)}),w_0\phi)$ are $0$. 

\smallskip\noindent
{\bf Case 2.} Suppose that either $o=2$ and $\scrL\in\{A_r|r\in\dbN\setminus\{1\}\}\cup\{D_r|r-3\in 2\dbN\}\cup\{E_6\}$ or $o=1$ and $\scrL\in\{D_r|r-2\in 2\dbN\}$. This Case is the same as Case 1.

\smallskip\noindent
{\bf Case 3.} Suppose that $o=3$. Thus $\scrL=D_4$ and $w_0$ takes $\grg_{1\alpha_1(i)}$ to $\grg_{1-\alpha_1(i)}$. For $\alpha\in\Phi_1$ we have $(w_0\phi)^{3n}(\grg_{1\alpha})\subseteq \grg_{1-\alpha}[{1\over p}]$. Thus this Case is the same as Case 1 but with $n_1:=3n$.

\smallskip\noindent
{\bf Case 4.} Suppose that $o=1$ and $\scrL=A_r$, with $r\Ge 2$. The simply connected semisimple group scheme cover $G_1^{\sc}$ of $G_1$ can be naturally identified with $\SL_{M_0}$, where $M_0:=W(k)^{r+1}$. We choose a $W(k)$-basis $\{e_1,\ldots,e_{r+1}\}$ for $M_0$ such that the inverse image of $T_1$ in $G_1^{\sc}$ normalizes $W(k)e_s$ for all $s\in S(1,r+1)$. Let $w$ be the image in $N(W(k))$ of an element $w^{\sc}\in G_1^{\sc}(W(k))=\SL_{M_0}(W(k))$ that maps $W(k)e_s$ onto $W(k)e_{s+1}$ for all $s\in S(1,r+1)$, where $e_{r+2}:=e_1$. As $o=1$ and due to the circular form of $w^{\sc}$, all Newton polygon slopes of $(\Lie(G)[{1\over p}],w\phi)$ are $0$.

\smallskip\noindent
{\bf Case 5.} Suppose that $o=1$ and $\scrL=E_6$. Let $G_0$ be the reductive subgroup scheme of $G$ whose Lie algebra is generated by $\Lie(T)$ and $\grg_{i,\alpha}$'s, where $i\in I$ and $\alpha\in\Phi_i\cap \sum_{j\in\{1,3,4,5,6\}} \dbZ\alpha_i(j)$, cf. [10, Vol. III, Exp. XXII, Thm. 5.4.7 and Prop. 5.10.1]. Thus $G_0^{\ad}$ is a product of $n$ absolutely simple adjoint group schemes of $A_5$ Lie type whose Lie algebras are permuted transitively by $\phi$. The group scheme $G_0$ is the centralizer of $Z^0(G_0)$ in $G$ and we have $T\leqslant G_0$. Thus $\mu$ factors through $G_0$. As $\Lie(G_{0B(k)})$ is normalized by $\phi$, we get that the triple $(M,\phi,G_0)$ is a Shimura $p$-divisible object. As $w_0$ normalizes $G_0\cap G_1$, the triple $(M,w_0\phi,G_0)$ is also a Shimura $p$-divisible object. As $o=1$, from Case 4 we deduce the existence of $w_1\in (G_0\cap G_1\cap N)(W(k))$ such that $(M,w\phi,G_0)$ is basic, where $w:=w_1w_0$. Thus the Newton quasi-cocharacter $\nu_{0w}$ of $(M,w\phi,G_0)$ factors through $Z^0(G_{0B(k)})$. Its composite with the epimorphism $e_0:G_{0B(k)}\twoheadrightarrow G^{\ab}_{0B(k)}$ is the composite $\nu_0$ of the Newton quasi-cocharacter of $(M,w_0\phi,G_0)$ with $e_0$. We have $w_0(\grg_{1\alpha_1(i)})=\grg_{1-\alpha_1(j)}$, where $(i,j)\in\{(1,6),(6,1),(3,5),(5,3),(2,2),(4,4)\}$ (cf. [5, plate V, item (IX)]). Thus the maximal subtorus of $T_1$ fixed by $w_0$  has rank $2$. A similar argument shows that $w_0$ fixes a subtorus of $T_1\cap G^{\der}_0$ of rank $2$. Thus the automorphism of the rank 1 torus $G_1\cap Z^0(G_0)$ induced by $w_0$ is non-trivial. Therefore $\nu_0$ is trivial. Thus $\nu_{0w}$ is trivial i.e., $\scrC_w$ is basic. 

\smallskip\noindent
{\bf Case 6.} Suppose that either $o=2$ and $\scrL\in\{D_r|r-2\in 2\dbN\}$ or $o=1$ and $\scrL\in\{D_r|r-3\in 2\dbN\}$. If $p>2$ let $u=0$, let $y=1$, and let $z\in W(\dbF_{p^n})$ be such that $-z$ mod $p$ is not a square. If $p=2$ let $u=1$ and let $y,z\in W(\dbF_{p^n})$ be such that the quadratic form $yx_{2r-1}^2+x_{2r-1}x_{2r}+zx_{2r}^2$ does not represent $0$. If $r$ is odd (resp. even), then the group scheme $J$ is split (resp. non-split and splits over $W(\dbF_{p^{2n}})$). Thus $J$ is the adjoint group scheme of the $\pmb{\text{SO}}$ group scheme of the quadratic form $x_1x_2+x_3x_4+\cdots+x_{2r-1}x_{2r}$ (resp. $x_1x_2+x_3x_4+\cdots+x_{2r-3}x_{2r-2}+yx_{2r-1}^2+ux_{2r-1}x_{2r}+zx_{2r}^2$) on $W(\dbF_{p^n})^{2r}$, cf. [3, Ch. V, Subsects. 23.4 to 23.6] for the case of special fibres. Let $T_J^0$ be a torus of $J$ of rank $r-3$ that is a product of non-split rank 1 tori and whose centralizer $C_J$ in $J$ is such that $C_J^{\ad}$ is split, absolutely simple of $A_3$ Lie type. The existence of $T_J^0$ is a consequence of the fact that the $\pmb{\text{SO}}$ group scheme of the quadratic form $x_{2i-3}x_{2i-2}+x_{2i-1}x_{2i}$, $i\in 2\dbN\cap S(1,r-3)$, has a maximal torus that is a product of two non-split rank 1 tori. Let $T_J$ be a maximal torus of $C_J$ that contains $T_J^0$. Let $h\in G(W(k))$ be such that $hTh^{-1}$ is $\tilde T:=(\Res_{W(\dbF_{p^n})/\dbZ_p}T_J)_{W(k)}$. Let $\tilde C:=(\Res_{W(\dbF_{p^n})/\dbZ_p}C_J)_{W(k)}$; it is a subgroup scheme of $G$. Let $\tilde \mu:=h\mu h^{-1}$; it is a cocharacter of $\tilde T$. The triple $(M,\sigma\tilde\mu({1\over p}),\tilde C)$ is a Shimura $p$-divisible object (we recall that $\sigma:M\arrowsim M$ fixes $M_{\dbZ_p}$, see Subsection 2.6). Let $\tilde w\in\tilde C(W(k))$ be an element that normalizes $\tilde T$ and such that $(M,\tilde w\sigma\tilde\mu({1\over p}),\tilde C)$ is basic, cf. Case 4. The Newton cocharacter $\tilde\nu_{\tilde w}$ of $(M,\tilde w\sigma\tilde\mu({1\over p}),\tilde C)$ factors through $\tilde T^0:=Z^0(\tilde C)=(\Res_{W(\dbF_{p^n})/\dbZ_p}T_J^0)_{W(k)}$, cf. Corollary 2.3.2. As $T_J^0$ is a product of non-split rank 1 tori, the product of the cocharacters of the orbit under $\sigma$ of each cocharacter of $\tilde T^0/Z(\tilde C^{\der})\cap\tilde T^0$, is the trivial cocharacter. Thus the Newton cocharacter $\tilde\nu_{\tilde w}$ is trivial. Thus $(M,\tilde w\sigma\tilde\mu({1\over p}),G)$ is also basic. But we can write $h^{-1}\tilde w\sigma\tilde\mu({1\over p})h=w\sigma\mu({1\over p})=w\phi$, where $w\in G(W(k))$ is an element that normalizes $\Lie(T)$. We have $w\in N(W(k))$ and moreover $\scrC_w$ is basic.

\medskip
We come back to the case when $k$ is arbitrary and $G$ is split. We prove Theorem 1.3.3 for $\scrC_g$ basic. Let $w\in N(W(k))$ be such that $\scrC_w$ is basic (cf. the six Cases). We have $\nu_g=\nu_w$, cf. Corollary 2.6.4. Thus $\scrC_w\otimes_k \bar k$ and $\scrC_g\otimes_k \bar k$ are rational inner isomorphic, cf. Proposition 2.7. Thus Theorem 1.3.3 holds if $\scrC_g$ is basis. To prove this we did not use that $b_L=1$ and moreover Case 1 works also for the $E_8$, $F_4$, and $G_2$ Lie types. Thus for each $p$-divisible object $\tilde\scrC=(\tilde M,\tilde\phi,\tilde G)$ with a group over a perfect field $\tilde k$ such that $\tilde G$ is a split, reductive group scheme, there exists $\tilde g\in\tilde G(W(\tilde k))$ such that $(\tilde M,\tilde g\tilde\phi,\tilde G)$ is basic.

\medskip\noindent
{\bf 4.2.3. Step 2.} We assume $\scrC_g$ is not basic.
Let $U_G^+(g\phi)$ be the unipotent radical of $P_G^+(g\phi)$. Let $\tilde T$ be a split, maximal torus of $P_G^+(g\phi)$. Up to an inner isomorphism we can assume $\tilde T$ is also a maximal torus of $P$, cf. Subsection 3.1. Two split, maximal tori of $P$ are $P(W(k))$-conjugate. Thus up to an inner isomorphism defined by an element of $P(W(k))$ we can assume $\tilde T=T\leqslant P^+(g\phi)$. Let $L^+$ be the Levi subgroup scheme of $P_G^+(g\phi)$ that contains $T$. We write $g\phi=u^+\phi_1$, where $u^+\in U_G^+(g\phi)(W(k))$ and $\grC:=(M,\phi_1,L^+)$ is a basic Shimura $p$-divisible object (cf. Subsubsection 3.1.1 and Theorem 3.1.2 (b)). Let $\tilde g_0\in L^+(W(k))$ be such that $\tilde g_0\phi_1$ normalizes $\Lie(T)$ and takes the Lie algebra of a Borel subgroup scheme of $L^+$ contained in $L^+\cap P$ to itself (see Subsection 2.5). Let $\tilde w_0:=\tilde g_0(u^+)^{-1}g\in G(W(k))$. We have $\tilde w_0\phi=\tilde g_0\phi_1$ and therefore $\tilde w_0$ normalizes $\Lie(T)$. Thus $\tilde w_0\in N(W(k))$. From the end of Subsubsection 4.2.2 we get the existence of $w^+\in (N\cap L^+)(W(k))$ such that if $w:=w^+\tilde w_0\in N(W(k))$, then $\grE:=(M,w^+\tilde g_0\phi_1,L^+)=(M,w\phi,L^+)$ is basic. As $\grC$ is also basic and $w^+\tilde g_0\in L^+(W(k))$, the Newton quasi-cocharacters of $\grE$ and $\grC$ coincide (cf. Corollary 2.6.4). From this and Theorem 3.1.2 (c) we get that $\nu_w$ and $\nu_g$ are $G(B(k))$-conjugate. Suppose that $k=\bar k$. We know that $\scrC_{(u^+)^{-1}g}$ and $\scrC_w$ (resp. and $\scrC_g$) are rational inner isomorphic, cf. Proposition 2.7 applied to $\grC$ and $\grE$ (resp. cf. Theorem 3.1.2 (c)). Thus $\scrC_w$ and $\scrC_g$ are rational inner isomorphic. This ends the proof of Theorem 1.3.3.\endproof

\medskip\smallskip\noindent
{\bf 4.3. Corollary.} {\it We assume $b_L\Le 1$ and we use the notations of Subsubsection 1.2.4 for $g\in G(W(k))$. Up to a rational inner isomorphism and up to an extension to a finite field extension of $k$, we can assume that we have a direct sum decomposition $M:=\oplus_{\gamma\in \grS(g)} M\cap M_{\gamma}(g)$ and that $\scrC_g\otimes_k \bar k$ is the extension to $\bar k$ of a Shimura $p$-divisible object over a finite field.}

\medskip
\proof We can assume that $G$ is split and (cf. Subsection 2.5) that $\phi$ normalizes the Lie algebra of a split, maximal torus $T$ of $G$ through which $\mu$ factors. Let $w\in N(W(k))$ and $h\in G(B(\bar k))$ be as in Theorem 1.3.3 (b). The $W(\bar k)$-lattice $h(M\otimes_{W(k)} W(\bar k))$ of $M\otimes_{W(k)} B(\bar k)$ is defined over the strict henselization of $W(k)$ and therefore also over $W(k_1)$ with $k_1$ a finite field extension of $k$. Let $M_1$ be the $W(k_1)$-lattice of $M\otimes_{W(k)} B(k_1)$ such that we have $h(M\otimes_{W(k)} W(\bar k))=M_1\otimes_{W(k_1)} W(\bar k)$. Let $G_1$ be the Zariski closure of $G_{B(k_1)}$ in $\GL_{M_1}$. For $h_1\in hG(W(\bar k))\cap G(B(k_1))$ we have $h_1(M\otimes_{W(k)} W(k_1))=M_1$. We choose $h_1$ such that $h^{-1}h_1$ is congruent to $1_{M\otimes_{W(k)} W(\bar k)}$ modulo a high power of $p$; this implies that $h_1\mu_{W(k_1)} h_1^{-1}$ is a Hodge cocharacter of $(M_1,g\phi\otimes\sigma_{k_1},G_1)$. Thus we can write $h_1^{-1}(g\phi\otimes\sigma_{k_1})h_1=h_2(\phi\otimes\sigma_{k_1})$, with $h_2\in G(W(k_1))$. The element $h_1$ is a rational inner isomorphism between $\scrC_{h_2}:=(M\otimes_{W(k)} W(k_1),h_2(\phi\otimes\sigma_{k_1}),G_{W(k_1)})$ and $\scrC_g\otimes_k k_1$. Moreover $h_1^{-1}h\in G(W(\bar k))$ is an inner isomorphism between $\scrC_w\otimes_k \bar k$ and $\scrC_{h_2}\otimes_{k_1} \bar k$. Thus to prove the Corollary, we can assume $g=w$. As $w\phi$ normalizes $\Lie(T)$, the Newton cocharacter of $\scrC_w$ factors through $T_{B(k)}$ and thus it extends to a cocharacter of $T$. Therefore we have a direct sum decomposition $M:=\oplus_{\gamma\in \grS(w)} M\cap M_{\gamma}(w)$ normalized by $T$ and thus also by $\mu$. To check the last part of the Corollary we can assume $k=\bar k$.

Let $M_{\dbZ_p}^w:=\{x\in M|w\phi\mu(p)(x)=x\}$. As in the proof of Claim 2.2.2 we argue that $T$ and $G$ are the extensions to $W(k)$ of a torus $T_{\dbZ_p}^w$ of $\GL_{M_{\dbZ_p}^w}$ and respectively of a reductive subgroup scheme of $\GL_{M_{\dbZ_p}^w}$ isomorphic to $G_{\dbZ_p}$ (see Subsection 2.6) and accordingly denoted also by $G_{\dbZ_p}$; we have  $T_{\dbZ_p}^w\leqslant G_{\dbZ_p}$ as this holds after pullback to $\Spec(W(k))$. Let $k_0$ be a finite field such that the cocharacter $\mu$ of $T^w$ is the pullback of a cocharacter $\mu_0$ of $T^w_{W(k_0)}$. We conclude that $\scrC_w$ is the extension to $k=\bar k$ of the Shimura $p$-divisible object $(M_{\dbZ_p}^w\otimes_{\dbZ_p} W(k_0),(1\otimes\sigma_{k_0})\mu_0({1\over p}),G_{W(k_0)})$ over $k_0$. \endproof

\medskip\smallskip\noindent
{\bf 4.4. Corollary.} {\it Suppose that $k=\bar k$ and $b_L\Le 1$. Let $o(W_G)$ be the order of the Weyl group $W_G$ of $G$. Let $k_1$ be an algebraically closed field that contains $k$. Let $\scrR$ be the equivalence relation on the set $\{(M\otimes_{W(k)} W(k_1),g(\phi\otimes\sigma_{k_1}),G_{W(k_1)})|g\in G(W(k_1))\}$ of Shimura $p$-divisible objects, defined by rational inner isomorphisms. Let $\grR$ be the cardinal of the set of equivalence classes of $\scrR$. Then $\grR\in\{1,\ldots,o(W_G)\}$ and it does not depend on $k_1$.} 

\medskip
\proof We can assume we are in the context of Theorem 1.3.3. The estimate $\grR\le o(W_G)$ follows from Theorem 1.3.3 (b) and Fact 4.2.1. All $\scrC_w$'s with $w\in N(W(k))$ are definable over finite fields (cf. end of Corollary 4.3) and each rational inner isomorphism between extensions to $k_1$ of two such $\scrC_w$'s is defined by an element of $G(B(k))=G(B(k_1))\cap \GL_M(B(k))$ (cf. [37, Lemma 3.9]). Thus the number $\grR$ is independent of $k_1$. \endproof

\medskip\smallskip\noindent
{\bf 4.5. On $\grR$.} Until Section 5 we assume $k=\bar k$. The estimate $\grR\le o(W_G)$ is gross. Often one can compute $\grR$ explicitly in two steps as follows. Let $T$, $\mu$, $B$, and $g_0$ be as in Subsection 2.5. To ease notations we will assume that $g_0=1_M$. Thus $\phi$ normalizes $\Lie(T)$ and takes $\Lie(B)$ to itself. Let $\sigma_0$, $M_{\dbZ_p}$, and $G_{\dbZ_p}$ be obtained as in the proof of Claim 2.2.2. We have $\sigma_0=1_{M_{\dbZ_p}}\otimes\sigma$ and $\phi=\sigma_0\mu(({1\over p})=(1_{M_{\dbZ_p}}\otimes\sigma)\mu({1\over p})$. Let $T^0$ (resp. $B^0$) be the maximal torus of $G_{\dbZ_p}$ whose extension to $W(k)$ is $T$ (resp. is $B$). As $\Lie(T)$ (resp. $\Lie(B)$) is normalized by $\sigma_0=\phi\mu(p)$, the existence of $T$ (resp. $B$) is argued in the same way we argued the existence of $G_{\dbZ_p}$ in the proof of Claim 2.2.2 (and it is also implied by [48, Prop. 3.2]). Let $N^0$ be the normalizer of $T^0$ in $G_{\dbZ_p}$. Let $m\in\dbN$ be the smallest number such that the torus $T^0_{W(\dbF_{p^m})}$ is split. Let $k_0:=\dbF_{p^m}$. We identify $W_G=(N^0/T^0)(W(k_0))$. For $w\in W_G$, let $n_w\in N^0(W(k_0))$ be an element that represents $w$. Let $\mu_0:\dbG_m\to T^0_{W(k_0)}$ be the unique cocharacter such that we have $\mu_{0W(k)}=\mu$. Let 
$$\scrC^0_{n_w}:=(M_{\dbZ_p}\otimes_{\dbZ_p} W(k_0),n_w(1_{M_{\dbZ_p}}\otimes\sigma_{k_0})\mu_0({1\over p}),G_{W(k_0)}).$$ 
\noindent
Let $o_w:\dbG_m\to T^0_{W(k_0)}$ be the inverse of the product of the cocharacters of $T^0_{W(k_0)}$ that form the orbit of $\mu_0$ under powers of $n_w\sigma_{k_0}$. Let $d_w\in\dbN$ be the number of elements of this orbit. The Newton quasi-cocharacter 
$$\nu_w^0\in \Xi(G_{B(k_0)})$$ 
of $\scrC^0_{n_w}$ is the equivalence class of the pair $(o_w,{1\over d_w})\in \Lambda(G_{B(k_0)})$. Due to Fact 4.2.1, $\nu_w^0$ does not depend on the choice of $n_w$ and this justifies our notation. Let $\nu(w)$ be the quasi-cocharacter of $T^0_{B(k_0)}$ that is $N^0(W(k_0))$-conjugate to $\nu_w^0$ and such that its action on $\Lie(B^0_{B(k_0)})$ is via non-negative rational powers of $p$. Let $L(w)$ be the $\dbQ_p$-form of the Levi subgroup of $\scrC^0_{n_w}\otimes_{k_0} k$ whose Lie algebra is $\{x\in\Lie(G_{B(k)})|n_w(1_{M_{\dbZ_p}}\otimes\sigma)\mu_0({1\over p})(x)=x\}$.

The first step is to list all distinct $\nu(w)$'s. Let $\grP$ be their number. Let $w_1$, $w_2\in W_G$. If $\scrC_{n_{w_1}}=\scrC_{n_{w_1}}^0\otimes_{k_0} k$ and $\scrC_{n_{w_2}}=\scrC_{n_{w_2}}^0\otimes_{k_0} k$ are rational inner isomorphic, then $\nu_{n_{w_1}}$ and $\nu_{n_{w_2}}$ are $G(B(k))$-conjugate. Thus $\nu_{w_1}^0$ and $\nu_{w_2}^0$ are $G_{B(k_0)}(B(k_0))$-conjugate and therefore we have $\nu(w_1)=\nu(w_2)$. This implies  
$$\grP\le \grR.\leqno (9)$$ 
\indent
Suppose that $\nu(w_1)=\nu(w_2)$. This implies that the two parabolic subgroups $\break P_{G_{W(k_0)}}^+(n_{w_1}\sigma_{k_0}\mu_0({1\over p}))_{B(k_0)}$ and $P_{G_{W(k_0)}}^+(n_{w_2}\sigma_{k_0}\mu_0({1\over p}))_{B(k_0)}$ of $G_{B(k_0)}$ are $G_{W(k_0)}(B(k_0))$-conjugate. As both parabolic subgroup schemes $P_{G_{W(k_0)}}^+(n_{w_1}\sigma_{k_0}\mu_0({1\over p}))$ and $P_{G_{W(k_0)}}^+(n_{w_2}\sigma_{k_0}\mu_0({1\over p}))$ of $G_{W(k_0)}$ contain $T_{W(k_0)}$, there exists $w_{12}\in N^0(W(k_0))$ such that we have an identity $w_{12}P_{G_{W(k_0)}}^+(n_{w_1}\sigma_{k_0}\mu_0({1\over p}))w_{12}^{-1}=P_{G_{W(k_0)}}^+(n_{w_2}\sigma_{k_0}\mu_0({1\over p}))$. In general, $w_{12}$ does not centralize $\mu_0$ and therefore in general the inequality (9) is strict (cf. also Example 4.6.4 below). 

Thus the second step is to decide for which pairs $(w_1,w_2)\in W_G\times W_G$ with the property that $\nu(w_1)=\nu(w_2)$, there exists a rational inner isomorphism $h_{12}$ between $\scrC_{n_{w_1}}$ and $\scrC_{n_{w_2}}$. The existence of $h_{12}$ is a problem of deciding if the torsor of $L(w_1)$ that defines $L(w_2)$ is trivial or not (see [23, Prop. 6.3] and [37, Prop. 1.17]); as the group $W_G$ is finite, one can compute all classes that define such torsors of $L(w_1)$.

\medskip\noindent
{\bf 4.5.1. Remark.} We use the language of [23] and [37]. Let $B(G_{B(k)})$ be the set of $\sigma$-conjugacy classes of $G(B(k))$. Let $\scrU_+:=\bar C^{\Gal(B(k_0)/\dbQ_p)}\cap (X_*(T^0_{B(k_0)})\otimes_{\dbZ}\dbQ)$, where $\bar C$ is the closed Weyl chamber of $X_*(T^0_{B(k_0)})\otimes_{\dbZ} \dbQ$ that corresponds to $B^0$. The set $\scrU_+$ has a natural partial ordering. Let $\scrW_G:W_G\to B(G_{B(k)})$ be the map  that takes $w\in W_G$ to the $\sigma$-conjugacy class of $n_w\sigma_{k_0}(\mu_0)({1\over p})\in G(B(k))$ (to be compared with Remark 1.2.3). Let $\scrN_G:B(G_{B(k)})\to \scrU_+$ be the Newton map. Let $\kappa_G:B(G_{B(k)})\to \pi_1(G_{B(k)})_{\sigma}$ be the Kottwitz map. We identify naturally $\nu(w)$ with $\scrN_G(\scrW_G(w))\in\scrU_+$. Let $-\bar\mu\in\scrU_+$ be the average of the orbit of $-\mu$ (or $-\sigma(\mu)$) under $\Gal(B(k_0)/\dbQ_p)$. From [37, Thm. 4.2] we get 
$$\im(\scrW_G)\subseteq \{x\in B(G_{B(k)})|\kappa_G(x)=\kappa_G(\sigma_{k_0}(\mu_0)({1\over p}))\}\cap \scrN_G^{-1}(\{x\in\scrU_+|x\Le -\bar\mu\}).\leqno (10)$$
One can combine Theorem 1.3.3 with [55] to show that the two sets of (10) are in fact equal. In general it is hard to compute the two sets of the intersection of (10). If $w\in W_G$, then in general the set $\scrN_G^{-1}(\nu(w))$ is not included in $\scrW_G(W_G)$ (see Example 4.6.3 below). The number of elements of $\im(\scrW_G)$ (resp. of $\im(\scrN_G\circ\scrW_G)$) is $\grR$ (resp. is $\grP$).

\medskip\smallskip\noindent
{\bf 4.6. Examples.} Suppose that $b_L=1$ and $G_{\dbZ_p}=\Res_{W(k_0)/\dbZ_p} J$, where $J$ is an absolutely simple adjoint group scheme over $W(k_0)$ of $B_r$ Lie type. Thus $J$ is split and it is the $\pmb{\text{SO}}$ group scheme of the following quadratic form $Q=x_1x_2+\cdots+x_{2r-1}x_{2r}+x_{2r+1}^2$ on $V_0:=W(k_0)^{2r+1}$. Therefore the number $m$ (of Subsection 4.5) is the number $n$ (of Subsubsection 4.2.2). As $\dbZ_p$-modules, we can identify $M_{\dbZ_p}=V_0$. We write 
$$M_{\dbZ_p}\otimes_{\dbZ_p} W(k_0)=V_0\otimes_{\dbZ_p} W(k_0)=\oplus_{s=1}^n V_s,$$ 
where each $V_s$ is a free $W(k_0)$-module of rank $2r+1$ normalized by $G_{W(k_0)}$. Let $l\in\dbN$ be such that the cocharacter $\mu_0$ acts non-trivially on precisely $l$ of the $V_s$'s. Let $B_Q$ be the symmetric bilinear form on $V_0$ defined by the rule $B_Q(x,y):=Q(x+y)-Q(x)-Q(y)$, where $Q$ is identified naturally with a quadratic function $V_0\to W(k_0)$ and where $x$, $y\in V_0$. We recall (cf. Subsection 2.1) that the scalar extensions of $B_Q$ are also denoted by $B_Q$. 

\medskip\noindent
{\bf 4.6.1. Proposition.} {\it Under the assumptions of the previous paragraph we have:

\medskip
{\bf (a)} If $l=1$, then $\grP=r+1$.

\smallskip
{\bf (b)} If $l=2$, then $\grP=r+(r-[{r\over 2}])([{r\over 2}]+1)$.}
\medskip

\proof For $s\in S(1,n)$ we choose a $W(k_0)$-basis $\scrB_s=\{e_1^s,e_2^s,\ldots,e_{2r+1}^s\}$ for $V_s$ such that $T^0_{W(k_0)}$ normalizes $W(k_0)e$ for all $e\in\scrB_s$. Let $\scrS:=\cup_{s=1}^n \scrB_s$. We can assume that $\scrS$ is such that each $B_Q(e_i^s,e_j^s)$ is $0$ or $1$ depending on the fact that $\{(i,j),(j,i)\}$ has or has not a trivial intersection with the set $\{(1,2),(3,4),\ldots,(2r-1,2r)\}$ and that $\mu_0(p)$ acts:

\medskip
-- trivially on $e_i^s$, if either $i\ge 3$ or $s>l$, and

\smallskip
-- as the multiplication by $p^{i}$ on $e_{{3-i}\over 2}^s$, if $s\in S(1,l)$ and $i\in\{-1,1\}$.

\medskip
We can also assume that $\scrS$ is such that there exists an $n$-th cycle $\pi$ which is a permutation of $S(1,n)$  with the property that the $\sigma_{k_0}$-linear automorphism of $M_{\dbZ_p}\otimes_{\dbZ_p} W(k_0)=V_0\otimes_{\dbZ_p} W(k_0)$ (which fixes $M_{\dbZ_p}=V_0$) takes each $e_i^s$ to $e_i^{\pi(s)}$. From [5, plate II] we get that $W_G$ as a set is in one-to-one correspondence with the set of those elements 
$$h\in \GL_{\oplus_{s=1}^n V_s}(W(k_0))=\GL_{M_{\dbZ_p}}(W(k_0))$$ 
which for each $s\in S(1,n)$, take $e^s_{2r+1}$ to $e^{\pi(s)}_{2r+1}$ and take each pair $(e_{2q-1}^s,e_{2q}^s)$ with $q\in S(1,r)$ to either $(e_{2q_1-1}^{\pi(s)},e_{2q_1}^{\pi(s)})$ or $(e_{2q_1}^{\pi(s)},e_{2q_1-1}^{\pi(s)})$ for some $q_1\in S(1,r)$. For each such $h$ we have $h(\scrS)=\scrS$ and there exists a unique $w\in W_G$ and a unique representative $n_w\in N_0(W(k_0))$ of $w$ such that the actions of $h$ and $n_w(1_{V_0}\otimes\sigma_{k_0})$ on $\scrS$ are the same.

Two cocharacters $\mu_1$ and $\mu_2$ of $G_{B(k_0)}$ are $G(B(k_0))$-conjugate if and only if for each $s\in S(1,n)$ the formal characters of the action of $\dbG_m$ on $V_s$ via $\mu_j$ does not depend on $j\in\{1,2\}$. Thus, as $\sigma_{k_0}$ permutes transitively the simple factors of $G_{W(k_0)}$, the number $\grP$ is the number of distinct Newton polygons of $F$-crystals of the form 
$$\scrD_{w}:=(V_0\otimes_{\dbZ_p} B(k_0),n_w(1_{V_0}\otimes\sigma_{k_0})\mu_0({1\over p})),$$ 
where $h$ and $w$ vary as described. Due to the existence of $B_Q$, the Newton polygon $N(w)$ of $\scrD_{w}$ is $0$-symmetric i.e., its multiplicity of a slope $\alpha$ is the same as its multiplicity of the slope $-\alpha$. Thus to describe $N(w)$'s it is enough to list (with multiplicities) their positive Newton polygon slopes. Below by orbits we mean the orbits of the action of $h$ on $\scrS$.

Suppose that $l=1$. If the orbit of $e^1_1$ contains $e^1_2$, then the Newton polygon $N(w)$ has only one slope $0$. If the orbit of $e^1_1$ does not contain $e^1_2$ and has $c\in S(1,r)$ elements, then the Newton polygon $N(w)$ has only one positive slope ${1\over {nc}}$ with multiplicity $nc$. Thus $\grP=r+1$ i.e., (a) holds.

If $l=2$, we split the computation of $\grP$ as follows. We count:

\medskip
{\bf (i)} $1$ for the Newton polygon that has all slopes $0$;

\smallskip
{\bf (ii)} $r$ for Newton polygons that have only one positive slope $2\over {nc}$ with multiplicity $nc$, $c\in S(1,r)$; they correspond to the case when $e_1^1$ and $e_1^2$ are in the same orbit but this orbit does not contain $e_2^1$ (and thus does not contain $e_2^2$);

\smallskip
{\bf (iii)} $r-1$ for Newton polygons that have only one positive slope $1\over {nc}$ with multiplicity $nc$, $c\in S(1,r-1)$; they correspond to the case when there exists $i\in\{1,2\}$ such that the orbit of $e_1^i$ does not contain $e_2^i$, $e_1^{3-i}$, or $e_2^{3i-1}$ but $e_1^{3-i}$ and $e_2^{3-i}$ belong to the same orbit;

\smallskip
{\bf (iv)} $1$ for each pair $(c,d)\in S(1,r)\times S(1,r)$, with $c+d\le r$ and $c<d$. 

\medskip
For each pair as in (iv) we get a Newton polygon that has two positive slopes ${1\over {nc}}$ and ${1\over {nd}}$ with multiplicities $nc$ and respectively $nd$; it corresponds to the cases when the orbit of $e_1^j$ under $h$ does not contain $e^j_2$, $e_1^{3-j}$, or $e^{3-j}_2$ but has $v_j$ elements, $j=\overline{1,2}$, with $(v_1,v_2)\in\{(c,d),(d,c)\}$. We note that the case $c=d$ gives rise to Newton polygons that were already counted by (ii). Thus the number of distinct Newton polygons is $1+r+(r-1)+\sum_{c=1}^{[{r\over 2}]} (r-2c)=r([{r\over 2}]+2)-[{r\over 2}]([{r\over 2}]+1)=r+(r-[{r\over 2}])([{r\over 2}]+1)$. Thus  (b) holds.\endproof

\medskip\noindent
{\bf 4.6.2. Example.} We assume that $l\in\{1,2\}$ and we check that $\grP=\grR$. The set $\grS(w)$ of Newton polygon slopes of $\scrD_{w}$ has 1, 3, or 5 elements, cf. proof of Proposition 4.6.1. The Newton polygon slope decomposition $V_0\otimes_{\dbZ_p} B(k_0)=\oplus_{j\in \grS(w)} M_j(w)$ of $\scrD_{w}$ is such that $B_Q(M_j(w),M_u(w)))=0$ if $j+u\neq 0$ and we have $M_j(w)=\oplus_{s=1}^n V_s[{1\over p}]\cap M_j(w)$; here $j$, $u\in \grS(w)$. We have $\dim_{B(k_0)}(V_s[{1\over p}]\cap M_j(w))=\dim_{B(k_0)}(V_s[{1\over p}]\cap M_{-j}(w))$. Thus the number $r_0(w):=\dim_{B(k_0)}(V_s[{1\over p}]\cap M_0(w))$ is odd. Let $B_{Q_s}(w)$ be the symmetric bilinear form on $V_s[{1\over p}]\cap M_0(w)$ induced by $B_Q$. Let $\grI_s(w)$ be the isomorphism class of the triple $(V_s[{1\over p}]\cap M_0(w),(n_w(1_{V_0}\otimes\sigma_{k_0})\mu_0({1\over p}))^n,B_{Q_s}(w))\otimes_{k_0} k$. Maps defined by the rule $e_i^s\to e_i^{\pi(s)}$ allow us to identify $\grI_s(w)$ with $\grI_{\pi(s)}(w)$. Thus $\grI_s(w)$ does not depend on $s\in S(1,n)$ and therefore we denote it by $\grI(w)$. If $l=2$, then $\grI(w)$ is not determined in general by $r_0(w)$. 

To check that $\grP=\grR$, it suffices to show that the isomorphism class of the $\sigma^n$-$F$-isocrystal with a bilinear form $(M_j(w)+M_{-j}(w),(n_w(1_{V_0}\otimes\sigma_{k_0})\mu_0({1\over p}))^n,B_{Q_{\pm j}})\otimes_{k_0} k$ is uniquely determined by $j\in \grS(w)$ and by $\dim_{B(k_0)}(M_j(w))$, where $B_{Q_{\pm j}}$ is the restriction of $B_Q$ to $M_j(w)+M_{-j}(w)$. If $j\in \grS(w)\cap (0,\infty)$, then this holds due to Dieudonn\'e--Manin's classification of $\sigma^n$-$F$-isocrystals over $k$ and the fact that $B_{Q_{\pm j}}(M_i(w),M_i(w))=0$ if $i\in\{j,-j\}$. Thus the equality $\grP=\grR$ holds if and only if always the Newton polygon $N(w)$ determines $\grI(w)$. To check that $N(w)$ determines $\grI(w)$, we can assume $r_0(w)\le 2r-1$ (cf. Proposition 2.7). The case $l=1$ is a consequence of the fact that each absolutely simple adjoint group scheme of $B_{{r_0(w)-1}\over 2}$ Dynkin type over $W(k_0)$ is split. 

Next we consider the case when $l=2$. The fact that $\grI(w)$ is determined by $N(w)$ follows easily from the descriptions of items 4.6.1 (ii) to (iv). The only ``ambiguity" of the type of $h$ producing Newton polygons as in items 4.6.1 (ii) to (iv) is between item 4.6.1 (ii) and the variant of item 4.6.1 (iv) with $c=d$ (as $\grI_s(w)$ does not depend on $s$, it is irrelevant if in item 4.6.1 (iii) we have $i=1$ or $i=2$). But if $h$ is as in item 4.6.1 (ii) or as in the mentioned variant, then for $i,s\in\{1,2\}$ we have $e_i^s\notin M_0(w)$ and thus the fact that $N(w)$ determines $\grI(w)$ is argued as for $l=1$. 

\medskip\noindent
{\bf 4.6.3. Example.} We refer to Proposition 4.6.1 and Example 4.6.2 with $n=l=1$. The reductive group $L(w)$ of Subsection 4.5 is a product $L^0(w)\times_{\dbQ_p} L^+(w)$, where $L^0(w)$ is an absolutely simple group of $B_u$ Dynkin type with $u\in S(0,r)$ and where $L^+(w)$ is trivial if $u=r$ and is the group of invertible elements of a semisimple $\dbQ_p$-algebra if $u<r$. If $u<r$, then the group $L^+(w)$ corresponds to the unique Newton polygon slope in $\grS(w)\cap (0,\infty)$. If $u>0$, then we have $\Lie(L^0(w))=\Lie(L(w))\cap \End(M_0(w)\otimes_{B(k_0)} B(k))$. We have $u=0$ if and only if $N(w)$ has slope ${1\over r}$. The set $H^1(\dbQ_p,L(w))=H^1(\dbQ_p,L^+(w))\times H^1(\dbQ_p,L^0(w))$ has $1$ or $2$ classes depending on the fact that $u$ is or is not $0$. Thus from [37, Props. 1.15 and  1.17] we get that if $N(w)$ has (resp. does not have) slope ${1\over r}$, then $\scrN_G^{-1}(\nu(w))$ has $1$ element (resp. $2$ elements). Thus the set $\cup_{w\in W_G} \scrN_G^{-1}(\nu(w))$ has $2r+1$ elements, while $\scrW_G(W_G)$ has $\grP=r+1$ elements.  

\medskip\noindent
{\bf 4.6.4. Example.} Suppose that $r\in\dbN$, that $m=2n$, and that $G_{\dbZ_p}=\Res_{W(\dbF_{p^{n}})/\dbZ_p} J$, where $J$ is the $\pmb{\text{SO}}$ group scheme of the quadratic form $x_1x_2+x_3x_4+x_5x_6+\cdots+x_{2r-3}x_{2r-2}+yx_{2r-1}^2+ux_{2r-1}x_{2r}+zx_{2r}^2$ on $V_0:=W(\dbF_{p^{n}})^{2r}$ for $u,y,z\in W(\dbF_{p^{n}})$ as in Case 6 of Subsubsection 4.2.2. By repeating the constructions of Proposition 4.6.1 and Example 4.6.2 with the appropriate restrictions on $h$ (cf. the structure of the Weyl group in [5, plate IV]), the identity component of the group of automorphisms of $(V_s[{1\over p}]\cap M_0(w),(n_w(1_{V_0}\otimes\sigma_{k_0})\mu_0({1\over p}))^n,B_{Q_s}(w))\otimes_{k_0} k$ is not determined by the Newton polygon $N(w)$. 

For instance, if $n=1$, $r\ge 3$, and the orbit of $e_1^1$ under $h$ does (resp. does not) contain $e_{2r-1}^1$ and also does not contain $e_{2r}^1$, then this group is split (resp. is non-split); thus we have $\grR>\grP$. As $n=1$, we have $m=2$, $e_{2r-1}^1,e_{2r}^1\in V_0\otimes_{\dbZ_p} W(\dbF_{p^2})$ correspond naturally to the quadratic form $yx_{2r-1}^2+ux_{2r-1}x_{2r}+zx_{2r}^2$, and $h\in \GL_{V_0}(W(\dbF_{p^2}))=\GL_{M_{\dbZ_p}}(W(\dbF_{p^2}))$. 

\medskip\smallskip\noindent
{\bf 4.7. Quasi-polarizations.}
Suppose that there exists a perfect bilinear form $\lambda_M$ on $M$ normalized by $G$ and such that we have $\lambda_M(\phi(x),\phi(y))=p^{a+b}\sigma(\lambda_M(x,y))$ for all $x,y\in M$ (see Subsection 2.1 for $a$ and $b$). We refer to $\lambda_M$ as a principal bilinear quasi-polarization of $\scrC_{1_M}$. Let $G^0$ be the Zariski closure in $G$ of the identity component of the subgroup of $G_{B(k)}$ that fixes $\lambda_M$. Thus the quotient group scheme $G/G^0$ is either $\dbG_m$ or trivial. The group scheme $G^0$ is the extension of $G^{\der}$ by a subtorus of $G^{\ab}$ and thus it is a reductive group scheme over $W(k)$. We situate ourselves in the context of Theorem 1.3.3 (b). Thus $b_L\Le 1$, $k=\bar k$, and $\phi$ normalizes $\Lie(T)$. Let $g\in G^0(W(k))$. We check that in Theorem 1.3.3 we can choose $w$ and $h$ to fix $\lambda_M$. For this, we can assume $G^0\neq G$.

In Subsubsection 3.1.1 we can always choose $g_3$, $g_2$, and $g_1\in G^0(W(k))$ and in the proof of Theorem 3.1.2 we have $h_i\in G^0(B(k))$ for all $i\in S(1,q)$. Thus to check that we can choose $w$ and $h$ to fix $\lambda_M$, it suffices to consider the case when $\scrC_g$ is basic. In Subsubsection 4.2.2, we can always choose $w\in (G^{\der}\cap N)(W(k))$. Thus we only have to show that if in Proposition 2.7 we have $g_1$, $g_2\in G^0(W(k))$, then there exists $h\in G^0(B(k))$ such that we have $hg_1\phi=g_2\phi h$. If $G^\flat$ is as in Subsubsection 2.6.5, then $G^{\flat 0}:=G^\flat\times_G G^0$ is a reductive group scheme and we have $g\in\im(G^{\flat 0}(W(k))\to G^0(W(k)))$. Thus as in the proof of Proposition 2.7 we can assume that $G^{\der}$ is simply connected. The subtorus $G^0/G^{\der}$ of $G^{\ab}$ is the pullback of a subtorus of $G^{\ab}_{\dbZ_p}$ and $G^0(W(k))$ surjects onto $(G^0/G^{\der})(W(k))$. Based on these facts, an argument similar to the one that proved Fact 2.6.2 shows that up to inner isomorphisms defined by elements of $G^0(W(k))$ we can assume that $g_1$, $g_2\in G^{\der}(W(k))$. Therefore the image of $h$ in $G^{\ab}(B(k))$ belongs to $G^{\ab}_{1\dbQ_p}(\dbQ_p)$, where $G_{1\dbQ_p}$ is as in the proof of Proposition 2.7. But as $H^1(\dbQ_p,G^{\der}_{\dbQ_p})$ has only the trivial class (see [21]), $G_{1\dbQ_p}(\dbQ_p)$ surjects onto $G^{\ab}_{1\dbQ_p}(\dbQ_p)$. Thus by replacing $h$ with a right $G_{1\dbQ_p}(\dbQ_p)$-multiple of it we can assume that $h\in G^{\der}(B(k))\leqslant G^0(B(k))$.

\bigskip\smallskip
\noindent
{\boldsectionfont 5. Three geometric applications}

\bigskip
In this Section we apply Theorem 1.3.3 (b) to geometric contexts which pertain to Shimura varieties of Hodge type. Subsection 5.1 introduces standard Hodge situations. They are concrete ways of getting good moduli spaces of principally polarized abelian varieties endowed with (specializations of) Hodge cycles in mixed characteristic $(0,p)$. In Subsubsections 5.1.4 to 5.1.9 we list properties of standard Hodge situations we will require in Subsections 5.2 to 5.4. In Subsection 5.2 we formulate Manin problems for standard Hodge situations. The Main Theorem 5.2.3 solves these problems if $p\Ge 3$ or if $p=2$ and two mild conditions hold. In Subsection 5.3 we study rational stratifications of special fibres of the mentioned moduli spaces. In Subsection 5.4 we formulate integral Manin problems. Theorem 5.4.2 solves them for certain standard PEL situations. 

We use the terminology of [9] for Hodge cycles on an abelian scheme $A_Z$ over a reduced $\dbQ$--scheme $Z$. Thus each Hodge cycle $v$ of $A_Z$ has a de Rham component $v_{\text{dR}}$ and an \'etale component $v_{\acute et}$. The \'etale component $v_{\acute et}$ as its turn has an $l$-component $v_{\acute et}^l$, for each rational prime $l$. For instance, if $Z$ is the spectrum of a field $E$, then $v_{\acute et}^p$ is a suitable $\Gal(E)$-invariant tensor of the tensor algebra of $H^1_{\acute et}(A_{\bar Z},\dbQ_p)\oplus (H^1_{\acute et}(A_{\bar Z},\dbQ_p))^*\oplus \dbQ_p(1)$, where $\bar Z:=\Spec(\bar E)$ and where $\dbQ_p(1)$ is the usual Tate twist. If $E$ is a subfield of $\dbC$ we also use the Betti realization of $v$: it corresponds to $v_{\text{dR}}$ (resp. $v_{\acute et}^l$) via the standard isomorphism that relates the Betti cohomology with $\dbQ$--coefficients of $A_Z\times_Z \Spec(\dbC)$ with the de Rham (resp. $\dbQ_l$-\'etale) cohomology of $A_{\bar Z}$ (see [9]). 

Let $\dbS:=\Res_{\dbC/\dbR} \dbG_m$. A {\it Shimura pair} $(G_{\dbQ},X)$ consists of a reductive group $G_{\dbQ}$ over $\dbQ$ and a $G_{\dbQ}(\dbR)$-conjugacy class $X$ of homomorphisms $\dbS\to G_{\dbR}$ which satisfy Deligne's axioms of [8, Subsect. 2.1]: the Hodge $\dbQ$--structure on $\Lie(G_{\dbQ})$ defined by $x\in X$ is of type $\{(-1,1),(0,0),(1,-1)\}$, the composite of $x(i)$ with the adjoint representation $\theta:G_{\dbR}\to \GL_{\Lie(G^{\ad}_{\dbR})}$ is a Cartan involution of $\Lie(G^{\ad}_{\dbR})$, and no simple factor of $G^{\ad}_{\dbQ}$ becomes compact over $\dbR$. Thus $X$ has a canonical structure of a hermitian symmetric domain, cf. [8, Cor. 1.1.17]. For basic facts on Shimura pairs, on their reflex fields, types, and canonical models, and on injective maps between them see [7], [8], [30, Sect. 1], [31], and [45, Subsects. 2.1 to 2.10]. Let $k(\ast)$ be the residue field of a finite prime $\ast$ of a number field.

\medskip\smallskip\noindent
{\bf 5.1. Standard Hodge situations.} Let $(\GSp(W,\psi),S)$ be a Shimura pair that defines a Siegel modular variety. Thus $(W,\psi)$ is a symplectic space over $\dbQ$ and $S$ is the set of all homomorphisms $\dbS\to \GSp(W,\psi)_{\dbR}$ that define Hodge $\dbQ$--structures on $W$ of type $\{(-1,0),(0,-1)\}$ and that have either $2\pi i\psi$ or $-2\pi i\psi$ as polarizations. Let $r:={{\dim_\dbQ(W)}\over 2}$. Let $\dbA_f$ (resp. $\dbA_f^{(p)}$) be the $\dbQ$--algebra of finite ad\`eles (resp. of finite ad\`eles with the $p$-component omitted). We have $\dbA_f=\dbA_f^{(p)}\times\dbQ_p$. We start with an injective map 
$$f:(G_{\dbQ},X)\hookrightarrow (\GSp(W,\psi),S)$$ 
of Shimura pairs. Thus $f:G_{\dbQ}\hookrightarrow \GSp(W,\psi)$ is a monomorphism of reductive groups over $\dbQ$ such that we have $f_{\dbR}\circ x\in S$ for all $x\in X$. 

We have $\dbS(\dbR)=\dbC\setminus\{0\}$. We identify $\dbS(\dbC)=(\dbC\setminus\{0\})\times (\dbC\setminus\{0\})$ in such a way that the monomorphism $\dbS(\dbR)\hookrightarrow\dbS(\dbC)$ induces the map $z\to (z,\bar z)$. Let $x\in X$. Let $\mu_x:\dbG_m\to G_{\dbC}$ be the cocharacter given on complex points by the rule $z\to x_{\dbC}(z,1)$. The reflex field $E(G_{\dbQ},X)$ of $(G_{\dbQ},X)$ is the subfield of $\dbC$ that is the field of definition of the $G_{\dbQ}(\dbC)$-conjugacy class $[\mu_x]$ of (any) $\mu_x$ (see [8] and [30]). It is a number field.
Let $v$ be a prime of $E(G_{\dbQ},X)$ that divides $p$. Let $O_{(v)}$ be the localization of the ring of integers of $E(G_{\dbQ},X)$ at the prime $v$. Let $L$ be a
$\dbZ$-lattice of $W$ such that $\psi$ induces a perfect form 
$\psi:L\otimes_{\dbZ} L\to\dbZ$. Let $\dbZ_{(p)}$ be the localization of $\dbZ$ at the prime $p$. Let $L_{(p)}:=L\otimes_{\dbZ} \dbZ_{(p)}$. 

Until the end we assume that:

\medskip\noindent
{\bf ($\natural$)} the Zariski closure
$G_{\dbZ_{(p)}}$ of $G_{\dbQ}$ in $\GSp(L_{(p)},\psi)$ is a reductive group scheme over $\dbZ_{(p)}$. 

\medskip\noindent
Thus $v$ is unramified over $p$, cf. [31, Cor. 4.7 (a)]. Let $K_p:=\GSp(L_{(p)},\psi)(\dbZ_p)$. It is a hyperspecial subgroup of $\GSp(W\otimes_{\dbQ} \dbQ_p,\psi)(\dbQ_p)$ and the intersection $H:=G_{\dbQ_p}(\dbQ_p)\cap K_p$ is a hyperspecial subgroup of $G_{\dbQ_p}(\dbQ_p)$. We call the triple 
$$(f,L,v)$$  
a {\it potential standard Hodge situation}.

Let $C(G_{\dbQ})$ be the set of compact, open subgroups of $G_{\dbQ}(\dbA_f)$ with the inclusion relation. Let $\Sh(G_{\dbQ},X)$ be the canonical model over $E(G_{\dbQ},X)$ of the complex Shimura variety (see [7, Thm. 4.21 and Cor. 5.7]; see [8, Cor. 2.1.11] for the below identity):
$$\Sh(G_{\dbQ},X)_{\dbC}:=\text{proj.}\text{lim.}_{\tilde H\in C(G_{\dbQ})} G_{\dbQ}(\dbQ)\backslash X\times G_{\dbQ}(\dbA_f)/\tilde H=G_{\dbQ}(\dbQ)\backslash X\times G_{\dbQ}(\dbA_f).\leqno (11)$$ 
From (11) and [31, Prop. 4.11] we get that $\Sh(G_{\dbQ},X)_{\dbC}/H=G_{\dbZ_{(p)}}(\dbZ_{(p)})\backslash X\times G_{\dbQ}(\dbA_f^{(p)})$. From this identity and its analogue for $\Sh(\GSp(W,\psi),S)_{\dbC}/K_p$, we get that $\Sh(G_{\dbQ},X)/H$ is a closed subscheme of $\Sh(\GSp(W,\psi),S)_{E(G_{\dbQ},X)}/K_p$. Let $\scrM$ be the $\dbZ_{(p)}$-scheme which parametrizes isomorphism classes of principally polarized abelian schemes over $\dbZ_{(p)}$-schemes that are of relative dimension $r$ and that have in a compatible way level $s$ symplectic similitude structures for all numbers $s\in (\dbN\setminus p\dbN)$; the group $\GSp(W,\psi)(\dbA^{(p)}_f)$ acts naturally on $\scrM$. These symplectic structures and this action are defined naturally via $(L,\psi)$ (for instance, see [45, Subsect. 4.1]). We have a natural identification $\scrM_{\dbQ}=\Sh(\GSp(W,\psi),S)/K_p$ compatible with the $\GSp(W,\psi)(\dbA^{(p)}_f)$-actions, cf. [7, Ex. 4.16].  Let 
$$\scrN$$ 
be the normalization of the Zariski closure of $\Sh(G_{\dbQ},X)/H$ in $\scrM_{O_{(v)}}$. The natural actions of $G_{\dbQ}(\dbA_f^{(p)})$ on $\scrM_{\dbQ}=\Sh(\GSp(W,\psi),S)/K_p$ and on $\Sh(G_{\dbQ},X)/H$ give rise to a natural action of $G_{\dbQ}(\dbA_f^{(p)})$ on $\scrN$. Let $(\scrA,\scrP_{\scrA})$ be the pullback to $\scrN$ of the universal principally polarized abelian scheme over $\scrM$. 

Let $(v_{\alpha})_{\alpha\in\scrJ}$ be a family of tensors of $\scrT(L_{(p)}^*)$ such that $G_{\dbQ}$ is the subgroup of $\GL_W$ that fixes $v_{\alpha}$ for all $\alpha\in\scrJ$, cf. [9, Prop. 3.1 c)]. The choices of $L$ and $(v_{\alpha})_{\alpha\in\scrJ}$ allow a moduli interpretation of $\Sh(G_{\dbQ},X)(\dbC)$ (see [7], [8], [31,] and [45]). For instance,
$\Sh(G_{\dbQ},X)(\dbC)=G_{\dbQ}(\dbQ)\backslash X\times G_{\dbQ}(\dbA_f)$ is the set of isomorphism classes of principally polarized abelian varieties over $\dbC$ that are of dimension $r$, that carry a family of Hodge
cycles indexed by $\scrJ$, that have level $s$ symplectic similitude structures for all $s\in\dbN$, and that satisfy some additional conditions (for instance, see [45, Subsect. 4.1]). Thus the abelian scheme $\scrA_{E(G_{\dbQ},X)}$ is naturally endowed with a family $(w_{\alpha}^{\scrA})_{\alpha\in\scrJ}$ of Hodge cycles (the Betti realizations of pullbacks of $w_{\alpha}^{\scrA}$ via $\dbC$-valued points of $\scrN_{E(G_{\dbQ},X)}$ correspond to $v_{\alpha}$).

Let $G^0_{\dbZ_{(p)}}$ be the maximal reductive, closed subgroup scheme of $G_{\dbZ_{(p)}}$ that fixes $\psi$. The maximal compact subtorus of $\dbS$ is connected and, when viewed as a subgroup of $G_{\dbR}$ via the image of an arbitrary element $x\in X$, it has $-1_{W\otimes_{\dbQ} \dbR}$ as an $\dbR$-valued point. Thus $Z(\Sp(L_{(p)},\psi))\leqslant G_{\dbZ_{(p)}}^0$ and therefore we have an identity
$$G_{\dbZ_{(p)}}^0=\Sp(L_{(p)},\psi)\cap G_{\dbZ_{(p)}}.\leqno (12)$$ 
\noindent
{\bf 5.1.1. Standard Shimura $F$-crystals.}  The isomorphism $L_{(p)}\arrowsim L^*_{(p)}$ induced by $\psi$ allows us to naturally identify $G_{\dbZ_{(p)}}$
with a closed subgroup scheme of $\GL_{L^*_{(p)}}$. Let $L_p:=L_{(p)}\otimes_{\dbZ_{(p)}} \dbZ_p$. Let $T^0$ be a maximal torus of a Borel subgroup scheme $B^0$ of $G_{\dbZ_p}$. Let $\mu_0:\dbG_m\to T^0_{W(k(v))}$ be a cocharacter such that the following two axioms hold:
\medskip
{\bf (a)} over an embedding $e_{k(v)}:W(k(v))\hookrightarrow\dbC$ which extends the composite inclusion $O_{(v)}\subseteq E(G_{\dbQ},X)\subseteq\dbC$, it is
$G_{\dbC}(\dbC)$-conjugate to the cocharacters $\mu_x:\dbG_m\to G_\dbC$, $x\in X$;
\smallskip
{\bf (b)} if $L^*_p\otimes_{\dbZ_p} W(k(v))=F^1_0\oplus F^0_0$ is the direct sum decomposition normalized by $\mu_0$ and such that $\mu_0$ acts trivially on $F^0_0$, then $B_{W(k(v))}$ normalizes $F^1_0$.  
\medskip
The existence of $\mu_0$ is implied by [31, Prop. 4.6 and Cor. 4.7]. Let $M_0:=L^*_p\otimes_{\dbZ_p} W(k(v))$, let $\phi_0:=(1_{L^*_p}\otimes\sigma_{k(v)})\circ\mu_0({1\over p})$, and let
$$
\scrC_0:=(M_0,\phi_0,G_{W(k(v))}).
$$
We first introduced $\scrC_0$ in [44]. Up to isomorphisms of $M_0$ defined by elements of $G_{\dbZ_p}(\dbZ_p)$, $\scrC_0$ does not depend on the choice of either $T^0$ or $B^0$.

Let $k$ be an algebraically closed field of characteristic $p$. Let $G:=G_{W(k)}$. For $g\in G(W(k))$ let $\scrC_{g}:=(M_0\otimes_{W(k(v))} W(k),g(\phi_0\otimes\sigma),G,(v_{\alpha})_{\alpha\in\scrJ})$. We refer to the family $(\scrC_{g})_{g\in G(W(k))}$ as the {\it standard Shimura $F$-crystals with tensors} over $k$ of $(f,L,v)$.

\medskip\noindent
{\bf 5.1.2. Definition.} The potential standard Hodge situation $(f,L,v)$ is called a {\it standard Hodge situation} if the following two axioms hold:
\medskip
{\bf (a)} the scheme $\scrN$ is regular and formally smooth over $O_{(v)}$;
\smallskip
{\bf (b)} for each field  $k$ as above and for every point $z\in\scrN(W(k))$, the quadruple $(M,F^1,\phi,\tilde G)$ is a Shimura filtered $F$-crystal; here $(M,\phi)$ is the Dieudonn\'e module of the special fibre of the $p$-divisible group of $A:=z^*(\scrA)$, $F^1$ is the Hodge filtration of $M=H^1_{\text{dR}}(A/W(k))$ defined by $A$, and $\tilde G$ is the Zariski closure in $\GL_M$ of the subgroup of $\GL_{M[{1\over p}]}$ that fixes the de Rham component $t_{\alpha}\in\scrT(M[{1\over p}])=\scrT(H^1_{\text{dR}}(A_{B(k)}/B(k)))$ of $z^*(w_{\alpha}^{\scrA})$ for all $\alpha\in\scrJ$. 
\medskip
Let $\psi_M$ be the perfect bilinear form on $M$ that is the de Rham (crystalline) realization of $z^*(\scrP_{\scrA})$. We have $\psi_M(\phi(x),\phi(y))=p\sigma(\psi_M(x,y))$ for all $x$, $y\in M$.

\medskip\noindent
{\bf 5.1.3. Remark.} If the axiom 5.1.2 (a) holds, then $\scrN$ is an integral canonical model of $(G_{\dbQ},X,H,v)$ in the sense of [45, Subsubsects. 3.2.3 to 3.2.6] (cf. [45, Cor. 3.4.4]). If $p\ge 5$, then the axiom 5.1.2 (a) holds (cf. [45, Subsubsects. 3.2.12, 3.4.1, and 6.4.2]). More recently, in [50] (resp. in [20]) it is claimed that each potential standard situation is a standard Hodge situation for all primes $p\Ge 2$ (resp. provided $p\Ge 3$).

\medskip\noindent
{\bf 5.1.4. Example.} We check that if $G_{\dbQ}$ is a torus, then $(f,L,v)$ is a standard Hodge situation. Standard class field theory implies that $\scrN$ is a pro-\'etale cover of $\Spec(O_{(v)})$ (see [45, Ex. 3.2.8]); thus axiom 5.1.2 (a) holds. Let $T_{\dbZ_{(p)}}^{\text{big}}$ be a maximal torus of $\GL_{L_{(p)}}$ that contains the torus $G_{\dbZ_{(p)}}$. We can assume that $\scrJ$ is such that each element of $\Lie(T_{\dbZ_{(p)}}^{\text{big}})$ is a $v_{\alpha}$ for some $\alpha\in\scrJ$. The $W(k)$-span of endomorphisms of $M$ generated by those elements $t_{\alpha}$ which correspond to elements of $\Lie(T_{\dbZ_{(p)}}^{\text{big}})$, is a $W(k)$-subalgebra of $\End(M)$ isomorphic to $W(k)^{2d}$ and therefore it is the Lie algebra of a maximal torus of $\GL_M$ that contains $\tilde T$. This implies that $\tilde T$ is a torus and thus the axiom 5.1.2 (b) holds as well. 

\medskip\noindent
{\bf 5.1.5. Fontaine comparison theory.} If $p\ge 3$ (resp. if $p=2$) let $B^+(W(k))$ be the Fontaine ring used in [13, Sect. 4] (resp. in [13, Sect. 8] and obtained using the $2$-adic completion). We recall that $B^+(W(k))$ is an integral, local $W(k)$-algebra endowed with a separated and decreasing filtration $\bigl(F^i(B^+(W(k))\bigr)_{i\in\dbN\cup\{0\}}$, with a Frobenius lift, and with a $\Gal(B(k))$-action (see also [15, Sect. 2]). 
\smallskip
Until the end we will assume $(f,L,v)$ is a standard Hodge situation and we will use the notations of the axiom 5.1.2 (b). Let $H^1:=H^1_{\acute et}(A_{\overline{B(k)}},\dbZ_p)$. There exists a standard but non-canonical identification $H^1=L_{p}^*$ under which: (i) the $p$-component of the \'etale component of $z^*(w_{\alpha}^{\scrA})$ is $v_{\alpha}$ for all $\alpha\in\scrJ$, and (ii) the $\dbQ_p$-\'etale realization $\lambda_{H^1}$ of $z^*(\scrP_{\scrA})$ is a $\dbG_m(\dbZ_{(p)})$-multiple $\gamma_p$ of the perfect form $\psi^*$ on $L_{p}^*$ defined by $\psi$ via duality (see [45, top of p. 473]). As the complex $1\to G^0_{\dbZ_{(p)}}(\dbZ_p)\to G_{\dbZ_{(p)}}(\dbZ_p)\to\dbG_m(\dbZ_p)\to 1$ is exact, we can assume $\gamma_p=1$. From Fontaine comparison theory we get the existence of a $B^+(W(k))$-linear monomorphism
$$i_A:M\otimes_{W(k)} B^+(W(k))\hookrightarrow H^1\otimes_{\dbZ_p} B^+(W(k))\leqno (13)$$
that respects the tensor product Frobenius endomorphisms (the Frobenius endomorphism of $H^1$ being $1_{H^1}$). The existence of $i_A$ is a particular case of [13, Thm. 7] (for $p=2$ cf. also end of [13, Sect. 8]); see also [49, Subsubsect. 2.2.2] which works with Fontaine rings $B^+(W(k))$ that are defined for all primes $p\ge 2$ in the same way. 

For de Rham rings of periods and for the de Rham conjecture we refer to [12], [13], [15], and [2]. For de Rham cycles we refer to [2, Subsect. (1.3)]. 

\medskip\noindent
{\bf 5.1.6. Theorem.} {\it Let $\scrV$ be a complete discrete valuation ring of mixed characteristic $(0,p)$ and perfect residue field. Let $\scrF:=\scrV[{1\over p}]$. Let $\scrZ\in\scrN(\scrF)$. Let $\scrE:=\scrZ^*(\scrA)$. Then each Hodge cycle $w_0$ that involves no \'etale Tate twist on $\scrE$, is a de Rham cycle. In other words, under the isomorphism of the de Rham conjecture, the de Rham component of $w_0$ is mapped to the $p$-component of the \'etale component of $w_0$.}

\medskip
If $\scrE$ is definable over a number subfield of $\scrF$, this was known since long time (for instance, see [2, Thm. (0.3)]). The general case follows from loc. cit. and [45, Principle B of Subsubsect. 5.2.16] (in the part of [45, Subsect. 5.2] preceding the Principle B an odd prime is used; however the proof of loc. cit. applies to all primes). 
\medskip\noindent
{\bf 5.1.7. Corollary.} {\it We have $\phi(t_{\alpha})=t_{\alpha}$ for all $\alpha\in\scrJ$.}

\medskip
\proof The cycle $w_{\alpha}^{\scrA}$ involves no Tate twist. Thus the tensorization of (13) with the de Rham ring $B_{\text{dR}}(W(k))$ that contains $B^+(W(k))$, is an isomorphism that takes $t_{\alpha}$ to $v_{\alpha}$ (cf. Theorem 5.1.6). If $K^+(W(k))$ is the field of fractions of $B^+(W(k))$, then the isomorphism $\scrT(M[{1\over p}])\otimes_{B(k)} K^+(W(k))\arrowsim\scrT(H^1[{1\over p}])\otimes_{\dbQ_p} K^+(W(k))$ induced naturally by $i_A$, takes $t_{\alpha}$ to $v_{\alpha}$ and is compatible with Frobenius endomorphisms. Thus each $t_{\alpha}$ is fixed by $\phi$. \endproof

\medskip
By multiplying each $v_{\alpha}$ with a power of $p$ we can assume that we have $t_{\alpha}\in\scrT(M)$ for all $k$ and $z$ as in the axiom 5.1.2 (b). Let $H_0$ be a compact, open subgroup of $G_{\dbQ}(\dbA_f^{(p)})$ that acts freely on $\scrM$ (cf. Serre's Lemma). The group $H_0$ also acts freely on $\scrN$. Thus $\scrN$ is a pro-\'etale cover of $\scrN/H_0$, cf. proof of [45, Prop. 3.4.1]. Let $\Spec(Y_0)$ be an affine, open subscheme of $\scrN/H_0$ such that the $k(v)$-algebra has a finite $p$-basis and the point $y\in\scrN(k)$ defined by the special fibre of $z$ factors through the affine, open subscheme $\Spec(Y):=\Spec(Y_0)\times_{\scrN/H_0} \scrN$ of $\scrN$. Let $\Phi_Y$ be a Frobenius lift of the $p$-adic completion $Y^\wedge$ of $Y$. Let $\scrA_Y$ and $w_{\alpha}^{Y}$ be the pullbacks of $\scrA$ and $w_{\alpha}^{\scrA}$ (respectively) to $\Spec(Y)$. The de Rham component $t_{\alpha}^Y$ of $w_{\alpha}^Y$ is annihilated by the Gauss--Manin connection on $H^1_{\text{dR}}(\scrA_Y/Y)$ (see [9, Prop. 2.5]) and it is fixed by $\Phi_Y$ (cf. Corollary 5.1.7 applied to dominant Teichm\"uller lifts of $\Spec(Y^\wedge)$). Thus, as each $F$-crystal over $Y/pY$ is defined by its evaluation at the natural thickening attached to the closed embedding $\Spec (Y/pY)\hookrightarrow \Spec(Y)$, each $t_{\alpha}^Y$ is a crystalline tensor of $\scrT(H^1_{\text{crys}}(\scrA_{Y/pY}/Y^\wedge))=\scrT(H^1_{\text{dR}}(\scrA_Y/Y)\otimes_Y Y^\wedge)$. Thus $t_{\alpha}\in\scrT(M)$ depends only on $y$ and not on $z$. This justifies the following terminology and notations.

\medskip\noindent
{\bf 5.1.8. Notations.} Let $y\in\scrN(k)$ be defined by the special fibre of $z\in\scrN(W(k))$. We refer to $\scrC_y:=(M,\phi,\tilde G,(t_{\alpha})_{\alpha\in\scrJ})$ (resp. to $(M,\phi,\tilde G)$) as the Shimura $F$-crystal with tensors (resp. as the Shimura $F$-crystal) attached to $y$. Let $(M_{\dbZ_p},\tilde G_{\dbZ_p},(t_{\alpha})_{\alpha\in\scrJ})$ be the $\dbZ_p$ structure of $(M,\tilde G,(t_{\alpha})_{\alpha\in\scrJ})$ defined by $\phi\mu(p)$ (cf. Subsection 2.6), where $\mu:\dbG_m\to \tilde G$ is an arbitrary Hodge cocharacter of $(M,\phi,\tilde G)$. 
\medskip\noindent
{\bf 5.1.9. Lemma.} {\it We fix a Hodge cocharacter $\mu$ of $(M,\phi,\tilde G)$. There exists a $B(k)$-linear isomorphism 
$$\rho:M[{1\over p}]\arrowsim M_0\otimes_{W(k(v))} B(k)=L_p^*\otimes_{\dbZ_p} B(k)=H^1\otimes_{\dbZ_p} B(k)$$
 that takes $t_\alpha$ to $v_\alpha$ for all $\alpha\in\scrJ$ and such that $\rho\mu_{B(k)}\rho^{-1}=\mu_{0B(k)}$. In particular, we have $\rho(F^1[{1\over p}])=F^1_0\otimes_{W(k(v)} B(k)$ and there exists $g_y\in G(B(k))$ such that $\rho\phi\rho^{-1}=g_y(\phi_0\otimes\sigma_k)$.}

\medskip
\proof
 As $k=\bar k$, the set $H^1(B(k),G_{B(k)})$ has only one class. As $i_A\otimes_{B^+(W(k))} 1_{B_{\text{dR}}(W(k))}$ is an isomorphism that takes $t_\alpha$ to $v_\alpha$ for all $\alpha\in\scrJ$, we get that there exists a $B(k)$-linear isomorphism $\rho:M[{1\over p}]\arrowsim M_0\otimes_{W(k(v))} B(k)$
 that takes $t_\alpha$ to $v_\alpha$ for all $\alpha\in\scrJ$. We check that we can choose $\rho$ such that the cocharacters $\rho\mu_{B(k)}\rho^{-1}$ and $\mu_{0B(k)}$ of $G_{B(k)}$ coincide. To check this we can assume that there exists an embedding $e_k:W(k)\hookrightarrow\dbC$ that extends the embedding $e_{k(v)}:W(k(v))\hookrightarrow\dbC$ of the axiom 5.1.1 (a). Thus $\rho\mu_{B(k)}\rho^{-1}$ and $\mu_{0B(k)}$ are $G(\dbC)$-conjugate, cf. axiom 5.1.1 (a). As $G_{B(k)}$ is split, we get that the cocharacters $\rho\mu_{B(k)}\rho^{-1}$ and $\mu_{0B(k)}$ are $G(B(k))$-conjugate. Thus by composing $\rho$ with an automorphism of $M_0\otimes_{W(k(v))} B(k)$ which is an element of $G(B(k))$, we can assume that we have $\rho\mu_{B(k)}\rho^{-1}=\mu_{0B(k)}$. The last part of the Lemma is obvious. \endproof

\medskip\noindent
{\bf 5.1.10. Definition.} Let $y$, $y_1\in\scrN(k)$. By an isomorphism (resp. a rational isomorphism) between $\scrC_y=(M,\phi,\tilde G,(t_{\alpha})_{\alpha\in\scrJ})$ and $\scrC_{y_1}=(M_1,\phi_1,\tilde G_1,(t_{1\alpha})_{\alpha\in\scrJ})$ we mean an isomorphism $h:M\arrowsim M_1$ (resp. $h:M[{1\over p}]\arrowsim M_1[{1\over p}]$) that satisfies $h\phi=\phi_1 h$ and that takes $t_{\alpha}$ to $t_{1\alpha}$ for all $\alpha\in\scrJ$. Similarly we define (rational) isomorphisms  between $\scrC_y$ and $\scrC_{g}$ or between two Shimura $F$-crystals with tensors indexed by the same set $\scrJ$.

\medskip\noindent
{\bf 5.1.11. Remark.} The notion of rational isomorphisms between $\scrC_y$'s and $\scrC_{g}$'s does not depend on the choice of the family of tensors $(v_{\alpha})_{\alpha\in\scrJ}$ in Subsection 5.1. Also, as the right translations of $y\in\scrN(k)$ by elements of $G_{\dbQ}(\dbA_f^{(p)})$ correspond at the level of abelian varieties to passages to isogenies of $A_k$ prime to $p$ (see [45, Subsect. 4.1] for the case of complex points of $\scrN$), the isomorphism class of $\scrC_y$ depends only on the $G_{\dbQ}(\dbA_f^{(p)})$-orbit of $y$. Under such isogenies the cycle $\psi_M$ gets replaced by suitable $\dbG_m(\dbZ_p)$-multiples of it.

\medskip\smallskip\noindent
{\bf 5.2. Manin problem for $(f,L,v)$ and $k$.} Show that the following two things hold:

\medskip
{\bf (a)} For each point $y\in\scrN(k)$, there exists an element $g\in G(W(k))$ such that $\scrC_g$ is rational isomorphic to $\scrC_y$.

\smallskip
{\bf (b)} For each element $g\in G(W(k))$, there exists a point $y\in\scrN(k)$ such that $\scrC_y$ is rational isomorphic to $\scrC_g$.

\medskip\noindent
{\bf 5.2.1. Assumptions.} Below we often assume that of the following two conditions hold:

\medskip
{\bf (i)} the $GAL$ property holds for $(f,L,v)$ and $k$ i.e., for each $z\in\scrN(W(k))$ there exists a $\dbQ_p$-linear isomorphism $\rho_1:M_{\dbZ_p}[{1\over p}]\arrowsim L_p^*[{1\over p}]=H^1[{1\over p}]$ that takes $t_{\alpha}$ to $v_{\alpha}$ for all $\alpha\in\scrJ$;  

\smallskip
{\bf (ii)} the $GAL$ property of (i) holds for each standard Hodge situation $(f^\prime,L,v_1)$  and $k$, where $f^\prime$ is the composite of an injective map $(T_{\dbQ}^\prime,\{h^\prime\})\hookrightarrow (G_{\dbQ},X)$ with $f$, for $T_{\dbQ}^\prime$ a maximal torus of $G_{\dbQ}$, and where $v_1$ is an arbitrary prime of $E(T_{\dbQ}^\prime,\{h^\prime\})$ that divides $v$.

\medskip\noindent
{\bf 5.2.2. Remarks.} {\bf (a)} Let $z\in\scrN(W(k))$. The existence of a $\dbQ_p$-linear isomorphism $\rho_1:M_{\dbZ_p}[{1\over p}]\arrowsim L_p^*[{1\over p}]=H^1[{1\over p}]$ that takes $t_{\alpha}$ to $v_{\alpha}$ for all $\alpha\in\scrJ$ is encoded in the fact that a suitable class $\gamma_z\in H^1(\dbQ_p,G_{\dbQ_p})$ is trivial, cf. Lemma 5.1.9. If $G^{\der}$ is simply connected, then such classes were computed in [38, Subsect. 1.19 and Prop. 1.20].

We get that the condition 5.2.1 (i) holds if the set $H^1(\dbQ_p,G_{\dbQ_p})$ has only one class. For instance, this is so if $G_{\dbQ_p}^{\der}$ is simply connected and $G^{\ab}_{\dbQ_p}$ is split. 

\smallskip
{\bf (b)} As $G_{\dbZ_p}$ is smooth and has a connected special fibre, each torsor of $G_{\dbZ_p}$ in the flat topology is trivial (see Subsection 2.6). Thus there exists an isomorphism $M_{\dbZ_p}\arrowsim L_p^*=H^1$ that takes $t_{\alpha}$ to $v_{\alpha}$ for all $\alpha\in\scrJ$ if and only if there exists an isomorphism $M\arrowsim L_p^*\otimes_{\dbZ_p} W(k)=H^1\otimes_{\dbZ_p} W(k)$ that takes $t_{\alpha}$ to $v_{\alpha}$ for all $\alpha\in\scrJ$ (cf. also [49, Lemma 2.5.2 (a)]). 

\smallskip
{\bf (c)} Either a particular case of [49, Thm. 1.2] or [20, Cor. (1.4.3)] implies that if either $G_{\dbQ}$ is a torus or $p\ge 3$, then there exists an isomorphism $M\arrowsim L_p^*\otimes_{\dbZ_p} W(k)=H^1\otimes_{\dbZ_p} W(k)$ that takes $t_{\alpha}$ to $v_{\alpha}$ for all $\alpha\in\scrJ$. From this and (b) we get that the condition 5.2.1 (i) holds provided $p\Ge 3$ and that the condition 5.2.1 (ii) always holds for all $p\Ge 2$.

\medskip\noindent
{\bf 5.2.3. Main Theorem.} {\it Let $(f,L,v)$ be a standard Hodge situation. We have:

\medskip
{\bf (a)} If the condition 5.2.1 (i) holds, then the property 5.2 (a) holds. 

\smallskip
{\bf (b)} We assume either (i) that the condition 5.2.1 (ii) holds  or (ii) that the condition 5.2.1 (i) holds, that $Z(G_{\dbQ})$ is connected, and that the group $G_{\dbQ_p}(\dbQ_p)$ surjects onto $G^{\ad}_{\dbQ_p}(\dbQ_p)$. Then the property 5.2 (b) holds.

\smallskip
{\bf (c)} If $p\Ge 3$, then properties 5.2 (a) and (b) hold.}

\medskip
\proof Part (c) follows from (a) and (b) and from Remark 5.2.2 (c). 

We prove (a). Let $y\in\scrN(k)$. Let $z\in\scrN(W(k))$ be a point that lifts $y$. We use the notations of the axiom 5.1.2 (b) and of Notations 5.1.8. Two maximal tori of Borel subgroups of $\tilde G_{\dbQ_p}$ are $\tilde G_{\dbQ_p}(\dbQ_p)$-conjugate (see [3, Ch. V, Thms. 19.2 and 20.9 (i)]). Thus we can choose the $\dbQ_p$-linear isomorphism $\rho_1$ of the condition 5.2.1 (i) such that the maximal subtorus $\rho_1^{-1}T^0_{\dbQ_p}\rho_1$ of $\tilde G_{\dbQ_p}$ extends to a maximal torus $\tilde T^0:=\rho_1^{-1}T^0\rho_1$ of $\tilde G_{\dbZ_p}$. Let $N^0$ be the normalizer of $T^0$ in $G_{\dbZ_p}$. Thus $\tilde N^0:=\rho_1^{-1}N^0\rho_1$ is the normalizer of $\tilde T^0$ in $\tilde G_{\dbZ_p}$. From Theorem 1.3.3 (b) and Subsection 4.5 we get the existence of an element $\tilde w\in \tilde N^0(W(k))$ and of a $\tilde G(W(k))$-conjugate $\tilde\mu:\dbG_m\to\tilde T^0_{W(k)}$ of $\mu:\dbG_m\to \tilde G$, such that $\scrC_y$ is rational isomorphic to $(M,\tilde w(1_{M_{\dbZ_p}}\otimes\sigma)\tilde\mu({1\over p}),\tilde G,(t_{\alpha})_{\alpha\in\scrJ})$. Under $\rho_{1B(k)}$, this last Shimura $F$-crystal with tensors is rational isomorphic to $\scrC_{1y}:=(M_0\otimes_{W(k(v))} W(k),w(1_{L^*_p}\otimes\sigma)\tilde\mu_0({1\over p}),G,(v_{\alpha})_{\alpha\in\scrJ})$, where $w:=\rho_{1B(k)}\tilde w\rho_{1B(k)}^{-1}\in N^0(W(k))$ and $\tilde\mu_0:=\rho_{1B(k)}\tilde\mu\rho_{1B(k)}^{-1}:\dbG_m\to T^0_{W(k)}$. As $\rho\rho_{1B(k)}^{-1}\in G(B(k))$, from Lemma 5.1.9 we get that the generic fibres of $\tilde\mu_0$ and $\mu_{0W(k)}$ are $G(B(k))$-conjugate. Thus, as $G$ is split, the cocharacters $\tilde\mu_0$ and $\mu_{0W(k)}$ are also conjugate under an element $h\in G(W(k))$. But $hw\sigma(h^{-1})\in G(W(k))$ and thus $\scrC_{1y}$ is rational isomorphic to $\scrC_{hw\sigma(h^{-1})}$ under an isomorphism defined by $h$. Therefore $\scrC_y$ is rational isomorphic to $\scrC_{hw\sigma(h^{-1})}$ i.e., the property 5.2 (a) holds. This proves (a).

We prove (b). Let $\dbF:=\overline{\dbF_p}$. To show that the property 5.2 (b) holds we can assume that $g\in N^0(W(\dbF))$, cf. Theorem 1.3.3 (b). Let $w:=g\in G(W(k))$. Let $T^w$ be a maximal torus of $G_{\dbZ_p}$ such that the $\dbZ_p$ structure of the following quadruple $(M_0\otimes_{W(k(v))} W(\dbF),T^0_{W(\dbF)},G_{W(\dbF)},(v_{\alpha})_{\alpha\in\scrJ})$ defined by the $\sigma$-linear automorphism $w(1_{L_p^*}\otimes\sigma)$ of $M_0\otimes_{W(k(v))} W(\dbF)=L_p^*\otimes_{\dbZ_p} W(\dbF)$ (see Subsection 2.6), is isomorphic to the quadruple $(L_p^*,T^w,G_{\dbZ_p},(v_{\alpha})_{\alpha\in\scrJ})$ under an isomorphism defined by an element  $h_w\in G_{\dbZ_p}(W(k))$. Thus $\mu^w_{0W(\dbF)}:=h_w\mu_{0W(\dbF)}h_w^{-1}$ is a cocharacter of $T^w_{W(\dbF)}$. As the proof of (b) is quite lengthy, we divide it into six boldfaced parts as follows. 

\smallskip
{\bf A Shimura pair of dimension $0$.} Let $T_w$ be a maximal $\dbZ_{(p)}$-torus of $G_{\dbZ_{(p)}}$ such that the following two things hold (cf. [17, Lemma 5.5.3]):
\medskip
{\bf (i)} over $\dbR$ it is the extension of a compact torus by $Z(\GL_{W\otimes_{\dbQ} \dbR})$;
\smallskip
{\bf (ii)} there exists $g_w\in G_{\dbZ_p}(\dbZ_p)$ such that $g_wT^wg_w^{-1}=T_{w\dbZ_p}$.
\medskip 

Let $\mu_w:\dbG_m\to T_{w\dbC}$ be the cocharacter obtained from the cocharacter $g_w\mu^w_{0W(\dbF)}g_w^{-1}$ of $g_wT^w_{W(\dbF)}g_w^{-1}=T_{wW(\dbF)}$ by extension of scalars under a fixed embedding $e_{\dbF}:W(\dbF)\hookrightarrow\dbC$ that extends the embedding $e_{k(v)}:W(k(v))\hookrightarrow\dbC$ of the axiom 5.1.1 (a). From the axiom 5.1.1 (a) we get:
\medskip
{\bf (iii)} $\mu_w$ as a cocharacter of $G_{\dbC}$ is $G_{\dbQ}(\dbC)$-conjugate to $\mu_x$, where  $x\in X$.
\medskip
Let $S_{w\dbC}$ be the subtorus of $T_{w\dbC}$ generated by $Z(\GL_{W\otimes_{\dbQ} \dbC})$ and by the image of $\mu_w$. As $T_{w\dbR}/Z(\GL_{W\otimes_{\dbQ} \dbR})$ is compact, each subtorus of the extension of $T_{w\dbR}/Z(\GL_{W\otimes_{\dbQ} \dbR})$ to $\dbC$ is defined over $\dbR$. Thus $S_{w\dbC}$ is the extension to $\dbC$ of a subtorus $S_{w\dbR}$ of $T_{w\dbR}$. From (iii) we get that $Z(\GL_{W\otimes_{\dbQ} \dbC})$ and the image of $\mu_w$ have the finite, \'etale group $\pmb{\mu}_2$ as their intersection. Thus we can identify naturally $S_{w\dbR}=\dbS$. Let  
$$x_w:\dbS\hookrightarrow T_{w\dbR}$$
be the resulting monomorphism. The Shimura pair $(T_w,\{x_w\})$ has dimension $0$.

\smallskip
{\bf A possible Shimura pair.} Let $X_w$ be the $G_{\dbQ}(\dbR)$-conjugacy class of $x_w$. The Hodge $\dbQ$--structure on $\Lie(G_{\dbQ})$ defined by $x_w$ is of type $\{(-1,1),(0,0),(1,-1)\}$, cf. (iii). Thus $(G_{\dbQ},X_w)$ is a Shimura pair if and only if $\theta\circ x_w(i)$ defines a Cartan involution of $\Lie(G^{\ad}_{\dbR})$. In general $\theta\circ x_w(i)$ is not a Cartan involution of $\Lie(G^{\ad}_{\dbR})$. Here is an example. If $G^{\ad}_{\dbR}$ is an $\pmb{\text{SO}}(2,2n+1)_{\dbR}$ group with $n\Ge 2$, then it has two $\pmb{\text{SO}}(2)_{\dbR}$ subgroups $F_1$ and $F_2$ which over $\dbC$ are $G^{\ad}_{\dbR}(\dbC)$-conjugate and whose centralizers in $G^{\ad}_{\dbR}$ are isomorphic to $\pmb{\text{SO}}(2)_{\dbR}\times_{\dbR} \pmb{\text{SO}}(2n+1)_{\dbR}$ and $\pmb{\text{SO}}(2)_{\dbR}\times_{\dbR} \pmb{\text{SO}}(2,2n-1)_{\dbR}$ (respectively); if the image of $x_w$ in $G^{\ad}_{\dbR}$ is $F_2$, then $\theta\circ x_w(i)$ does not define a Cartan involution of $\Lie(G^{\ad}_{\dbR})$. 

\smallskip
{\bf Twisting process.} We use a {\it twisting process} to show that we can choose $T_w$ such that $(G_{\dbQ},X_w)$ is a Shimura pair. Let $T^{\star}_w$ be the image of $T_w$ in $G^{\ad}_{\dbZ_{(p)}}$. The functorial map $H^1(\dbR,T^{\star}_{w\dbR})\to H^1(\dbR,G^{\ad}_{\dbR})$ is surjective (cf. [22, Lemma 10.1]) and the complex 
$$H^1(\dbQ,T^{\star}_{w\dbQ})\to H^1(\dbR,T^{\star}_{w\dbR})\times \prod_{l\,\text{a}\,\text{prime}} H^1(\dbQ_l,T^{\star}_{w\dbQ_l})\to X_*(T^{\star}_{w\overline{\dbQ}})_{\Gal(\dbQ),\text{tors}}$$
of abelian groups is exact (cf. [30, Thm. B.24]). If the prime $l$ is different from $p$ and such that the torus $T^{\star}_{w\dbQ_l}$ is split, then the group $H^1(\dbQ_l,T^{\star}_{w\dbQ_l})$ surjects onto $X_*(T^{\star}_{w\overline{\dbQ}})_{\Gal(\dbQ),\text{tors}}$ (cf. [30, Prop. B.22]). Therefore we have a natural epimorphism 
$$H^1(\dbQ,T^{\star}_{w\dbQ})\twoheadrightarrow H^1(\dbR,T^{\star}_{w\dbR})\times H^1(\dbQ_p,T^{\star}_{w\dbQ_p}).\leqno (14)$$ 
We recall that $G^{\ad}_{\dbR}$ is an inner form of its compact form, cf. [8, Subsubsect. 2.3.4]. From this, the surjectivity of $H^1(\dbR,T^{\star}_{w\dbR})\to H^1(\dbR,G^{\ad}_{\dbR})$, and (14) we get that there exists a class $\gamma_c\in H^1(\dbQ,T^{\star}_{w\dbQ})$ whose image in  $H^1(\dbQ_p,T^{\star}_{w\dbQ_p})$ is trivial and such that the inner twist $G^c_{\dbQ}$ of $G_{\dbQ}$ via $\gamma_c$ is compact over $\dbR$. The group $T_{w\dbQ}$ is naturally a maximal torus of $G^c_{\dbQ}$ and we have an identification $G_{\dbQ_p}=G^c_{\dbQ_p}$ that extends the identity automorphism of $T_{w\dbQ_p}$.

Let $\grN$ be the set of non-compact, simple factors of $G^{\ad}_{\dbR}$. For $j\in \grN$ let $I^{\star j}_{w\dbR}$ be the image of $S_{w\dbR}$ in the non-compact, simple factor $j$ of $G^{\ad}_{\dbR}$. Let $I^{\star}_{w\dbR}:=\prod_{j\in \grN} I^{\star j}_{w\dbR}$; it is a compact subtorus of $T^{\star}_{w\dbR}$ and thus also of $G^c_{\dbR}$. Let $\gamma_{s\dbR}^\star\in H^1(\dbR,I^{\star}_{w\dbR})=\prod_{j\in\grN} H^1(\dbR,I^{\star j}_{w\dbR})$ be the unique class such that all its components are non-trivial classes. The inner twist of $G^{c\ad}_{\dbR}$ via $\gamma_{s\dbR}^\star$ is isomorphic to $G^{\ad}_{\dbR}$, cf. [8, Subsubsect. 1.2.3]. Thus the inner twist of $G^c_{\dbR}$ via $\gamma_{s\dbR}^\star$ is isomorphic to $G_{\dbR}$. Let $\gamma_s\in H^1(\dbQ,T^{\star}_{w\dbQ})$ be a class such that its image in $H^1(\dbQ_p,T^{\star}_{w\dbQ_p})$ is trivial and its image in $H^1(\dbR,T^{\star}_{w\dbR})$ is the image of $\gamma_{s\dbR}^\star$ in $H^1(\dbR,T^{\star}_{w\dbR})$, cf. (14). Let $G^\prime_{\dbQ}$ be the inner twist of $G^c_{\dbQ}$ via $\gamma_s$. We have natural isomorphisms $j_{\dbR}:G^\prime_{\dbR}\arrowsim G_{\dbR}$ and $G^\prime_{\dbQ_p}\arrowsim G_{\dbQ_p}$ and moreover $T^{\star}_{w\dbQ}$ and $T_{w\dbQ}$ are maximal tori of $G^{\prime\ad}_{\dbQ}$ and $G^\prime_{\dbQ}$ (respectively). Let $X_w^\prime$ be the $G^\prime_{\dbQ}(\dbR)$-conjugacy class of the composite of $x_w$ with the inclusion $T_{w\dbR}\hookrightarrow G^\prime_{\dbR}$. Let $X_w^{\prime\ad}$ be the $G^{\prime\ad}_{\dbQ}(\dbR)$-conjugacy class of the composite of $x_w$ with the homomorphism $T_{w\dbR}\to G^{\prime\ad}_{\dbR}$. The pair $(G^{\prime\ad}_{\dbQ},X_w^{\prime\ad})$ is a Shimura pair, cf. [8, Subsubsect. 1.2.3]. Thus $(G^{\prime}_{\dbQ},X_w^{\prime})$ is also a Shimura pair. 

The class $\gamma_0:=\gamma_c^{-1}\gamma_s^{-1}\in H^1(\dbQ,T^{\star}_{w\dbQ})$ is trivial over $\dbQ_p$. The inner twist of $G^{\prime\ad}_{\dbR}$ via the image $\gamma_{0\dbR}$ of $\gamma_0$ in $H^1(\dbR,G^{\prime\ad}_{\dbR})$ is $G^{\ad}_{\dbR}$ and therefore (via $j_{\dbR}$) it is naturally isomorphic to $G^{\prime\ad}_{\dbR}$. Thus $\gamma_{0\dbR}$ is the trivial class. Let $N^{\star\prime}_{w\dbQ}$ be the normalizer of $T^{\star}_{w\dbQ}$ in $G^{\prime\ad}_{\dbQ}$. Let $\gamma_n\in H^1(\dbQ,N^{\star\prime}_{w\dbQ})$ be the class that parametrizes inner isomorphisms between the two pairs $(G^{\prime\ad}_{\dbQ},T^{\star}_{w\dbQ})$ and $(G^{\ad}_{\dbQ},T^{\star}_{w\dbQ})$. The use of word inner makes sense here as the class $\gamma_{0\dbR}$ is trivial; in other words we consider only isomorphisms between extensions of $(G^{\prime\ad}_{\dbQ},T^{\star}_{w\dbQ})$ and $(G^{\ad}_{\dbQ},T^{\star}_{w\dbQ})$ to some field $K$ of characteristic $0$ that have the property that over a larger field $K_1$ which contains both $K$ and $\dbR$, are defined by isomorphisms $G^{\prime\ad}_{K_1}\arrowsim G^{\ad}_{K_1}$ that are composites of the extension to $K_1$ of the adjoint of $j_{\dbR}$ with inner automorphisms of $G^{\ad}_{K_1}$. The class $\gamma_n$ is trivial over $\dbQ_p$. The class $\gamma_n$ is also trivial over $\dbR$ (as $\gamma_{0\dbR}$ is trivial and as two maximal compact tori of $G^{\ad}_{\dbR}$ are $G^{\ad}_{\dbR}(\dbR)$-conjugate). Let $T_{w\dbQ}^\prime$ be the torus of $G_{\dbQ}$ that is the twist of $T_{w\dbQ}$ via $\gamma_n$. The tori $T_{w\dbQ_p}^\prime$ and $T_{w\dbQ_p}$ are $G^{\ad}_{\dbQ}(\dbQ_p)$-conjugate. Let $x_w^\prime:\dbS\hookrightarrow T^\prime_{w\dbR}$ be the natural twist of $x_w$ via $\gamma_n$. We need an extra property:

\medskip
{\bf (iv)} Let $T_1$ and $T_2$ be two maximal tori of $G_{\dbZ_p}$ whose images in $G^{\ad}_{\dbZ_p}$ are $G^{\ad}_{\dbZ_p}(\dbZ_p)$-conjugate. Then $T_1$ and $T_2$ are $G_{\dbZ_p}(\dbZ_p)$-conjugate. 

\medskip
To check (iv), let $T_{1\dbF_p}^{\star}$ be the image of $T_{1\dbF_p}$ in $G^{\ad}_{\dbF_p}$. It is easy to check that we have an identity $G_{\dbZ_p}^{\ad}(\dbF_p)=\im(G_{\dbZ_p}(\dbF_p)\to G_{\dbZ_p}^{\ad}(\dbF_p))T^{\star}_{1\dbF_p}(\dbF_p)$. Thus the $G_{\dbZ_p}^{\ad}(\dbF_p)$-conjugates of $T_{1\dbF_p}$ are the same as the $G_{\dbZ_p}(\dbF_p)$-conjugates of $T_{1\dbF_p}$. Thus $T_{1\dbF_p}$ and $T_{2\dbF_p}$ are $G_{\dbZ_p}(\dbF_p)$-conjugate. From this and the infinitesimal liftings of [10, Vol. II, Exp. IX, Thm. 3.6], we get that for all $u\in\dbN$ there exists $\ell_u\in G_{\dbZ_p}(\dbZ_p)$ such that $\ell_uT_1\ell_u^{-1}$ and $T_2$ coincide mod $p^u$. Loc. cit. assures us that we can assume that $\ell_{u+1}$ and $\ell_u$ are congruent mod $p^u$. Thus the $p$-adic limit $\ell_{\infty}$ of the sequence $(\ell_u)_{u\in\dbN}$ exists and is an element of $G_{\dbZ_p}(\dbZ_p)$ with the property that $\ell_{\infty}T_1\ell_{\infty}^{-1}=T_2$. Thus (iv) holds.

We have an identity $G^{\ad}_{\dbQ}(\dbQ_p)=G^{\ad}_{\dbZ_p}(\dbZ_p)G^{\ad}_{\dbQ}(\dbQ)$, cf. [31, Lemma 4.9]. Thus, as the tori $T_{w\dbQ_p}^\prime$ and $T_{w\dbQ_p}$ are $G^{\ad}_{\dbQ}(\dbQ_p)$-conjugate, up to a $G^{\ad}_{\dbQ}(\dbQ)$-conjugation of $T_{w\dbQ}^\prime$ we can assume that the Zariski closure $T_w^{\prime}$ of $T_{w\dbQ}^\prime$ in $\GL_{L_{(p)}}$ is a torus with the property that $T_{w\dbZ_p}^{\prime}$ and $T_{w\dbZ_p}$ are $G^{\ad}_{\dbZ_p}(\dbZ_p)$-conjugate. Thus $T_{w\dbZ_p}^{\prime}$ and $T_{w\dbZ_p}$ are $G_{\dbZ_p}(\dbZ_p)$-conjugate, cf. (iv). Thus $T_{w\dbZ_p}^{\prime}$ and $T^w$ are $G_{\dbZ_p}(\dbZ_p)$-conjugate, cf. (ii). As $(G^{\prime}_{\dbQ},X_w^{\prime})$ is a Shimura pair, by replacing $(T_w,x_w)$ with $(T_w^{\prime},x_w^\prime)$ we can assume that $(G_{\dbQ},X_w)$ is a Shimura pair and that the properties (i) to (iii) continue to hold.  

\smallskip
{\bf Adjoint Shimura pairs.} The adjoint Shimura pair $(G_{\dbQ}^{\ad},X^{\ad})$ of $(G_{\dbQ},X)$ is defined by the property that $X^{\ad}$ is the $G^{\ad}_{\dbR}(\dbR)$-conjugacy class of the composite of any $x\in X$ with the epimorphism $G_{\dbR}\twoheadrightarrow G_{\dbR}^{\ad}$. We check that the adjoint Shimura pairs $(G_{\dbQ}^{\ad},X^{\ad})$ and $(G_{\dbQ}^{\ad},X_w^{\ad})$ coincide i.e., we have $X^{\ad}=X^{\ad}_w$. Let $G_{0\dbR}\in\grN$ be a simple, non-compact factor of $G^{\ad}_{\dbR}$. Let $X_0^{\ad}$ (resp. $X_{0w}^{\ad}$) be the direct factor of $X^{\ad}$ (resp. of $X^{\ad}_{w}$) that is the $G_{0\dbR}(\dbR)$-conjugacy class of the composite $x_{\star}$ (resp. $x_{0w}$) of any $x\in X$ (resp. of $x_w\in X_w$) with the epimorphism $G_{\dbR}\twoheadrightarrow G_{0\dbR}$. We know that $x_{\star\dbC}$ and $x_{0w\dbC}$ are $G_{0\dbR}(\dbC)$-conjugate, cf. (iii).  Thus $x_{\star}$ and $x_{0w}$ are $G_{0\dbR}(\dbR)$-conjugate, cf. [8, Prop. 1.2.2] and the fact that each element of $G_{0\dbR}(\dbC)$ that normalizes $G_{0\dbR}$ does belong to $G_{0\dbR}(\dbR)$. Thus we have $X_{0}^{\ad}=X^{\ad}_{0w}$. This implies that $X^{\ad}=X_w^{\ad}$. 

Both $X$ and $X_w$ are disjoint unions of connected components of $X^{\ad}$. The connected components of $X^{\ad}$ are permuted transitively by $G^{\ad}_{\dbZ_{(p)}}(\dbZ_{(p)})$, cf. [45, Cor. 3.3.3]. Thus by replacing the injective map $i_w:(T_{w\dbQ},\{x_w\})\hookrightarrow (G_{\dbQ},X_w)$ with its composite with an isomorphism $(G_{\dbQ},X_w)\arrowsim (G_{\dbQ},X)$ 
defined by an element of $G^{\ad}_{\dbZ_{(p)}}(\dbZ_{(p)})$, we can assume that $X_w=X$. The fact that under such a replacement (ii) still holds, is implied by (iv). 

\smallskip
{\bf A special Shimura pair.} As a conclusion, $(T_{w\dbQ},\{x_w\})$ is a special Shimura pair of $(G_{\dbQ},X)=(G_{\dbQ},X_w)$. The reflex field $E(T_{w\dbQ},\{x_w\})$ is a finite field extension of $E(G_{\dbQ},X)$. As $T_{wW(\dbF)}$ is a split torus, $E(T_{w\dbQ},\{x_w\})$ is naturally a subfield of $e_{\dbF}(W(\dbF))[{1\over p}]$. Let $v_w$ be the prime of $E(T_{w\dbQ},\{x_w\})$ such that the localization $O_{(v_w)}$ of the ring of integers of $E(T_{w\dbQ},\{x_w\})$ with respect to it is $E(T_{w\dbQ},\{x_w\})\cap e_{\dbF}(W(\dbF))$. The prime $v_w$ divides $v$ and the cocharacter $\mu_{0w}$ of ${T_w}_{W(k(v_w))}$ we get as in Subsubsection 5.1.1 but for the potential standard Hodge situation $(f\circ i_w,L,v_w)$, is such that its extension to $\dbC$ via $e_{\dbF}$ is $\mu_w$ itself. Thus the extension of $\mu_{0w}$ to $W(\dbF)$ is the cocharacter $g_w\mu^w_{0W(\dbF)}g_w^{-1}$ of $T_{wW(\dbF)}$. The normalization $\scrT_w$ of $\Spec(O_{(v_w)})$ in $\Sh(T_{w\dbQ},\{x_w\})/T_w(\dbZ_p)$ is a pro-\'etale cover of $\Spec(O_{(v_w)})$ (see [45, Ex. 3.2.8]) and thus its local rings are discrete valuation rings. From N\'eron--Ogg--Shafarevich criterion we get that the pullbacks of $(\scrA,\scrP_{\scrA})_{E(G_{\dbQ},X)}$ and of its symplectic similitude structures via the natural morphism $\Sh(T_{w\dbQ},\{x_w\})/T_w(\dbZ_p)\hookrightarrow \Sh(G_{\dbQ},X)_{E(T_{w\dbQ},\{x_w\})}/H$, extend to $\scrT_w$. Thus we have a natural morphism $\scrT_w\to\scrM_{O_{(v_w)}}$ which (by the very definitions of $\scrN$ and $\scrT_w$) factors through $\scrN_{O_{(v_w)}}$. Next we use the notations of the axiom 5.1.2 (b) and Notations 5.1.8 for a point $z\in\im(\scrT_w(W(k))\to\scrN(W(k)))$. The triple $(f\circ i_w,L,v_w)$ is a standard Hodge situation, cf. Example 5.1.4. Let $(v_{\alpha})_{\alpha\in\scrJ_w}$ be a family of tensors of $\scrT(L^*_{(p)})$ such that $\scrJ\subseteq\scrJ_w$ and $T_{w\dbQ}$ is the subgroup of $\GL_W$ that fixes $v_{\alpha}$ for all $\alpha\in\scrJ_w$. For $\alpha\in\scrJ_w\setminus\scrJ$ let $t_{\alpha}\in\scrT(M[{1\over p}])$ be the de Rham component of the Hodge cycle on $A_{B(k)}$ whose $p$ component of its \'etale component is $v_{\alpha}$ (thus $t_{\alpha}\in\scrT(M[{1\over p}])$ is well defined for all $\alpha\in\scrJ_w$). 

\smallskip
{\bf End of the proof of (b).} To end the proof of (b) we only have to show the existence of a rational isomorphism between $\scrC_y$ and $\scrC_w$. Let $\phi_{0w}:=(1_{L_p^*}\otimes \sigma_{k(v_w)})\circ\mu_{0w}({1\over p})$; it is a $\sigma_{k(v_w)}$-linear endomorphism of $M_{0w}:=L_p^*\otimes_{\dbZ_p} W(k(v_w))$. But $g_wh_w\in G_{\dbZ_p}(W(k))$ defines an isomorphism between $\scrC_w$ and $C_w:=(M_{0w},\phi_{0w},G_{W(k(v_w))},(v_{\alpha})_{\alpha\in\scrJ})\otimes_{k(v_w)} k$. Thus we are left to check that there exists a rational isomorphism $\scrI$ between $\scrC_y$ and $C_w$.

By applying Lemma 5.1.9 to the context of $z$ and $\scrJ_w$, we get the existence of an isomorphism $\rho:M[{1\over p}]\arrowsim M_{0w}\otimes_{W(k(v_w))} B(k)$ that takes $t_{\alpha}$ to $v_{\alpha}$ for all $\alpha\in\scrJ_w$.  Thus we also have $\rho\mu_{B(k)}\rho^{-1}=\mu_{0wB(k)}$ and moreover the element $g_y$ of Lemma 5.1.9 belongs to $T_w(B(k))$. This implies that $\Lie(\tilde T_w)$ is normalized by $\phi$ and $\mu$ thus (as in the proof of Claim 2.2.2 we argue that) we can speak about the maximal torus $\tilde T_{w\dbZ_p}$ of $\tilde G_{\dbZ_p}$ whose extension to $W(k)$ is $\tilde T_w$. Let $\gamma\in H^1(\dbQ_p,T_{w\dbQ_p})$ be the class that defines the torsor of $T_{w\dbQ_p}$ which parametrizes isomorphisms $M_{\dbZ_p}[{1\over p}]\arrowsim L_p^*[{1\over p}]$ that take $t_{\alpha}$ to $v_{\alpha}$ for all $\alpha\in\scrJ_w$. If the condition 5.2.1 (ii) holds, then $\gamma$ is the trivial class. If $\gamma$ is the trivial class, then we can assume $\rho$ is the tensorization with $B(k)$ of an isomorphism $M_{\dbZ_p}[{1\over p}]\arrowsim L_p^*[{1\over p}]$ and thus $g_y$ is the identity element; therefore $\scrI$ exists. Thus to end the proof of (b), we only have to show that $\gamma$ is the trivial class if the condition 5.2.1 (i) holds, $Z(G_{\dbQ})$ is connected, and the group $G_{\dbQ_p}(\dbQ_p)$ surjects onto $G^{\ad}_{\dbQ_p}(\dbQ_p)$.

Let $g_y^{\ad}$ be the image of $g_y$ in $T^{\star}_{w\dbQ_p}(B(k))$. The two hyperspecial subgroups of $G^{\ad}(B(k))$ that normalize $\Lie(G^{\ad}_{W(k)})$ and $\rho(\Lie(\tilde G^{\ad}_{W(k)})\rho^{-1}$ are inner conjugate by an element $g_0\in G^{\ad}_{B(k)}(B(k))$, cf. [43, p.  47]. As $T^{\star}_{wW(k)}(W(k))$ is a subgroup of both hyperspecial subgroups and as two maximal tori of $G^{\ad}_{W(k)}$ are $G^{\ad}_{W(k)}(W(k))$-conjugate, we can assume that $g_0\in T^{\star}_{wW(k)}(B(k))$. As $\rho\mu_{B(k)}\rho^{-1}=\mu_{0wB(k)}$, $g_0g_y^{\ad}\sigma(g_0)^{-1}$ normalizes $\Lie(G^{\ad}_{W(k)})$ and therefore it belongs to $G^{\ad}_{W(k)}(W(k))$ and thus also to $T^{\star}_{wW(k)}(W(k))$. Let $t^{\star}\in T^{\star}_{wW(k)}(W(k))$ be such that $g_0g_y^{\ad}\sigma(g_0)^{-1}=t^{\star}\sigma(t^{\star})^{-1}$, cf. Fact 2.6.2. As $Z(G_{\dbQ})=\Ker(T_{w\dbQ}\to T^{\star}_{w\dbQ})$ is connected, there exists an element $t\in T_w(B(k))$ that maps to $g_0^{-1}t^{\star}\in T^{\star}_w(B(k))$. Thus up to a replacement of $\rho$ by $t^{-1}\circ\rho$, we can assume that $g_y\in Z(G_{\dbQ})(B(k))$. This implies that the image $\gamma^{\star}$ of $\gamma$ in $H^1(\dbQ_p,T^{\star}_{w\dbQ_p})$ is the trivial class. As the condition 5.2.1 (i) holds, the image of $\gamma$ in $H^1(\dbQ_p,G_{\dbQ_p})$ is also the trivial class. Thus there exists $h\in G(B(k))$ such that $h\rho(M_{\dbZ_p}[{1\over p}])=L_p^*[{1\over p}]$. As $\gamma^{\star}$ is the trivial class, there exists $h_1^{\ad}\in G^{\ad}_{\dbQ_p}(\dbQ_p)$ such that $h_1^{\ad}$ and $h$ act (via inner conjugation) in the same way on $\Lie(T^{\star}_{w\dbQ_p})$. Let $h_1\in G_{\dbQ_p}(\dbQ_p)$ be such that it maps to $h_1^{\ad}$. By replacing $h$ with $h_1^{-1}h$, we can assume that $h$ fixes $\Lie(T^{\star}_{w\dbQ_p})$. Thus $h\in T_w(B(k))$ and therefore $\gamma$ is the trivial class. \endproof

\medskip\noindent
{\bf 5.2.4. Simple properties.} {\bf (a)} If $G_{\dbQ}=\GSp(W,\psi)$, then the condition 5.2.1 (i) holds (cf. Remark 5.2.2 (a)), $Z(G_{\dbQ})=\dbG_m$ is connected, and $G_{\dbQ_p}(\dbQ_p)$ surjects onto $G^{\ad}_{\dbQ_p}(\dbQ_p)$. Thus properties 5.2 (a) and (b) hold even if $p=2$, cf. Main Theorem 5.2.3. By combining this with Fact 1.1 we get a new, self-contained proof of the original Manin problem.

\smallskip
{\bf (b)} We refer to the proof of Theorem 5.2.3 (b). As there exists an isomorphism $M_{\dbZ_p}\arrowsim L_p^*=H^1$ that takes $t_{\alpha}$ to $v_{\alpha}$ for all $\alpha\in\scrJ_w$ (cf. Remark 5.2.2 (c)), $g_y\in T_w(B(k))$ is the identity and thus $\scrC_w$ and $\scrC_y$ are in fact isomorphic.

\medskip\smallskip\noindent
{\bf 5.3. Rational stratifications.} Let $\scrS_{NP}$ be the stratification of $\scrN_{k(v)}$ in locally closed, reduced subschemes defined by Newton polygons of pullbacks of the $p$-divisible group $\scrD$ of $\scrA$ via geometric points of $\scrN_{k(v)}$. It is the Newton polygon stratification of $\scrN_{k(v)}$ associated to the $F$-crystal over $\scrN_{k(v)}$ defined by $\scrD$, cf. [19, Thm. 2.3.1]. 

\medskip\noindent
{\bf 5.3.1. Theorem.} {\it Let $(f,L,v)$ be a standard Hodge situation. 

\medskip
{\bf (a)} There exists a stratification $\scrS_{\text{rat}}$ of $\scrN_{k(v)}$ in locally closed, reduced subschemes such that two points of $\scrN_{k(v)}$ with values in the same algebraically closed field factor through the same stratum of $\scrS_{\text{rat}}$ if and only if there exists a rational isomorphism between their attached Shimura $F$-crystals with tensors. 

\smallskip
{\bf (b)} Each stratum of $\scrS_{\text{rat}}$ is an open closed subscheme of a stratum of $\scrS_{NP}$ (thus $\scrS_{\text{rat}}$ refines $\scrS_{NP}$).}

\medskip
\proof We use left lower indices to denote pullbacks of $F$-crystals. Let $x$ be an independent variable. Let $S_0$ be a stratum of $\scrS_{NP}$. Let $S_1$ be an irreducible component of $S_0$. To prove the Theorem it is enough to show that for each two geometric points $y_1$ and $y_2$ of $S_1$ with values in the same algebraically closed field $k_1$, there exists a rational isomorphism between $\scrC_{y_1}$ and $\scrC_{y_2}$. We can assume that $k_1=\overline{k((x))}$ and that $y_1$ and $y_2$ factor through the generic point and respectively the special point of a $k$-morphism $m:\Spec(k[[x]])\to\scrN_k$. We denote also by $y_1$ and $y_2$, the $k_1$-valued points of either $\Spec(k[[x]])$ or its perfection $\Spec(k[[x]]^{\text{perf}})$ defined naturally by the factorizations of $y_1$ and $y_2$ through $m$. 

Let $\Phi$ be the Frobenius lift of $W(k)[[x]]$ that is compatible with $\sigma$ and that takes $x$ to $x^p$. Let $\grC=(V,\phi_V,\nabla_V)$ be the $F$-crystal over $k[[x]]$ of $m^*(\scrD)$. Thus $V$ is a free $W(k)[[x]]$-module, $\phi_V:V\to V$ is a $\Phi$-linear endomorphism, and $\nabla_V:V\to Vdx$ is a connection. Let $t_{\alpha}^V\in\scrT(V)$ be the de Rham realization of $m^*(w_{\alpha}^{\scrA})$ (see paragraph before Notations 5.1.8). Let $\scrC_{y_1}=(M_1,\phi_1,\tilde G_{1W(k_1)},(t_{1\alpha})_{\alpha\in\scrJ})$ (see Definition 5.1.10). Let $g_1\in\tilde G_{1W(k_1)}(W(k_1))$ be such that $(M_1,g_1\phi_1,\tilde G_{1W(k_1)},(t_{1\alpha})_{\alpha\in\scrJ})$ is the extension to $k_1$ of a Shimura $F$-crystal with tensors $\scrC_1$ over a finite field $k_{01}$ and there exists a rational isomorphism $l_1$ between $\scrC_{y_1}$ and $\scrC_1\otimes_{k_{01}} k_1=(M_1,g_1\phi_1,\tilde G_{1W(k_1)},(t_{1\alpha})_{\alpha\in\scrJ})$ defined by an element $h_1\in\tilde G_{1W(k_1)}(B(k_1))$, cf. Corollary 4.3. Let $\scrC_1^-$ be $\scrC_1$ but viewed only as an $F$-crystal. We can identify $\scrC_{1k_1}^-=(h_1^{-1}(M_1),\phi_1)$.

From [19, Thm. 2.7.4] we get the existence of an isogeny $i_0:\grC_0\to\grC$, where $\grC_0$ is an $F$-crystal over $k[[x]]$ whose extension to the $k[[x]]$-subalgebra $k[[x]]^{\text{perf}}$ of $k_1$ is constant (i.e., is the pullback of an $F$-crystal over $k$). Let $i_1:M_2\to M_1$ be the $W(k_1)$-linear map that defines $y_1^*(i_0)$. We can assume that $i_1(M_2)$ is contained in $h_1^{-1}(M_1)$. We get a morphism $\grC_{0k_1}\to\scrC^-_{1k_1}$. It is the extension to $k_1$ of a morphism $i_2:\grC_{0k[[x]]^{\text{perf}}}\to\scrC^-_{1k[[x]]^{\text{perf}}}$, cf. [37, Lemma 3.9] and the fact that $\grC_{0k[[x]]^{\text{perf}}}$ and $\scrC^-_{1k[[x]]^{\text{perf}}}$ are constant $F$-crystals over $k[[x]]^{\text{perf}}$. Let $i_3:\scrC^-_{1k[[x]]^{\text{perf}}}\to\grC_{0k[[x]]^{\text{perf}}}$ be a morphism such that $i_2\circ i_3=p^q1_{\scrC^-_{1k[[x]]^{\text{perf}}}}$, where $q\in\dbN$. By composing $i_3$ with $i_{0k[[x]]^{\text{perf}}}$ we get an isogeny 
$i_4:\scrC^-_{1k[[x]]^{\text{perf}}}\to\grC_{k[[x]]^{\text{perf}}}$
whose extension to $k_1$ is defined by the inclusion $h_1^{-1}(p^qM_1)\subseteq M_1$. The isomorphism of $F$-isocrystals over $\Spec(k[[x]]^{\text{perf}})$ defined by $p^{-q}$ times $i_4$ takes $t_{1\alpha}$ to $t_{\alpha}^V$ for all $\alpha\in\scrJ$, as this is so generically. Thus $y_2^*(i_4)$ is an isogeny which when viewed as an isomorphism of $F$-isocrystals is $p^q$ times a rational isomorphism $l_2$ between $\scrC_1\otimes_{k_{01}} k_1$ and $\scrC_{y_2}$. Therefore $l_2l_1$ is a rational isomorphism between $\scrC_{y_1}$ and $\scrC_{y_2}$. \endproof

\medskip\noindent
{\bf 5.3.2. Remarks.}  {\bf (a)} The stratification $\scrS_{\text{rat}}$ was introduced in [37]: it is only a concrete example of the stratifications in characteristic $p$ one gets based on [37, Thm. 3.6 (ii)].  

\smallskip
{\bf (b)} Theorem 5.3.1 (b) is a slight refinement of a concrete example of [37, Thm. 3.8] as it weakens the hypotheses of loc. cit. More precisely, it considers the ``Newton point'' of only one faithful
representation (which in the case when condition 5.2.1 (i) holds, is the representation of $G_{\dbQ_p}$ on $W^∗\otimes_{\dbQ}\dbQ_p$) while in loc. cit. one has to consider the ``Newton points'' of all finite dimensional representations of $G_{\dbQ_p}$. 

\smallskip
{\bf (c)} The stratifications $\scrS_{\text{rat}}$ and $\scrS_{NP}$ are $G_{\dbQ}(\dbA_f^{(p)})$-invariant (cf. Remark 5.1.11) and therefore are pullbacks of stratifications $\scrS_{\text{rat},H_0}$ and $\scrS_{NP,H_0}$ of the $k(v)$-scheme $\scrN_{k(v)}/H_0$ of finite type, where $H_0$ is as before Notations 5.1.8. Thus $\scrS_{NP,H_0}$ has a finite number of strata and each stratum of it has a finite number of connected components. Therefore $\scrS_{\text{rat},H_0}$ and thus also $\scrS_{\text{rat}}$ has a finite number of strata, cf. Theorem 5.3.1 (b). If $\grR$ is as in Corollary 4.4 and if properties 5.2 (a) and (b) hold, then $\scrS_{\text{rat}}$ has precisely $\grR$ strata.

\medskip\noindent
{\bf 5.3.3. Example.} We assume that $G^{\ad}_{\dbQ}$ is absolutely simple of $B_n$ Dynkin type, that $Z(G_{\dbQ})=\dbG_m$, and that $(f,L,v)$ is a standard Hodge situation. The faithful representation of $G^{\der}_{\dbC}$ on $W\otimes_{\dbQ} \dbC$ is a direct sum of trivial and spin representations (see [8]) and thus $G^{\der}_{\dbQ}$ is simply connected. Thus the condition 5.2.1 (i) holds, cf. Remark 5.2.2 (a). As $Z(G_{\dbQ})=\dbG_m$ is connected and as $G_{\dbQ_p}(\dbQ_p)$ surjects onto $G^{\ab}_{\dbQ_p}(\dbQ_p)$, the hypotheses of the Main Theorem 5.2.3 hold. Thus properties 5.2 (a) and (b) hold even if $p=2$, cf. Main Theorem 5.2.3. From this and Proposition 4.6.1 (a) we get that $\scrS_{\text{rat}}$ has $n+1$ strata (cf. also end of Remark 5.3.2 (c)). 

\medskip\smallskip\noindent
{\bf 5.4. Integral Manin problem for $(f,L,v)$ and $k$.} Show that the following two things hold:

\medskip
{\bf (a)} For each point $y\in\scrN(k)$, there exists an element $g\in G^0(W(k))$ such that $\scrC_g$ is isomorphic to $\scrC_y$.

\smallskip
{\bf (b)} For each element $g\in G^0(W(k))$, there exists a point $\tilde y\in\scrN(k)$ such that $\scrC_{\tilde y}$ is isomorphic to $\scrC_g$.

\medskip\noindent
{\bf 5.4.1. Some standard PEL situations.} Let $\iota$ be the involution of $\End(L_{(p)})$ defined by the identity $\psi(b(x),y)=\psi(x,\iota(b)(y))$,
where $b\in\End(L_{(p)})$ and $x$, $y\in L_{(p)}$. Let $\scrB:=\{b\in\End(L_{(p)})|b\,{\text{is}}\,{\text{fixed}}\,{\text{by}}\, G_{\dbZ_{(p)}}\}$.
We list four conditions:
\medskip
{\bf (i)} we have $\iota(\scrB)=\scrB$ and $\scrB[{1\over p}]$ is a simple $\dbQ$--algebra;
\smallskip
{\bf (ii)} the $W(\dbF)$-algebra $\scrB\otimes_{\dbZ_{(p)}} W(\dbF)$ is a product of matrix $W(\dbF)$-algebras;
\smallskip
{\bf (iii)} the group $G_{\dbQ}$ is the subgroup of $\GSp(W,\psi)$ that fixes all $b\in\scrB$;
\smallskip
{\bf (iv)} the Hasse principle holds for the reductive group $G_{\dbQ}$. 

\medskip\noindent
If (i) to (iii) hold, then the triple $(f,L,v)$ is a standard Hodge situation (see [28] and [24]); we refer to $(f,L,v)$ as a standard PEL situation. Condition (iv) holds if (i) to (iii) hold and all simple factors of $G^{\ad}_{\dbC}$ are of $C_n$ or $A_{2n-1}$ Lie type $(n\in\dbN$), cf. [24, pp. 393--394]. 

\medskip\noindent
{\bf 5.4.2. Theorem.} {\it If conditions 5.4.1 (i) to (iv) hold, then 5.4 (a) and (b) also hold.}

\medskip
\proof 
We use the notations (such as $y$ and $z$) of the axiom 5.1.2 (b), Notations 5.1.8, and Lemma 5.1.9. The isomorphism $i_A\otimes_{B^+(W(k))} 1_{B_{\text{dR}}(W(k))}$ takes $\psi_M$ to a non-zero scalar multiple of $\psi^*$. As $\tilde G$ does not fix $\psi_M$, we get the existence of a field extension $L_{\text{dR}}(W(k))$ of the field of fractions of $B_{\text{dR}}(W(k))$ such that there exists a symplectic isomorphism $(M\otimes_{W(k)} L_{\text{dR}}(W(k)),\psi_M)\arrowsim (L_p^*\otimes_{\dbZ_p} L_{\text{dR}}(W(k)),\psi^*)$ that takes $t_{\alpha}$ to $v_{\alpha}$ for all $\alpha\in\scrJ$. As $k=\bar k$ and as $G^0_{\dbQ}=\Sp(W,\psi)\cap G_{\dbQ}$ is connected (cf. (12)), the set $H^1(B(k),G^0_{B(k)})$ has only one class. Thus referring to Lemma 5.1.9, we can choose $\rho$ such that it defines a symplectic isomorphism $(M[{1\over p}],\psi_M)\arrowsim (L_p^*\otimes_{\dbZ_p} B(k),\psi^*)$. 

Let $\scrJ_{\scrB}:=\{\alpha\in\scrJ|v_{\alpha}\in\scrB\}$. We identify $\scrB\otimes_{\dbZ_{(p)}} \dbZ_p$ with the $\dbZ_p$-span of $v_{\alpha}$'s (resp. $t_{\alpha}$'s) with $\alpha\in\scrJ_{\scrB}$. Let $\rho_1:M\arrowsim M_0\otimes_{W(k(v))} W(k)$ 
be an isomorphism that takes $\psi_M$ to $\psi^*$ and $t_{\alpha}$ to $v_{\alpha}$, for all $\alpha\in\scrJ_{\scrB}$. The existence of $\rho_1$ after inverting $p$ follows from the previous paragraph and thus [24, Lemma 7.2] implies that $\rho_1$ exists. Strictly speaking, loc. cit. is stated over $\dbZ_p$ but its arguments apply entirely over $W(k)$. As $\rho\circ\rho_1^{-1}$ fixes $\psi^*$ and each $v_{\alpha}$ with $\alpha\in\scrJ_{\scrB}$, we have $\rho\circ\rho_1^{-1}\in G^0_{\dbQ}(B(k))$ (cf. condition 5.4.1 (iii)). Thus $\rho_1$ also takes $t_{\alpha}$ to $v_{\alpha}$ for all $\alpha\in\scrJ\setminus\scrJ_{\scrB}$. This implies that we can choose $\rho$ such that moreover we have $\rho(M)=M_0\otimes_{W(k(v))} W(k)$ and $\rho\mu\rho^{-1}=\mu_{0W(k)}$. Thus the representations of $\scrB/p\scrB$ on $\Lie(A_k)=F^1/pF^1$ and $F^1_0/pF^1_0\otimes_{k(v)} k$ are isomorphic and moreover we have $g_y\in G^0_{W(k)}(W(k))$. Thus the property 5.4 (a) holds.

We show that the property 5.4 (b) holds. Let $g\in G^0_{W(k)}(W(k))$. Let $N^0$ be the normalizer of $T^0$ in $G_{\dbZ_p}$. Let $w\in N^0(W(k))\cap G^0_{W(k)}(W(k))$ be such that there exists an element $h\in G^0_{W(k)}(B(k))$ that defines a rational isomorphism between $\scrC_w$ and $\scrC_g$, cf. Subsection 4.7. We take $z\in\scrN(W(k))$ to factor through $\scrT_w$ of the proof of Theorem 5.2.3 (b) and we use the notations of the mentioned proof. Let $T^0_w$ be the subgroup scheme of $T_w$ that fixes $\psi$. We have $T^0_w=G^0_{\dbZ_{(p)}}\cap T_w$ (cf. (12)) and therefore $T^0_w$ is a maximal torus of $G^0_{\dbZ_{(p)}}$. If $\scrB_w$ is the centralizer of $T_w$ in $\End(L_{(p)})$, then $T_w$ is the subgroup scheme of $\GSp(L_{(p)},\psi)$ that fixes $\scrB_w$ and $\scrB_w\otimes_{\dbZ_{(p)}} W(\dbF)$ is a product of matrix $W(\dbF)$-algebras. Thus as in the previous paragraph we argue that can assume that $g_y\in T^0_w(W(k))$. By composing $\rho$ with an element of $T^0_w(W(k))$, we can assume that $g_y$ is the identity element. Therefore $\scrC_y$ is isomorphic to $C_w$ (of the proof of Theorem 5.2.3 (b)). Thus $\scrC_y$ is isomorphic to $\scrC_w$ under an isomorphism $\rho:M\arrowsim M_0\otimes_{W(k(v))} W(k)$ that takes $\psi_M$ to $\psi^*$ and $t_{\alpha}$ to $v_{\alpha}$ for all $\alpha\in\scrJ_w$ and therefore also for all $\alpha\in\scrJ$.

Let $M_1:=(h\rho)^{-1}(L_p^*\otimes_{\dbZ_{p}} W(k))$. It is a $W(k)$-lattice of $M[{1\over p}]$. Let $\tilde G_1$ be the Zariski closure of $G_{B(k)}$ in $\GL_{M_1}$. The quadruple $(M_1,\phi,\tilde G_1,(t_{\alpha})_{\alpha\in\scrJ})$
is a Shimura $F$-crystal with tensors isomorphic to $\scrC_g$ (via $h\rho$). Moreover, $M_1$ is self dual with respect to $\psi_M$ and we have $t_{\alpha}(M_1)\subseteq M_1$ for all $\alpha\in\scrJ_{\scrB}$. Let $A_{1k}$ be the abelian variety over $k$ that is $\dbZ[{1\over p}]$-isogenous to $A_k$ and whose Dieudonn\'e module is (under this $\dbZ[{1\over p}]$-isogeny) $(M_1,\phi)$. Let $\lambda_{A_{1k}}$ be the principal polarization of $A_{1k}$ defined naturally by $\psi_M$ and the principal polarization $z^*(\scrP_{\scrA})$ of $A$. We endow $A_{1k}$ with the level $s$ symplectic similitude structures induced naturally by those of $A_k$, for all $s\in (\dbN\setminus p\dbN)$. To these structures and to $(A_{1k},\lambda_{A_{1k}})$ corresponds naturally a morphism $\tilde y:\Spec(k)\to\scrM$. Moreover, for $\alpha\in\scrJ_B$ the tensor $t_{\alpha}$ is the crystalline realization of a $\dbZ_{(p)}$-endomorphism of $A_{1k}$. 
To end the proof of property 5.4 (b) we only have to check that there exists a point $\tilde z\in\scrM(W(k))$ which lifts $\tilde y$ and which factors through $\scrN$ in such a way that $t_{\alpha}$ is the de Rham realization of the Hodge cycle $\tilde z^*(w_{\alpha}^{\scrA})$, for all $\alpha\in\scrJ$. Due to condition 5.4.1 (iii) it is enough to work in the last sentence only with indices $\alpha\in\scrJ_{\scrB}$. But due to condition 5.4.1 (iv), the existence of $\tilde z$ follows from the above part referring to representations of $\scrB/p\scrB$ and from the well known moduli considerations of [24, pp. 390 and 399]. \endproof 

\medskip\noindent
{\bf 5.4.3. Examples.} {\bf (a)} The principally quasi-polarized Dieudonn\'e module of a principally quasi-polarized $p$-divisible group over $k$ of height $2r$ is isomorphic to $(M,g\phi,\psi_M)$, where $g\in \Sp(M,\psi_M)(W(k))$. Thus Theorem 5.4.2 for Siegel modular varieties (i.e., for when $\scrN=\scrM$) says that each principally quasi-polarized $p$-divisible group over $k$ of height $2r$ is the one of a principally polarized abelian variety over $k$ of dimension $r$. This result was first obtained in [46, Prop. 5.3.3].

\smallskip
{\bf (b)} Suppose that the conditions 5.4 (i) to (iv) hold, that the $\dbQ$--algebra $\scrB[{1\over p}]$ is simple, that the group scheme $G_{\dbZ_p}$ is split, and that $G^{\der}_{\dbZ_p}$ is an $\SL_n$ group scheme. As $G_{\dbZ_p}$ is split we have $k(v)=\dbF_p$. As $G^{\ad}_{\dbC}$ is simple, the condition 5.4.1 (iv)  holds even if $n$ is odd (see [24, top of p. 394]). We also assume that we have a direct sum decomposition 
$$M_0=L_{p}^*=L_0\oplus L_1\leqno (15)$$ 
of $G_{\dbZ_p}$-modules such that the representation of $G^{\der}_{\dbZ_p}$ on $L_0$ is the standard faithful representation of rank $n$. Let $a\in S(1,n-1)$ be such that $\mu_0$ acts trivially (resp. non-trivially) on a direct summand of $L_0$ of rank $a$ (resp. rank $n-a$). It is well known that the derived group $G^{\der}_{\dbR}$ is isomorphic to $\pmb{SU}(a,n-a)_{\dbR}$. Let $\scrD$ be as in the beginning of Subsection 5.3. For $y\in\scrN(W(k))$, let $(M,\phi)$ be as in the axiom 5.1.2 (b). Let $\scrD_{W(k)}=\scrD_0\oplus\scrD_1$ and $M=N_0\oplus N_1$ be the direct sum decompositions that corresponds naturally to (15). The pair $(N_0,\phi)$ is a Dieudonn\'e module of height $n$ and dimension $n-a$. The Dieudonn\'e module of another $p$-divisible group over $k$ of height $n$ and dimension $n-a$ is isomorphic to $(N_0,g\phi)$ for some $g\in \SL_{N_0}(W(k))$, cf. Fact 2.6.3 applied to $(N_0,\phi,\GL_{N_0})$. From this and Theorem 5.4.2 we get that each $p$-divisible group over $k$ of height $n$ and dimension $n-a$ is isomorphic to $y^*(\scrD_0)$ for some point $y\in\scrN_{W(k)}(k)$.

\medskip\noindent
{\bf Acknowledgments.} We would like to thank University of Arizona and Binghamton University for good working conditions and the referee for several comments.

\medskip\smallskip
\noindent
\references{37}
{\nspace{

\medskip
\Ref[1] 
P. Berthelot, L. Breen, and W. Messing, 
\sl Th\'eorie de Dieudonn\'e cristalline II, 
\rm Lecture Notes in Math., Vol. {\bf 930}, Springer-Verlag, Berlin, 1982.

\Ref[2] D. Blasius, 
\sl A p-adic property of Hodge cycles on abelian varieties, 
\rm  Motives (Seattle, WA, 1991),  293--308, Proc. Sympos. Pure Math., {\bf 55}, Part 2, Amer. Math. Soc., Providence, RI, 1994.

\Ref[3]
A. Borel,
\sl Linear algebraic groups, 
\rm Grad. Texts in Math., Vol. {\bf 126}, Springer-Verlag, New York, 1991.

\Ref[4]
S. Bosch, W. L\"utkebohmert, and M. Raynaud,
\sl N\'eron models,
\rm Ergebnisse der Mathematik und ihrer Grenzgebiete (3). [Results in Mathematics and Related Areas (3)], Vol. {\bf 21}, Springer-Verlag, Berlin, 1990.

\Ref[5]
N. Bourbaki,
\sl Lie groups and Lie algebras, Chapters 4--6, 
\rm Springer-Verlag, Berlin, 2002.

\Ref[6] 
C.-L. Chai, 
\sl Newton polygons as lattice points, 
\rm Amer. J. Math. {\bf 122} (2000), no. 5,  967--990.

\Ref[7]
P. Deligne,
\sl Travaux de Shimura,
\rm S\'eminaire  Bourbaki, Exp. 389, Lecture Notes in Math., Vol. {\bf 244},  123--163, Springer-Verlag, Berlin, 1971.

\Ref[8]
P. Deligne,
\sl Vari\'et\'es de Shimura: interpr\'etation modulaire, et
techniques de construction de mod\`eles canoniques,
\rm (French) Automorphic forms, representations and $L$-functions (Oregon State Univ., Corvallis, OR, 1977), Part 2,   247--289, Proc. Sympos. Pure Math., {\bf 33}, Amer. Math. Soc., Providence, RI, 1979.

\Ref[9]
P. Deligne,
\sl Hodge cycles on abelian varieties,
\rm Hodge cycles, motives, and Shimura varieties, Lecture Notes in Math., Vol. {\bf 900},  9--100, Springer-Verlag, Berlin-New York, 1982. 

\Ref[10]
M. Demazure, A. Grothendieck, et al. 
\sl Sch\'emas en groupes, 
\rm Vols. I--III, Lecture Notes in Math., Vols. {\bf 151--153}, Springer-Verlag, Berlin, 1970.

\Ref[11] 
J. Dieudonn\'e, 
\sl Groupes de Lie et hyperalg\`ebres de Lie sur un corps de caract\'erisque $p>0$ (VII), 
\rm Math. Annalen {\bf 134} (1957),  114--133. 

\Ref[12] 
G. Faltings, 
\sl Crystalline cohomology and $p$-adic Galois representations,
\rm Algebraic analysis, geometry, and number theory (Baltimore, MD, 1988),  25--80, Johns Hopkins Univ. Press, Baltimore, MD, 1989.

\Ref[13] 
G. Faltings, 
\sl Integral crystalline cohomology over very ramified valuation rings, 
\rm J. Amer. Math. Soc. {\bf 12} (1999), no. 1,  117--144.

\Ref[14] 
J.-M. Fontaine and G. Laffaille, 
\sl Construction de repr\'esentations $p$-adiques, 
\rm Ann. Sci. \'Ecole Norm. Sup. {\bf 15} (1982), no. 4,  547--608.

\Ref[15] 
J.-M. Fontaine, 
\sl Le corps des p\'eriodes $p$-adiques, 
\rm J. Ast\'erisque {\bf 223},  59--101, Soc. Math. de France, Paris, 1994.

\Ref[16]
J.-M. Fontaine and M. Rapoport,
\sl Existence de filtrations admissibles sur des isocristaux,
\rm Bull. Soc. Math. de France {\bf 133} (2005),  73--86.

\Ref[17] 
G. Harder, 
\sl \"Uber die Galoiskohomologie halbeinfacher Matrizengruppen II, 
\rm Math. Z. {\bf 92} (1966),  396--415.

\Ref[18]
M. Harris and R. Taylor,
\sl The geometry and cohomology of some simple Shimura varieties,
\rm Annals of Mathematics Studies, Vol. {\bf 151}, Princeton Univ. Press, Princeton, NJ, 2001. 

\Ref[19] N. Katz, 
\sl Slope filtration of $F$-crystals, 
\rm Journ\'ees de G\'eom\'etrie alg\'ebrique de Rennes 1978, J. Ast\'erisque {\bf 63}, Part 1,  113--163, Soc. Math. de France, Paris, 1979.

\Ref[20]
M. Kisin,
\sl Integral canonical models of Shimura varieties of abelian type, 
\rm J. Amer. Math. Soc. {\bf 23} (2010), no. 4,  967--1012.

\Ref[21]
M. Knesser,
\sl Galois-Kohomologie halbeinfacher algebraischer Gruppen \"uber $p$-adischen K\"orpern I., 
\rm Math. Zeit. {\bf 88} (1965),  44--47; II., ibid. {\bf 89} (1965),  250--272. 

\Ref[22]
R. E. Kottwitz, 
\sl Stable trace formula: cuspidal tempered terms,
\rm Duke Math. J. {\bf 51} (1984), no. 3,  611--650.

\Ref[23]
R. E. Kottwitz, 
\sl Isocrystals with additional structure,
\rm Compositio Math. {\bf 56} (1985), no. 2,  201--220.

\Ref[24]
R. E. Kottwitz, 
\sl Points on some Shimura varieties over finite fields, 
\rm J. Amer. Math. Soc. {\bf 5} (1992), no. 2,  373--444.

\Ref[25] 
R. E. Kottwitz,
\sl On the Hodge--Newton decomposition for split groups,
\rm Int. Math. Res. Not. {\bf 26} (2003),  1433--1447. 

\Ref[26] 
G. Laffaille, 
\sl Groupes $p$-divisibles et modules filtr\'es: le cas peu ramifi\'e, 
\rm Bull. Soc. Math. de France {\bf 108} (1980),  187--206.

\Ref[27]
R. Langlands,
\sl Some contemporary problems with origin in the Jugendtraum,
\rm Mathematical developments arising from Hilbert's problems,  401--418, Amer. Math. Soc., Providence, RI, 1976.

\Ref[28]
R. Langlands and M. Rapoport,
\sl Shimuravariet\"aten und Gerben, 
\rm J. reine angew. Math. {\bf 378} (1987),  113--220.

\Ref[29] 
J. I. Manin, 
\sl The theory of formal commutative groups in finite characteristic, 
\rm Russian Math. Surv. {\bf 18} (1963), no. 6,  1--83. 

\Ref[30]
J. S. Milne,
\sl The points on a Shimura variety modulo a prime of good
reduction,
\rm The Zeta functions of Picard modular surfaces,  153--255, Univ. Montr\'eal Press, Montreal, Quebec, 1992.

\Ref[31]
J. S. Milne,
\sl Shimura varieties and motives,
\rm Motives (Seattle, WA, 1991),  447--523, Proc. Sympos. Pure Math., {\bf 55}, Part 2, Amer. Math. Soc., Providence, RI, 1994.

\Ref[32]
J. S. Milne and K. -y. Shih,
\sl Conjugates of Shimura varieties,
\rm  Hodge cycles, motives, and Shimura varieties, Lecture Notes in Math., Vol. {\bf 900},  280--356, Springer-Verlag, Berlin-New York, 1982.

\Ref[33]
M.-H. Nicole and A.~Vasiu,
\sl Minimal truncations of supersingular $p$-divisible groups,
\rm Indiana Univ. Math. J. {\bf 56} (2007), no. 6, 2887--2898.

\Ref[34] 
F. Oort, 
\sl Newton polygons and formal groups: conjectures by Manin and Grothendieck, 
\rm Ann. of Math. (2) {\bf 152} (2000), no. 1,  183--206.

\Ref[35] 
F. Oort, 
\sl Newton polygon strata in the moduli of abelian varieties, 
\rm Moduli of abelian varieties (Texel Island, 1999),  417--440, Progr. Math., {\bf 195}, Birkh\"auser, Basel, 2001.

\Ref[36]
R. Pink,
\sl $l$-adic algebraic monodromy groups, cocharacters, and the Mumford--Tate conjecture,
\rm J. reine angew. Math. {\bf 495} (1998),  187--237.

\Ref[37] 
M. Rapoport and M. Richartz, 
\sl On the classification and specialization of $F$-isocrystals with additional structure, 
\rm Compositio Math. {\bf 103} (1996), no. 2,  153--181.

\Ref[38]
M. Rapoport and T. Zink,
\sl Period spaces for $p$-divisible groups, 
\rm Annals of Mathematics Studies, Vol. {\bf 141}, Princeton Univ. Press, Princeton, NJ, 1996.

\Ref[39] 
J.-P. Serre, 
\sl Galois Cohomology, 
\rm Springer-Verlag, Berlin, 1997.

\Ref[40] 
G. Shimura, 
\sl Moduli of abelian varieties and number theory, 
\rm 1966 Algebraic Groups and Discontinuous Subgroups (Boulder, CO, 1965), Proc. Sympos. Pure Math., {\bf 9},  312--332, Amer. Math. Soc., Providence, RI.

\Ref[41] 
J. Tate, 
\sl Classes d'isog\'enie des vari\'et\'es sur un corps fini (d'apr\`es J. Honda), 
\rm S\'eminaire Bourbaki 1968/69, Exp. 352, Lecture Notes in Math., Vol. {\bf 179},  95--110, 1971. 

\Ref[42] 
J. Tits, 
\sl Classification of algebraic semisimple groups, 
\rm 1966 Algebraic Groups and Discontinuous Subgroups (Boulder, CO, 1965), Proc. Sympos. Pure Math., {\bf 9},  33--62, Amer. Math. Soc., Providence, RI.

\Ref[43]
J. Tits,
\sl Reductive groups over local fields, 
\rm Automorphic forms, representations and $L$-functions (Oregon State Univ., Corvallis, OR, 1977), Part 1,  29--69, Proc. Sympos. Pure Math., {\bf 33}, Amer. Math. Soc., Providence, RI, 1979.

\Ref[44]
A. Vasiu,
\sl Integral canonical models for Shimura varieties of Hodge type,
\rm Ph.D. Thesis, Princeton University, 1994.

\Ref[45]
A. Vasiu, 
\sl Integral canonical models of Shimura varieties of preabelian
 type, 
\rm Asian J. Math. {\bf 3} (1999), no. 2,  401--518.

\Ref[46]
A. Vasiu, 
\sl Crystalline boundedness principle,
\rm Ann. Sci. \'Ecole Norm. Sup. {\bf 39} (2006), no. 2,  245--300.

\Ref[47]
A. Vasiu,
\sl Level $m$ stratifications of versal deformations of $p$-divisible groups,
\rm J. Alg. Geom. {\bf 17} (2008), no. 4,  599--641.

\Ref[48]
A. Vasiu, 
\sl Mod $p$ classification of shimura $F$-crystals,
\rm Math. Nachr. {\bf 283} (2010), no. 8, 1068--1113.

\Ref[49]
A. Vasiu, 
\sl A motivic conjecture of Milne,
\rm math.NT/0308202. 

\Ref[50]
A. Vasiu, 
\sl Good reductions of Shimura varieties of Hodge type in arbitrary unramified mixed characteristic, Parts I and II,
\rm available at http://arxiv.org/abs/0707.1668 and http://arxiv.org/pdf/0712.1572

\Ref[51]
T. Zink,
\sl Isogenieklassen von Punkten von Shimuramannigfaltigkeiten mit Werten in einem endlichen K\"orper,
\rm Math. Nachr. {\bf 112} (1983),  103--124.

\Ref[52] 
T. Zink, 
\sl On the slope filtration, 
\rm Duke Math. J. {\bf 109} (2001), no. 1,  79--95.

\Ref[53]
W. Waterhouse, 
\sl Introduction to affine group schemes, 
\rm Grad. Texts in Math., Vol. {\bf 66}, Springer-Verlag, New York-Berlin, 1979.

\Ref[54] J.-P. Wintenberger, 
\sl Un scindage de la filtration de Hodge pour certaines
vari\'et\'es alg\'ebriques sur les corps locaux, 
\rm Ann. of Math. (2) {\bf 119} (1984), no. 3,  511--548.

\Ref[55] J.-P. Wintenberger,
\sl Existence de $F$-cristaux avec structures suppl\'ementaires, 
\rm Adv. Math. {\bf 190}  (2005),  no. 1,  196--224.

}}

\medskip
\hbox{Adrian Vasiu}
\hbox{Department of Mathematical Sciences,}
\hbox{Binghamton University, Binghamton, P. O. Box 6000,},
\hbox{New York 13902-6000, U.S.A.}
\hbox{e-mail: adrian\@math.binghamton.edu} 

\enddocument